\documentclass[11pt]{article}

\usepackage{amssymb,amsmath,amsthm,amsfonts} 
\usepackage{amscd}
\usepackage{bm} 
\input xy
\xyoption{all}

\usepackage[twoside]{geometry}
\newtheorem{theorem}{Theorem}[section]

\newtheorem{proposition}[theorem]{Proposition}
\newtheorem{lemma}[theorem]{Lemma}
\newtheorem{corollary}[theorem]{Corollary}

\theoremstyle{definition}

\newtheorem{remark}[theorem]{Remark}
\newtheorem{definition}[theorem]{Definition}

\newcommand{\X}{{\cal X}}
\newcommand{\cH}{{\cal H}}
\newcommand{\R}{\mathbb{R}}
\newcommand{\ra}{\rightarrow}
\newcommand{\lra}{\longrightarrow}
\newcommand{\sT}{\mathrm{T}}
\newcommand{\T}{\mathrm{T}}
\newcommand{\pa}{\partial}
\newcommand{\p}{\operatorname{p}}
\newcommand{\cl}{\operatorname{cl}}

\newcommand{\Image}{\operatorname{Im}}
\newcommand{\eps}{\varepsilon}
\newcommand{\dd}{\mathrm{d}}
\newcommand{\del}{\delta}

\newcommand{\B}{\bm}
\newcommand{\m}{\mathfrak{m}}
\newcommand{\BB}{\mathfrak{v}}
\newcommand{\xra}{\xrightarrow} 
\newcommand{\g}{\mathfrak{g}}
\newcommand{\so}{\mathfrak{so}(3)}
\newcommand{\SO}{\mathrm{SO}(3)}
\newcommand{\A}{\mathfrak{w}} 

\newcommand{\Sec}{\operatorname{Sec}}
\newcommand{\id}{\operatorname{id}}

\newcommand{\wt}{\widetilde}
\newcommand{\wh}{\widehat}

\newcommand{\Uadm}{\mathcal{U}_{adm}}
\newcommand{\BM}{\mathcal{BM}}
\newcommand{\GG}{\mathcal{G}}
\newcommand{\ol}{\overline}
\newcommand{\ul}{\underline}
\newcommand{\PP}{\operatorname{P}}
\newcommand{\AG}{\mathcal{A}(\mathcal{G})}
\def\<#1>{\big\langle #1\big\rangle}

\title{Pontryagin Maximum Principle -- a generalization\thanks{Research partially supported by the Polish
Ministry of Science and Higher Education under the grant N N201
416839.}}

\author{Janusz Grabowski\thanks{Institute of Mathematics, Polish Academy of Sciences, \'Sniadeckich 8, 00-956 Warszawa, Poland ({\tt jagrab@impan.pl}).}
\and Micha\l \ J\'o\'zwikowski\thanks{Institute of Mathematics, Polish Academy of Sciences, \'Sniadeckich 8, 00-956 Warszawa, Poland ({\tt mjozwikowski@gmail.com}).}}
\begin{document}
\maketitle

\begin{abstract}
The fundamental theorem of the theory of optimal control, the
Pontryagin maximum principle (PMP), is extended to the setting of
almost Lie (AL) algebroids, geometrical objects generalizing Lie
algebroids. This formulation of the PMP yields, in particular, a
scheme comprising reductions of optimal control problems similar
to the reduction for the rigid body in analytical mechanics. On
the other hand, in the presented approach the reduced and
unreduced PMPs are parts of the same universal formalism. The
framework is based on a very general concept of homotopy of
measurable paths and the geometry of AL algebroids.
\end{abstract}

{\small\textbf{Key words:}
optimal control, Lie algebroids, reduction of
control systems, variation, homotopy\\
}

{\small\textbf{MSC 2000:}
49K15, 49K99, 93C25 (Primary); 17B66, 22A22
(Secondary)}

\pagestyle{myheadings}
\thispagestyle{plain}
\markboth{J.~GRABOWSKI AND M.~J\'O\'ZWIKOWSKI}{PONTRYAGIN MAXIMUM PRINCIPLE}

\section{Introduction}\label{sec:intro}
\subsection{Reduction in optimal control theory}
It is a well-known phenomenon in analytical mechanics and control
theory that symmetries of a system lead to reductions of its
degrees of freedom. It is also well understood that such a
reduction procedure is not purely computational but is associated
with a reduction of the geometrical structures hidden behind.

A typical situation considered in control theory is a control
system $F:P\times U\lra\sT P$ on a manifold $P$ (with $U$ being
the set of control parameters) which is equivariant w.r.t. the
action of a Lie group $G$ on $P$ and the induced action on $\sT
P$. If this action is free and proper,  we in fact deal with a
$G$-invariant control system on a principal bundle $P\ra P/G$.
Introducing a $G$-invariant cost function $L:P\times U\lra\R$, one
ends up with a $G$-invariant optimal control problem (OCP) on a
principal bundle $P\ra P/G$.

There are basically two ways of obtaining optimality necessary
conditions for such a problem. In the first, one takes the
Pontryagin maximum principle (PMP) for the unreduced system on $P$
and performs the Poisson reduction of the associated Hamiltonian
equations. For the simple case of an invariant system on a Lie
group $P=G$ (see, e.g., \cite{jurdjevic}) one obtains a system on the Lie
algebra $\g=\sT G/G$, and the reduced Hamilton equations are
the Hamilton equations obtained by means of the Lie--Poisson
structure on $\g^\ast$. The best known example of this type is
probably the reduction for the rigid body in analytical mechanics:
from the cotangent bundle $\sT^\ast\,\SO$ of the group $\SO$
playing the role of the configuration space to the linear Poisson
structure on $\so^\ast$, the dual of the Lie algebra $\so$.
A similar situation appears for homogeneous spaces \cite{jurdjevic} and
general principal bundles \cite{martinez_red_opt_ctr, martinez_lie_classs_mech}. The reduced system
lives on the bundle $\sT P/G$ which is canonically a Lie
algebroid, called the {\em Atiyah algebroid of $P$}, and the
reduced Hamilton equations are associated with the linear
Poisson structure on $\sT^\ast P/G$ (equivalent to the presence of
a Lie algebroid structure on $\sT P/G$). In this approach one
obtains a version of the PMP, yet the Hamiltonian reduction seems
to be purely computational, and a big part of the geometry of the
problem remains hidden.

The second approach, called the \emph{Lagrangian reduction}, was
introduced by Marsden and his collaborators (see, for example,
\cite{cendra_holm_marsden}) in the context of analytical mechanics. Here, one
uses the reduced data $f:P/G\times U\lra \sT P/G$ and $l:P/G\times
U\lra\R$ and the reduced variations (homotopies) to obtain a
reduced version of the Euler-Lagrange equations. In this approach
it becomes clear that the reduction of the variational principle
is not only the reduction of data and geometrical structure but
also a reduction of variations (homotopies) --- this is most
clearly stated in \cite{cendra_holm_marsden} for the case of an invariant system
on a Lie group. By means of the Lagrangian reduction one can
obtain various results, such as Euler-Poincare equations and
Hammel equations. Despite these advantages, the Lagrangian
reduction seems to be more useful in mechanics than in control
theory, as one requires that the geometry of the set of controls
$U$ be very nice (it should be at least an affine subspace of
$\R^n$), so it accepts no discontinuity, switch-on-switch-off
controls, etc.
\subsection{Geometrical mechanics and algebroids}
On the other hand, in geometrical mechanics the suitability of Lie
algebroids as objects which allow a unified treatment of both
reduced and non-reduced systems was recognized some time ago. Lie
algebroids were introduced to mechanics by Weinstein \cite{weinstein} and
Libermann \cite{libermann} and since then studied in different contexts
by many authors \cite{cortes_leon_marrero, cortes_martinez, leon_marrero_martinez, martinez_lagr_mech, martinez_geom_form, martinez_red_opt_ctr,martinez_cft, martinez_lie_classs_mech, martinez_var_calc}. It
was observed a little bit later, following the approach to
analytical mechanics proposed by Tulczyjew, \cite{tulczyjew_ham_lagr,tulczyjew_urbanski_slow}, that
geometrical mechanics, together with the Euler-Lagrange and
Hamilton equations, constraints, etc., can be developed on more
general objects than Lie algebroids \cite{GGU_geom_mech,GG_var_calc}. They were
introduced in \cite{GU_algebroids} under the name \emph{(general)
algebroids}. This generalization turns out to be of practical use,
as the systems of mechanical type with nonholonomic constraints
allow a nice geometrical description in terms of
\emph{skew-algebroids} \cite{grabowski_nonholonomic}, which do not have to satisfy
the Jacobi identity in general.

\subsection{PMP on almost Lie algebroids}
The aim of our work is to extend the fundamental theorem of
optimal control -- the Pontryagin maximum principle (PMP) -- to
the setting of almost Lie (AL) algebroids, geometrical objects
generalizing Lie algebroids. Some differential versions of the PMP
in similar settings have been already indicated in
\cite{cortes_martinez,martinez_red_opt_ctr,GG_var_calc}, but the full version we present here is, to our
knowledge, completely novel.

Since Lie algebroids are infinitesimal (reduced) objects of
(local) Lie groupoids (like Lie algebras are for Lie groups), we
are motivated mostly by the Lie groupoid $\GG$ -- Lie algebroid
$A(\GG)$ reduction. Obviously, a reduction of an invariant control
system on a Lie groupoid should lead to a system on the associated
Lie algebroid. What is more, similarly to the scheme of the
Lagrangian reduction, we should also reduce the variations
(homotopies) from $\GG$ to $A(\GG)$. This will motivate the
abstract definition of the homotopy of admissible paths on an AL
algebroid (algebroid homotopy). Finally, reducing an invariant OCP
on the Lie groupoid $\GG$ would not be complete without reducing
the boundary conditions as well. The natural generalization of a
fixed-end-point boundary condition to a fixed-homotopy-class
condition on $\GG$ (see Section \ref{sec:ocp} for details) allows
us to express the reduced boundary conditions in $A(\GG)$ as
fixing the algebroid homotopy class of the trajectory of the
reduced control system. At the end, for a general AL algebroid
$E$, we can formulate an analogue of the OCP which, in the case of
an integrable algebroid $E=A(\GG)$, turns out to be an  invariant
OCP reduced from $\GG$. Let us note that our  understanding of
algebroid homotopies and homotopy classes is closely related to
that of Crainic and Fernandes \cite{crainic_fernandes}, where similar techniques
were used to generalize Lie's third theorem and integrate Lie
algebroids. However, our framework is much more general, as we no
longer remain in the smooth category.

To briefly explain the result, let us note that an AL algebroid is
a vector bundle $\tau:E\ra M$ together with a vector bundle map
$\rho:E\ra\sT M$ (\emph{anchor}) and a skew-symmetric bilinear
bracket $[\cdot,\cdot]$ on the space of sections of $E$ which
satisfy certain compatibility conditions. The algebroid structure
on $E$ is equivalent to the presence of a certain linear bi-vector
field $\Pi$ on the dual bundle $E^\ast$. Standard examples of an
AL (in fact Lie) algebroid are the following: the tangent bundle $E=\sT M\ra M$
with $\rho=\id_{\sT M}$ and the Lie bracket of vector fields, and
a finite-dimensional real Lie algebra $E=\g$ with the trivial anchor
map ($M$ is a single point in this case) and the Lie bracket on
$\g$. In the first case, $\Pi$ is the canonical Poisson tensor on
$\sT^\ast M$, whereas in the second it is the Lie--Poisson structure
on $\g^\ast$. An example of an AL algebroid which is not a Lie
algebroid is given by any real vector bundle with a smooth family
of skew-symmetric bilinear (but not Lie) operations on its fibers.
Note that the bivector field $\Pi$ defines the \emph{Hamiltonian
vector field} $\X_H$ associated with any $C^1$-function $H$ on
$E^\ast$, defined in the standard way as the contraction
$\X_H=\iota_{\dd H}\Pi$.

On the bundle $E$ we can consider \emph{admissible paths}, i.e.
bounded measurable  maps $a:[t_0,t_1]\ra E$ such that the
projection $x(t)=\tau(a(t))$ of $a(t)$ onto $M$ is absolutely
continuous (AC) and $\dot x(t)=\rho(a(t))$ a.e. On admissible
paths we have an equivalence relation $a\sim b$ interpreted as a
reduction of homotopy equivalence (with fixed end-points). Note
that equivalent paths need not be defined on the same time
interval. For an admissible path $\sigma$, we denote with
$[\sigma]$ the equivalence class of $\sigma$.

A \emph{control system} is defined by a continuous map $f:M\times
U\ra E$, where $U$ is a topological space of control parameters,
such that for each $u\in U$ the function $f(\cdot,u)$ is a section
of class $C^1$ of the bundle $E$. Every \emph{admissible control},
i.e. a bounded measurable path $u(t)$ in $U$, gives rise to an
absolutely continuous path in $M$ defined by the differential
equation
\begin{equation*}
\dot x(t)=\rho\left(f(x(t),u(t))\right)
\end{equation*}
and to an admissible path $a(t)=f(x(t),u(t))$ covering $x(t)$ (we
will call $a(t)$ the \emph{trajectory} of the control system and
the pair $(x(t),u(t))$ the \emph{controlled pair}). An
\emph{optimal control problem} (OCP) for this control system is
associated with a fixed homotopy class $[\sigma]$ of an admissible
path $\sigma$ and a \emph{cost function} $L:M\times U\ra\R$. The
problem is to find a controlled pair $(x(t),u(t))$ with
$t\in[t_0,t_1]$ (the time interval is to be found as well) such
that
\begin{equation}\tag{P}\label{eqn:PPP} \begin{split}\text{the integral $\int_{t_0}^{t_1}L\big(x(t),u(t)\big)\dd t$ is
minimal among all controlled pairs $(x,u)$ for}\\\text{ which the
$E$-homotopy class of the trajectory $f(x(t),u(t))$ equals
$[\sigma]$.}\end{split}
\end{equation}

\begin{theorem}
Let $(x(t),u(t))$, with $t\in[t_0,t_1]$, be a controlled pair
solving the OCP \eqref{eqn:PPP}. Then there exist a
curve $\xi:[t_0,t_1]\lra E^\ast$ covering $x(t)$ and a constant
$\ul\xi_0\leq 0$ such that the following hold
\begin{itemize}
    \item the curve $\xi(t)$ is a trajectory of the time-dependent family of Hamiltonian vector fields
    $\X_{H_t}$,
    $H_t(x,\xi):=H(x,\xi,u(t))$, where
    $$H(x,\xi,u)=\< f\left(x,u\right), \xi>+\ul\xi_0 L\left(x,u\right);$$
    \item the control $u$ satisfies the ``maximum principle''
    $$H(x(t),\xi(t),u(t))=\sup_{v\in U}H(x(t),\xi(t),v)$$
and $H(x(t),\xi(t),u(t))=0$ at every regular point $t$ of $u$;
    \item if $\ul\xi_0=0$, then the covector $\xi(t)$ is nowhere-vanishing.
\end{itemize}
\end{theorem}

\subsection{Discussion of the result and the proof}
The above result looks quite similar to the standard one. Indeed,
in the case $E=\sT M$ we obtain the standard PMP. The only
difference is that the fixed-end-point boundary condition is
substituted by the fixed-homotopy-class condition. However, this
makes no essential difference, as is discussed in detail in
Section \ref{sec:ocp}. For the case of an integrable Lie algebroid
$E=A(\GG)$ our version of PMP can be understood as a general
reduction scheme for invariant OCPs on Lie groupoids. In
particular, the theorem covers the known results on Hamiltonian
reduction of Jurdjevic \cite{jurdjevic} and Martinez \cite{martinez_red_opt_ctr, martinez_lie_classs_mech}
and Lagrangian reduction \cite{cendra_holm_marsden} (see Section
\ref{sec:exmples} for details). It is, however, worth mentioning
that in our approach the reduced and unreduced PMPs are parts of
the same universal formalism. Roughly speaking, we have
generalized the geometrical context in which the PMP can be used.
The technical setting remains quite general -- we work with
general bounded measurable controls and AC base trajectories.
Moreover, since AL algebroids do not come, in general, from
reductions, our result admits a wider spectrum of possible
applications. An attempt in this direction can be found in  the
last example of Section \ref{sec:exmples}.

The strategy of the proof mirrors the classical approach of
Pontryagin and his collaborators \cite{pontryagin}. We follow similar
steps, construct needle variations, and study the geometry of the
set of infinitesimal displacements obtained by these variations,
yet we have to overcome serious technical difficulties. These are
caused mainly by the measurable function setting and the general
approach to algebroid homotopy. In short, as AL algebroids are in
general not integrable, we cannot represent the homotopy classes
as curves on a finite-dimensional manifold and, to make the final
topological argument, we have to pass through infinite-dimensional
Banach spaces.

\subsection{Organization of the paper}
Sections \ref{sec:ala} and \ref{sec:meas} have an introductory
character. We give the definition of AL algebroids,
study the Hamiltonian vector fields and morphism in this context,
and recall basic results from the theory of measurable functions
and measurable ODEs. In Section \ref{sec:htp} we introduce the
notion of homotopy of admissible paths on AL algebroids. Two
important results are: Theorem \ref{thm:int_htp}, which gives an
interpretation of an algebroid homotopy as a reduction of the
standard homotopy by the groupoid action for integrable Lie
algebroids $A(\GG)$, and Lemma \ref{lem:gen_E_htp}, which states
that one-parameter families (homotopies) of admissible paths are
algebroid homotopies.  This is exactly the structure of an AL
algebroid which is essential for the former result. In Sections
\ref{sec:ocp} and \ref{sec:pmp} we introduce the OCP on an AL
algebroid and state our main result. Quite important is the
comparison of the standard OCP and the algebroid OCP given in
Remark \ref{rem:ocp_red} and on the preceding page. In Section
\ref{sec:ocp} we also introduce the algebroid homotopy associated
with a control system. The last three sections are devoted to the
proof of the PMP. The part in the last section mirrors the standard
proof from \cite{pontryagin}, whereas the important partial results from
Theorems \ref{thm:1st_main} and \ref{thm:separation_K_Lambda} use much of the
structure of the AL algebroid. Finally, Section \ref{sec:exmples}
is devoted to examples (including an elegant derivation of the
falling cat theorem of Montgomery \cite{montgomery_isohol} and an application
to nonholonomic systems) as well as relations of our work to other
results.

\section{Almost Lie algebroids}\label{sec:ala}

\subsection{Local coordinates and notation}
Let $M$ be a smooth manifold, and let $(x^a)$, $a=1,\dots,n$, be a
coordinate system in $M$. We denote by $\tau_M: \sT M
\rightarrow M$ the tangent vector bundle and by $\pi_M \colon
\sT^\ast M\rightarrow M$ the cotangent vector bundle.  We have the induced (adapted) coordinate systems $(x^a, {\dot x}^b)$ in $\T
M$ and $(x^a, p_b)$ in $\T^\ast M$.\index{local coordinates}

More generally, let $\tau: E \lra M$ be a vector bundle, and let $\pi: E^\ast \ra M$ be the dual bundle. Choose $(e_1,\dots,e_m)$ --- a basis of local sections of $\tau: E\ra M$, and let $(e^{1}_*,\dots, e^{m}_*)$ be the dual
basis of local sections of $\pi: E^\ast\lra M$. We have the
induced coordinate systems: $(x^a, y^i),  y^i=\iota(e^{i}_*)$ in $E$, and $(x^a, \xi_i), \xi_i = \iota(e_i)$ in $E^\ast$, where the linear functions  $\iota(e)$ are given by the canonical pairing $\iota(e)(v_x)=\< e(x),v_x>$.
The null section of $\tau:E\lra M$ will be denoted by $\theta$, and $\theta_x$ will stand for the null vector at point $x\in M$.

In this work the summation convention is assumed. 

\subsection{Algebroids}
\begin{definition} A \emph{skew-algebroid structure} on a vector bundle $\tau:E\lra M$ is given by a skew-symmetric bilinear bracket $[\cdot,\cdot]$ on the space $\Sec(E)$ of (local) sections of $\tau$, together with a vector bundle morphism $\rho:E\lra \sT M$ (the anchor map), such that the Leibniz rule
\begin{equation}\label{eqn:lieb_rule}
[X,f\cdot Y]=f[X,Y]+\rho(X)(f)Y
\end{equation}
is satisfied for every $X,Y\in\Sec(E)$ and $f\in C^\infty(M)$. If,
additionally, the anchor map is an algebroid morphism, i.e.
\begin{equation}\label{eqn:ala}
\rho\left([X,Y]\right)=[\rho(X),\rho(Y)]_{\sT M},
\end{equation}
we will speak of \emph{almost-Lie algebroids} (\emph{AL
algebroids} shortly).

If, in addition to \eqref{eqn:lieb_rule} and \eqref{eqn:ala}, the
bracket  satisfies the Jacobi identity (in other words the pair
$(\Sec(E)$,$[\cdot,\cdot])$ is a Lie algebra), we speak of a
\emph{Lie algebroid}.
\end{definition}

\smallskip
In the context of mechanics it is convenient to think about an
algebroid as a generalisation of the tangent bundle.
An element $a\in E$ has the interpretation of a generalized
velocity with actual velocity $v\in\sT M$ obtained by applying the
anchor map $v=\rho(a)$. The kernel of the anchor map represents
inner degrees of freedom.

Two basic examples of skew-algebroids are the tangent bundle $\sT
M$ (with the standard Lie bracket and $\rho=\id_{\sT M}$) and a
finite-dimensional real Lie algebra $\mathfrak{g}$ considered as a
vector bundle over a single point, with the trivial anchor and its
Lie bracket. The Lie algebra of a Lie group $G$ can be seen as a
reduction of the tangent bundle $\sT G$ by the left (or right)
action of $G$. Later we will closely investigate another
example of this kind. Namely, for any principal bundle $G\lra
P\lra M$, the reduced bundle $\sT P/G\lra M$ has the structure of
an \emph{Atiyah algebroid}.  All the above are examples of Lie
algebroids. However, it has been shown recently \cite{grabowski_nonholonomic} that
some skew-algebroids (not necessarily Lie) can be naturally
associated with nonholonomic mechanical systems and can be used to
derive the nonholonomic Euler-Lagrange equations by means of a
general geometric scheme \cite{GGU_geom_mech}.

In local coordinates $(x^a,y^i)$, introduced at the beginning of
this section the structure of an algebroid on $E$\index{skew-algebroid!local description} can be
described in terms of local function $\rho^a_i(x)$ and $c^i_{jk}(x)$ on
$M$ given by
$$\rho(e_i)=\rho^a_i(x)\pa_{x^a}\quad \text{and} \quad [e_i,e_j]=c^k_{ij}(x)e_k.$$
The skew-symmetry of the algebroid bracket results in the
skew-symmetry of $c^i_{jk}$ in lower indices, whereas condition
\eqref{eqn:ala} reads as
$$\left(\frac{\pa}{\pa
x^b}\rho^a_k(x)\right)\rho^b_j(x)-\left(\frac{\pa}{\pa
x^b}\rho^a_j(x)\right)\rho^b_k(x)=\rho^a_i(x)c^i_{jk}(x).$$

\subsection{Hamiltonian vector fields and tangent lifts} Let us now describe some geometric constructions
associated with the structure of a skew-algebroid .

It can be shown (cf. \cite{GU_algebroids,GU_poiss_nijn}) that the presence of the
structure of a skew-algebroid on $E$ is equivalent to the
existence of a linear bivector field $\Pi_{E^\ast}$ on $E^*$\index{skew-algebroid!as a linear bi-vector}. In local
coordinates, $(x^a,\xi_i)$ on $E^\ast$, it is given by
\begin{equation*}
 \Pi_{E^\ast} =c^k_{ij}(x)\xi_k
\partial _{\xi_i}\wedge \partial _{\xi_j} + \rho^b_i(x) \partial _{\xi_i}
\wedge \partial _{x^b}.
\end{equation*}
The linearity of $\Pi_{E^\ast}$ means that the corresponding mapping
$\wt\Pi:\sT^\ast E^\ast\lra\sT E^\ast$ is a morphism of double
vector bundles (cf. \cite{KU_dvb, GR_higher}). The tensor $\Pi_{E^\ast}$ is well
recognised in the standard situations: for the tangent algebroid
structure on $\T M$, it is the canonical Poisson structure on
$\T^\ast M$ dual to the canonical symplectic structure,
whereas for a Lie algebra $\g$, it is the
Lie--Poisson structure $\Pi_{\g^\ast}$ on $\g^\ast$. Actually, $E$ is a Lie
algebroid if and only if $\Pi_{E^\ast}$ is a Poisson tensor.

Now we can introduce the notion of a Hamiltonian vector field on
$E^\ast$. Let, namely, $h:E^\ast\lra\R$ be any $C^1$-function. We
define the \emph{Hamiltonian vector field}\index{Hamiltonian vector field} $\X_h$ in an obvious
way: $\X_h=\iota_{\dd h}\Pi_{E^\ast}$. In local coordinates,
\begin{equation}\label{eqn:ham_vf}
\X_h(x,\xi)=\rho^a_i(x) \frac{\pa h}{\pa
{\xi_i}}(x,\xi)\pa_{x^a}+\left(c^k_{ji}(x)\xi_k\frac{\pa h}
{\partial {\xi_j}}(x,\xi)- \rho^a_i(x)\frac{\pa h}{\pa
{x^a}}(x,\xi)\right) \pa_{\xi_i}.
\end{equation}

Another geometrical construction in the skew-algebroid setting is
the \emph{complete lift of an algebroid section}\index{complete lift} (cf.
\cite{GU_algebroids,GU_poiss_nijn}). For every $C^1$-section $X=f^i(x)e_i\in\Sec(E)$
we can construct canonically a vector field $\dd_\T(X)\in\Sec(\T
E)$ which, in local coordinates, reads as
\begin{equation}\label{eqn:tan_lift}
\dd_\T(X)(x,y) = f^i(x)\rho^a_i(x)\pa_{x^a} + \left( y^i
\rho^a_i(x) \frac{\pa f^k}{\pa x^a}(x) + c^k_{ij}(x)y^if^j(x)
\right) \pa _{y^k}.
\end{equation}
The vector field $\dd_\T(X)$ is linear w.r.t. the vector bundle structure $\T\tau:\T E\ra\T M$ (the above equation is linear w.r.t. $y^i$).

Consider the Hamiltonian vector field $\X_{\iota(X)}$ associated
with a linear function $\iota(X)(\cdot)=\<X,\cdot>_\tau$ on
$E^*$. It turns out that fields $\X_{\iota(X)}$ and $\dd_\sT(X)$
are related by
\begin{equation}\label{eqn:hvf_tgl}
\<\dd_\sT(X),\X_{\iota(X)}>_{\sT\tau}=0,
\end{equation}
where $\<\cdot,\cdot>_{\sT\tau}:\sT E\times_{\sT M}\sT
E^\ast\lra\R$ is the canonical pairing, being the tangent map of
$\<\cdot,\cdot>_\tau:E\times_M E^\ast\lra\R$ (in local coordinates,
$\<(x,y,\dot x,\dot y),(x,\xi,\dot x,\dot\xi)>_{\sT\tau}=\dot
y^j\xi_j+y^j\dot \xi_j$).

\subsection{Cartan Calculus}
The existence of a skew-algebroid structure on $E$ is equivalent
to the existence of exterior differential (de Rham) operators\index{exterior differential}
$$d_E:\Sec(\Lambda^kE^*)\lra\Sec(\Lambda^{k+1}E^*),\quad k=0,1,\dots\,,$$
defined by a straightforward generalisation of the Cartan formula
\begin{equation*}
\begin{split}d_E\omega (a_0,a_1,\hdots,a_k)=\sum_{i=0}^k(-1)^i\rho(a_i)\omega(a_0,\hdots,\check{a}_i,\hdots,a_k)
\\ +\sum_{1\leq i<j\leq
k}(-1)^{i+j}\omega\left([a_i,a_j],a_0,\hdots,\check{a}_i,\hdots,\check{a}_j,\hdots,a_k\right),
\end{split}
\end{equation*}
for $\omega\in\Sec(\Lambda^kE^*)$ and
$a_0,a_1,\hdots,a_k\in\Sec(E)$.

These operators, in general, needs not be cohomological. In fact,
$E$ is a Lie algebroid if and only if $d_E^2=0$. AL algebroids, in
turn, can be characterized by the condition that $d_E^2f=0$ for
every $f\in C^\infty(M)=\Sec(\Lambda^0E^*)$.

\subsection{Morphisms}
The above concept of the de Rham derivative allows one to give
a simple definition of a morphism of skew-algebroids. Namely,
given skew-algebroids $\wt\tau:\wt E\lra\wt M$ and $\tau:E\lra M$,
a bundle map $\Phi:\wt E\lra E$ over $\varphi:\wt M\lra M$ is a
\emph{skew-algebroid morphism}\index{morphism of skew-algebroids} if it is compatible with the exterior derivative:
\begin{equation}\label{eqn:E_morph_forms}
\Phi^*d_E\theta=d_{\wt E}\Phi^*\theta,\quad \text{for every
$\theta\in\Sec(\Lambda^kE^*)$}.
\end{equation}
Note that a vector bundle map $\Phi$ does not, in general, induce
any map on sections of $\wt E$, while the pull-back $\Phi^*$ of
sections of $E^*$ is always well defined.

Introduce local coordinates $(\wt x^\alpha,\wt y^\iota)$ and
$(x^a,y^i)$ and structure functions $\wt\rho^\alpha_\iota(\wt x)$, $\wt
c^\iota_{\kappa\mu}(\wt x)$ and $\rho^a_i(x)$, $c^i_{km}(x)$ on $\wt E$ and
$E$, respectively. The condition that
$\Phi\sim(\Phi^i_\iota,\varphi^a)$ is an algebroid morphism reads as
\begin{equation}\label{eqn:alg_morph}
\begin{split}&\Phi^i_\kappa(\wt x)\rho^a_i(\varphi( x))=\wt \rho^\alpha_\kappa(\wt x)\frac{\pa\varphi^a(\wt x)}{\pa \wt x^\alpha},\\
&\wt\rho^\alpha_\kappa(\wt x)\frac{\pa\Phi^i_\lambda(\wt x)}{\pa \wt
x^\alpha}-\wt\rho^\alpha_\lambda(\wt x)\frac{\pa\Phi^i_\kappa(\wt
x)}{\pa \wt x^\alpha}=c^i_{jk}\left(\varphi(\wt
x)\right)\Phi^j_\kappa(\wt x)\Phi^k_\lambda(\wt x).
\end{split}
\end{equation}

\subsection{Admissible paths} Consider an algebroid morphism $\sT\R|_I\lra E$, where $I=[t_0,t_1]\subset\R$ is an interval.
Every such map is uniquely determined by the image of the
canonical section $(t,\pa_t)$ of $\sT\R$ being a smooth curve
$a(t)$ in $E$ over the base path $x(t)$ in $M$. Condition
\eqref{eqn:E_morph_forms} reads as
\begin{equation}\label{eqn:adm}
\rho\left(a(t)\right)=\dot{x}(t) \quad \text{for every $t\in I$}.
\end{equation}
This means that the anchor map coincides with the tangent
prolongation of the projection $x(t)=\tau\left(a(t)\right)$. The
curves which satisfy \eqref{eqn:adm} will be called
\emph{admissible}\index{admissible path}\index{E-path|see{admissible path}}. In fact, \eqref{eqn:adm} also makes sense for non-smooth maps. From now on, by an \emph{admissible path on $E$} (or briefly \emph{$E$-path}) we shall mean a bounded measurable map $a:I\lra E$ over an absolutely continuous (AC) base
path $x=\tau\circ a:I\lra M$ such that \eqref{eqn:adm} is
satisfied a.e. in $I$. In such a case we will speak
of \emph{measurable $E$-paths}. For more information on measurable
functions see Section \ref{sec:meas}. Observe that from \eqref{eqn:alg_morph} it follows that a morphism of algebroids maps admissible paths into admissible paths.

To explain the meaning of admissible curves, observe that in the case
of the tangent algebroid $\T M$  admissible curves are precisely
the tangent lifts of base curves. We will show later (cf. Theorem
\ref{thm:int_htp}) that if an algebroid $E$ is integrable,
admissible curves come from a reduction of real curves in a Lie
groupoid integrating $E$.

Finally, we can introduce the concept of \emph{composition of
measurable $E$-paths}\index{composition of admissible paths}. Let $a:[t_0,t_1]\lra E$ and $\ol
a:[t_1,t_2]\lra E$ be two measurable $E$-paths with base paths
$x=\tau\circ a$ and $\ol x=\tau\circ\ol a$, respectively. Assume
that $x(t_1)=\ol x(t_1)$ (such paths will be called
\emph{composable}\index{composable admissible paths}). Clearly, the map $\wt a:[t_0,t_2]\lra E$
defined by
$$\wt a(t):=\begin{cases}
a(t) &\text{for $t\leq t_1$},\\
\ol a(t) &\text{for $t>t_1$}
\end{cases}$$
is another measurable $E$-path covering the AC curve
$$\wt x(t):=\begin{cases}
x(t) &\text{for $t\leq t_1$},\\
\ol x(t) &\text{for $t>t_1$}.
\end{cases}$$
This new $E$-path will be called the \emph{composition} of $a$ and
$\ol a$ and will be denoted by $\wt a=a\circ \ol a$.

\subsection{The product of skew-algebroids}
\index{product of skew-algebroids}
Given two skew-algebroids $(\tau_1:E_1\lra
M_1,\rho_1,[\cdot,\cdot]_1)$ and $(\tau_1:E_2\lra
M_2,\rho_2,[\cdot,\cdot]_2)$ we can define a skew-algebroid
structure on the product bundle $\tau=\tau_1\times\tau_2:E_1\times
E_2\lra M_1\times M_2$. The anchor will simply be
$\rho=\rho_1\times\rho_2:E_1\times E_2\lra\T M_1\times\T
M_2\approx \T(M_1\times M_2)$. The bracket can be defined by
equalities
\begin{align*}
\left[\p_1^\ast X_1,\p_1^\ast Y_1\right]&=\p_1^\ast [X_1,Y_1]_1,\\
\left[\p_2^\ast X_2,\p_2^\ast Y_2\right]&=\p_2^\ast [X_2,Y_2]_2,\\
\left[\p_1^\ast X_1,\p_2^\ast Y_2\right]&=\theta,
\end{align*}
where $X_1,Y_1\in\Sec(E_1)$ and $X_2,Y_2\in\Sec(E_2)$ are sections,  $\p_1:E_1\times E_2\lra E_1$ and $\p_2:E_1\times E_2\lra E_2$ are canonical vector bundle projections, and $\theta$ is a null section of $\tau$. The above equalities can be extended to arbitrary
sections by linearity and the Leibniz rule \eqref{eqn:lieb_rule}.
Clearly, the canonical projections $E_1\times E_2\lra E_i$, with
$i=1,2$, are algebroid morphisms, and if $E_1$ and $E_2$ are almost
Lie, then so is their product.

The local coordinate description of the product $E_1\times E_2$ is
very simple. If $(x^a,y^i)$ and $(\wt x^\alpha,\wt y^\iota)$ are
local coordinates on $E_1$ and $E_2$, respectively, we can
introduce natural coordinates $(X^A,Y^I)=(x^a,\wt x^\alpha,
y^i,\wt y^\iota)$ on $E_1\times E_2$. The structure functions
$C^I_{JK}(X)$ and $R^A_I(X)$ in these coordinates are trivial on
mixed-type terms ($R^\alpha_i=C^\iota_{j\kappa}=0$, etc.) and the
same as the structure functions of $E_1$ and $E_2$ on simple-type
terms ($C^i_{jk}(x,\wt x)=c^i_{jk}(x)$, $R^\alpha_\iota(x,\wt
x)=\rho^\alpha_\iota(\wt x)$, etc.).

\section{Measurable maps}\label{sec:meas}
In this section we will briefly recall some basic properties of
measurable and AC maps as well as ODEs in the measurable setting. Our discussion
in based mostly on \cite{bressan}.

\subsection{Basic facts} 
When speaking about measure we will always have in mind \emph{Lebesgue measure}\index{Lebesgue measure} in $\R^n$ or subsets of $\R^n$. This measure will be denoted by $\mu_L(\cdot)$. 

Recall that a map $f:V\supset\R^n\ra\R^k$, defined on a subset $V\subset\R^n$,
is \emph{measurable}\index{measurable map} if the inverse image of every open set is
Lebesgue measurable in $V$. The measurable map $f$ will be called
\emph{bounded}\index{measurable map!bounded} if the closure of its image is a compact set.
Observe that every bounded (or locally bounded) measurable
function is locally integrable. 

\noindent Measurable maps can be characterised as follows.

\begin{theorem}[Luzin]\label{thm:luzin}\index{Luzin theorem}
The map $f:\R^n\supset V\ra\R$ defined on a measurable set $V$ is measurable iff, for every $\eps>0$, there exists a closed subset $F\subset V$ such that
the restriction $f|_F$ is continuous and $\mu_L(V\setminus F)<\eps$.
\end{theorem}
\noindent For the proof see \cite{lojasiewicz}.

\noindent In our considerations much attention will be payed to regular points of measurable maps.
\begin{definition}
Let $f:[a,b]\ra\R^m$ be a measurable map. A point $x\in [a,b]$ is called a
\emph{regular point} (also: \emph{Lebesgue} or \emph{density point})\index{regular point} of $f$, iff  
$$\lim_{t\to
0}\frac1{|t|}\int_0^t|f(x+s)-f(x)|\dd s=0.$$
\end{definition}

\noindent For bounded measurable (or more generally integrable) maps we have the following result.

\begin{theorem}[Lebesgue]\index{Lebesgue theorem}
For an integrable map $f:[a,b]\ra\R^m$ almost every point in $[a,b]$ is a regular point of $f$.
\end{theorem}
\noindent For the proof see \cite{lojasiewicz}

A map $x:[a,b]\ra\R^k$ is called \emph{absolutely
continuous}\index{absolutely continuous map}\index{AC map|see{absolutely continuous map}}  (AC) if it can be written in the form
$$x(t)=x(a)+\int_{a}^{t}v(\tau)\dd \tau,$$
where $v(\cdot)$ is an integrable map. As we see,
an AC map $x(t)$ is differentiable at all the regular points $t$
of $v$ (hence, by Lebesgue theorem, differentiable a.e.). Its
derivative at such a point is simply $v(t)$. In this work we concentrate our attention mostly on \emph{absolutely 
continuous maps with bounded derivative}\index{absolutely continuous map!with bounded derivative}\index{ACB map|see{absolutely continuous map with bounded derivative}}  (ACB maps); i.e., maps for which $v(\cdot)$ is bounded measurable.

In our considerations we will use the following lemma.
\begin{lemma}\label{lem:reg}
Let $a:[0,1]\ra\R$ be a bounded measurable map, let $h:[0,1]\ra[0,1]$ be a continuous function, and let $g:[0,1]\ra[0,1]$ be a $C^1$--map with a non-vanishing derivative. Then the map
$$s\longmapsto G(s):=\int_0^1\left|a(g(t)h(s))-a(g(t)h(s_0))\right|\dd t$$
is regular (in fact continuous) at every $s_0$ such that $h(s_0)\neq 0$.
\end{lemma}
\begin{proof}
Let $c$ be a number such that $0<c\leq\left|g^{'}(x)\right|$ for every $x\in[0,1]$. Now if $A\subset[0,1]$ is a measurable subset then
$\mu_L\left(g^{-1}( A)\right)\leq \frac 1c\mu_L(A)$. 

Choose $\eps>0$. By Luzin Theorem \ref{thm:luzin} there exists a closed set $F\subset[0,1]$ such that $a(\cdot)$ is continuous on $F$ and $\mu_L([0,1]\setminus F)<\eps$. Now $a$ is uniformly continuous on $F$, $g$ is bounded and $h$ continuous, hence there exists $\delta>0$ such that $\left|a(g(t)h(s))-a(g(t)h(s_0))\right|<\eps$ if only $|s-s_0|<\delta$ and $t$ and $s$ are such that $g(t)h(s)\in F$ and $g(t)h(s_0)\in F$. As a consequence for $|s-s_0|<\delta$, we can estimate
\begin{align*}
&\int_0^1|a(g(t)h(s))-a(g(t)h(s_0))|\dd t \leq
\int_{\{t:g(t)h(s)\notin F\}}2\|a\|\dd t+\int_{\{t:g(t)h(s_0)\notin F\}}2\|a\|\dd t+\\
&+\int_{\{t:g(t)h(s)\in F\}\cap\{t:g(t)h(s_0)\in F\}}|a(g(t)h(s))-a(g(t)h(s_0))|\dd t\\
&\leq 2\|a\|\cdot\mu_L\left(g^{-1}\left(\frac 1{h(s)}([0,1]\setminus F)\right)\right)+2\|a\|\cdot\mu_L\left(g^{-1}\left(\frac 1{h(s_0)}([0,1]\setminus F)\right)\right)+\int_0^1\eps\dd t\\
&\leq\eps\left(2\|a\|\frac 1 c\left(\frac 1{ h(s)}+\frac 1{h(s_0)}\right)+1\right)
\end{align*}
Since $h(s_0)\neq 0$, if $|s-s_0|$ is sufficiently small, the values s of $G(s)$ are arbitrarily close to $0=G(s_0)$, which finishes the proof. \end{proof}

\subsection{Uniform regularity} 
Regular points play an important role in our considerations, since the behaviour of a measurable map at a regular point is similar to the behaviour of a continuous map. To study behaviour of  the families of measurable maps we introduce a notion of \emph{uniform regularity}.

\begin{definition}\label{def:ur}
Let $P$ be a topological space and consider a map $f:[a,b]\times P\ra\R^m$ such that $t\mapsto f(t,p)$ is a measurable  for every $p\in P$. We call $f$ \emph{uniformly regular with respect to $p\in P$ at $x\in[a,b]$}\index{uniformly regular map} iff the following conditions are satisfied:
\begin{align}
&\frac 1{|t|}\int_0^t\left|f(x+s,p)-f(x,p)\right|\dd s\underset{t\to 0}\lra0\quad\text{locally uniformly w.r.t. $p$,}\label{eqn:ur1}\\
&\text{the map}\quad p\mapsto f(x,p) \quad\text{is continuous},\label{eqn:ur2}
 \intertext{and for every compact set $K\subset P$ there exists a number $t_0>0$ such that}
 & p\mapsto\Big([0,t_0]\ni s\mapsto f(x+s,p)\Big)\quad \text{is a continuous map from $K$ to $L^1([0,t_0],\R^m)$}.\label{eqn:ur3}
\end{align} 
\end{definition}

Usually in mathematics the word ''uniform'' means ''in the same way for all parameters''. In the context of regularity this can be expressed by the condition \eqref{eqn:ur1} itself. Therefore Definition \ref{def:ur} is more specific then what one could expect under the name ''uniform regularity''. The sense of this definition is, however, to abstract several technical properties of measurable  maps which are important from the point of view of this work. Since, according to our knowledge, the notion of uniform regularity is not a well established term, we hope that Definition \ref{def:ur} would not be confusing.

Let us now investigate some simple properties of uniformly regular maps. In what follows we will consider only uniform regularity at point $0\in\R$ and restrict our attention to parameter spaces $P$ which are metric spaces (we can think of $P$  as of  a subset of $\R^m$). 

A basic example of a uniformly regular map is just a continuous map.\index{uniformly regular map!properties}
\begin{proposition}\label{prop:ur_cont}
Let $F:\R\times P\lra\R$ be a continuous map. Then $F$ is uniformly regular w.r.t. $p\in P$ at $s=0$. 
\end{proposition}
\begin{proof} Condition \eqref{eqn:ur2} is obvious. Fix now a compact set $K_P\subset P$ and restrict $s$ to a fixed interval $[0,t_0]$. Since $F$ is uniformly continuous on $[0,t_0]\times K_P$, for every $\eps>0$ there exists $\delta>0$ such that $\left|F(s,p)-F(0,p)\right|<\eps$ if $|s|<\del$ and for all $p\in K_P$. Consequently, 
$$\int_0^t\left|F(s,p)-F(0,p)\right|\leq|t|\eps$$
for $|t|<\delta$ and all $p\in K_P$. 
This proves \eqref{eqn:ur2}, i.e., $\frac 1{|t|}\int_0^t\left|F(s,p)-F(0,p)\right|\underset{t\to 0}\lra 0$ uniformly w.r.t. $p\in K_P$. 

To check \eqref{eqn:ur3} observe that, by the uniform continuity of $F$ on $[0,t_0]\times K_P$, for every $\eps>0$ there exists $\del>0$ such that $\left|F(s,p)-F(s,p^{'})\right|\leq\eps$ for every $s\in[0,t_0]$ and all $p,p^{'}\in K_P$ such that $|p-p^{'}|<\del$. Consequently,
$$\int_0^{t_0}\left|F(s,p)-F(s,p^{'})\right|\dd s\leq\left|t_0\right|\eps$$
for all $p,p^{'}\in K_P$ such that $|p-p^{'}|<\del$. This proves \eqref{eqn:ur3}.
\end{proof}

\noindent Another simple example is the following.

\begin{proposition}\label{prop:ur_trivial}
Let $f:\R\ra\R^m$ be a measurable map regular at $s=0$. For $p\in P$ define $\wt f(s,p):=f(s)$. Then the map $\wt f$ is uniformly regular w.r.t. $p\in P$ at $s=0$.  
\end{proposition}
\begin{proof}
Conditions \eqref{eqn:ur1}--\eqref{eqn:ur3} are trivially satisfied. 
\end{proof}

\noindent Below we discuss several ways of generating uniformly regular maps from given ones. 

\begin{proposition}\label{prop:ur_sum}
Let $f,g:\R\times P\ra\R^m$ be two maps uniformly regular w.r.t. $p\in P$ at $s=0$. Then the sum $f+g$ is also uniformly regular w.r.t $p\in P$ at $s=0$.
\end{proposition}
\begin{proof} The proof is straightforward.
\end{proof}

\begin{proposition}\label{prop:ur_multiplication}
Let $f:\R\times P\ra\R^m$ be bounded and uniformly regular w.r.t. $p\in P$ at $s=0$. Let $h:P\lra\R$ be a continuous map. Then the map $\wt f(s,p)=h(p)f(s,p)$ is uniformly regular w.r.t. $p$ at $s=0$.
\end{proposition}
\begin{proof}
Condition \eqref{eqn:ur2} is obvious. 
Choose now a compact set $K\subset P$. For $p\in K$ we have
$$\frac 1{|t|}\int_0^t\left|h(p)f(s,p)-h(p)f(0,p)\right|\dd s\leq \sup_{p\in K}|h(p)|\cdot \frac 1{|t|}\int_0^t\left|f(s,p)-f(0,p)\right|\dd s\underset{t\to 0}\lra 0$$
uniformly w.r.t. $p\in K$, and hence \eqref{eqn:ur1} is satisfied.

Finally, note that for $p,p^{'}\in K$ we have
\begin{align*}
&\int_0^{t_0}\left|h(p)f(s,p)-h(p^{'})f(s,p^{'})\right|\dd s\\
&\leq\left|h(p)\right|\int_0^{t_0}\left|f(s,p)-f(s,p^{'})\right|\dd s+\left|h(p)-h(p^{'})\right|\int_0^{t_0}\left|f(s,p^{'})\right|\dd s\\
&\leq\sup_{p\in K}|h(p)|\int_0^{t_0}\left|f(s,p)-f(s,p^{'})\right|\dd s+\left|h(p)-h(p^{'})\right|\left\|f(\cdot,p^{'})\right\|_{L^1}\underset{p\to p^{'}}\lra 0+0;
\end{align*}
that is, \eqref{eqn:ur3} is satisfied.
\end{proof}

\begin{lemma}\label{lem:ur_rescal}
Let $f:\R\times P\ra\R^m$ be bounded and uniformly regular w.r.t. $p\in P$ at $s=0$. Consider $\wt f(s,p,c):=f(sc,p)$ where $c\in\R$. Then $\wt f:\R\times P\times\R\lra\R^m$ is uniformly regular w.r.t. $p\in P$ and $c\in\R$ at $s=0$. 
\end{lemma}
\begin{proof}
Since $\wt f(0,p,c)=f(0,p)$, condition \eqref{eqn:ur2} is obvious.

Consider now compact sets $K_C\subset\R$ and $K_P\subset P$. For 
 $c\in K_C$ and $p\in K_P$ we have
\begin{align*}
&\frac 1{|t|}\int_0^t\left|\wt f(s,p,c)-\wt f(0,p,c)\right|\dd s=\frac 1{|t|}\int_0^t\left|f(sc,p)-f(0,p)\right|\dd s\\
&=\frac 1{|t|c}\int_0^{tc}\left|f(s^{'},p)-f(0,p)\right|\dd s^{'}
\underset{t\to 0}\lra 0.
\end{align*}
Since $K_C$ is bounded and $\frac 1{|t^{'}|}\int_0^{t^{'}}\left|f(s,p)-f(0,p)\right|\dd s
\underset{t^{'}\to 0}\lra 0$ uniformly w.r.t. $p\in K_P$, the above convergence is uniform w.r.t. $p\in K_P$ and $c\in K_C$. 

We are left with the proof of property \eqref{eqn:ur3}.  We will check that $(p,c)\mapsto\wt f(\cdot,p,c)$; $K_P\times K_C\lra L^1\left([0,\wt t_0],\R^m\right)$ is continuous separately w.r.t. $p$ and w.r.t. $c$ for a suitably chosen $\wt t_0$. 

Let $t_0>0$ be a number from the property \eqref{eqn:ur3} for $f(s,p)$ and $K=K_P$. To prove the continuity w.r.t. $p$ fix $c\in K_C$. If $c\neq 0$, then
\begin{align*}\int_0^t\left|\wt f(s,p,c)-\wt f(s,p^{'},c)\right|\dd s&=\int_0^t\left|f(sc,p)-f(sc,p^{'})\right|\dd s\\
&=\frac 1{c}\int_0^{tc}\left|f(s^{'},p)-f(s^{'},p^{'})\right|\dd s^{'}\underset{p\to p^{'}, }\lra 0
\end{align*}
if only $tc\leq t_0$. 

For $c=0$ we have 
$$\int_0^t\left|\wt f(s,p,c)-\wt f(s,p^{'},c)\right|\dd s=\int_0^t\left|f(0,p)-f(0,p^{'})\right|\dd s=|t|\left|f(0,p)-f(0,p^{'})\right|\underset{p\to p^{'}}\lra 0$$
for every $t$. In particular, we proved continuity w.r.t. $p$  for $\wt t_0:=\frac{t_0}{\sup_{c\in K_C}|c|}$.

Now fix $p\in P$, fix $c^{'}\in K_C$ , choose $\eps>0$ and consider $c\in K_C$. 
If $c^{'}=0$, then
\begin{align*}
&\int_0^{\wt t_0}\left|\wt f(s,p,c)-\wt f(c,p,c^{'})\right|\dd s\\
&=\int_0^{\wt t_0}\left| f(sc,p)-f(0,p)\right|\dd s=
\left|\wt t_0\right|\frac 1{\left|\wt t_0\right|c}\int_0^{\wt t_0 c}\left|f(s^{'},p)-f(0,p)\right|\dd s^{'}\underset{c\to 0}\lra 0.
\end{align*}
If $c^{'}\neq 0$ consider a closed set $F\subset[0,\wt t_0]$ such that $\mu_L\left([0,\wt t_0]\setminus F\right)<\eps \left|\wt t_0\right|$ and $\wt f(\cdot,p,c^{'})$ is continuous on $F$ (note that $p$ and $c^{'}$ are fixed). Such a set exists by Luzin Theorem \ref{thm:luzin}. 

Since $\wt f(\cdot,p,c^{'})$ is uniformly continuous on $F$, there exists a number $\del>0$ such that\\ $\left|\wt f(s,p,c^{'})-\wt f(s^{'},p,c^{'})\right|<\eps$ if $\left|s-s^{'}\right|\leq\delta$ and $s,s^{'}\in F$. Now $\wt f(s,p,c)=\wt f(\frac c{c^{'}}s,p,c^{'})$ and $\left|s-\frac c{c^{'}}s\right|\leq\left|\frac{c^{'}-c}{c^{'}}\right|\left|\wt t_0\right|$ for $s\in[0,\wt t_0]$, so we have 
 $$\left|\wt f(s,p,c^{'})-\wt f(s,p,c)\right|<\eps\quad\text{if}\quad |c^{'}-c|\leq\frac{c^{'}}{\left|\wt t_0\right|}\delta\quad \text{and}\quad s\in F\cap\frac{c^{'}}cF.$$
Note that 
\begin{align*}
&\mu_L\left([0,\wt t_0]\setminus F\cap\frac {c^{'}}c F\right)\leq\mu_L\left([0,\wt t_0]\setminus F\right)+\mu_L\left([0,\wt t_0]\setminus\frac {c^{'}}c F\right)\\
&\leq \eps\left|\wt t_0\right|+\left(\left|1-\frac{c^{'}}c\right|+\frac{c^{'}}c\eps\right)\left|\wt t_0\right|\leq 4\eps\left|\wt t_0\right|
\end{align*} if $|c-c^{'}|$ is small enough. Consequently,
\begin{align*}
&\int_0^{\wt t_0}\left|\wt f(s, p,c)-\wt f(s,p,c^{'})\right|\dd s\\
&\leq\int_{F\cap\frac {c^{'}}c F}\left|\wt f(s, p,c)-\wt f(s,p,c^{'})\right|\dd s+\int_{[0,\wt t_0]\setminus F\cap\frac {c^{'}}c F}\left|\wt f(s, p,c)-\wt f(s,p,c^{'})\right|\dd s\\
&\leq\eps\left|\wt t_0\right|+\mu_L\left([0,\wt t_0]\setminus F\cap\frac {c^{'}}c F\right)\cdot 2\left\|f\right\|\leq \eps\left(\left|\wt t_0\right|+8\left|\wt t_0\right|\left\|f\right\|\right)
\end{align*}
if $|c-c^{'}|$ is small enough. Since $\eps$ is an arbitrary positive number, this proves the continuity of $(p,c)\mapsto\wt f(\cdot,p,c)$ w.r.t. $c$.\end{proof}

\begin{lemma}\label{lem:ur_comp}
Let $f:\R\times P\ra\R^m$ be bounded and uniformly regular w.r.t. $p\in P$ at $s=0$, and let $G:\R^m\times\R\times P\times Q\lra\R^m$ be a continuous map w.r.t all variables. Then the composition $G(f(s,p),s,p,q)$ is uniformly regular w.r.t. $p\in P$ and $q\in Q$ at $s=0$.
\end{lemma}
\begin{proof}
We will prove the assertion for $G$ trivially depending on $q\in Q$. This will suffice, since we can denote $G(f(s,p),s,p,q)$ as $G(\wt f(s,p,q),s,p,q)=G(\wt f(s,\wt p),s,\wt p)$, where $\wt f(s,p,q):=f(s,p)$ and $\wt p:=(p,q)\in P\times Q$. Clearly, $\wt f(s,\wt p)$ is uniformly regular w.r.t. $\wt p=(p,q)\in P\times Q$ at $s=0$ (cf. Proposition \ref{prop:ur_trivial}) and the investigated composition has a desired simpler form $G(\wt f(s,\wt p),s,\wt p)$.  

Property \eqref{eqn:ur2} is obvious. To prove \eqref{eqn:ur1} estimate
\begin{align*}
&\left|G(f(s,p),s,p)-G(f(0,p),0,p)\right|\\
&\leq\left|G(f(s,p),s,p)-G(f(0,p),s,p)\right|+\left|G(f(0,p),s,p)-G(f(0,p),0,p)\right|.
\end{align*}
Since $\wt G(s,p):=G(f(0,p),s,p)$ is continuous w.r.t. $p$ and $s$, it satisfies \eqref{eqn:ur1}. Consequently, it is enough to check if
$$\frac 1{|t|}\int_0^t\left|G(f(s,p),s,p)-G(f(0,p),s,p)\right|\dd s\underset{t\to 0}\lra 0$$
locally uniformly w.r.t. $p$. To prove it consider a compact set $K_P\subset P$ and restrict $s$ to the interval $[0,t_0]\subset\R$. Since $f$ is bounded, its image $\Image f$ is contained in a compact subset $K\subset\R^m$. Fix $\eps>0$. The map $G$ is uniformly continuous on $K\times[0,t_0]\times K_P$, so there exists a number $\delta>0$ such that
$$\left|G(x,s,p)-G(y,s,p)\right|<\eps$$
if $|x-y|<\delta$ and $x,y\in K$, $s\in[0,t_0]$, and $p\in K_P$.  
 Since $f(s,p)$ is uniformly regular there exists a number $0<\wt t_0\leq t_0$ such that 
\begin{equation}\label{eqn:ur4}
\int_0^t\left|f(s,p)-f(0,p)\right|\dd s<|t|\cdot\eps\cdot\delta
\end{equation}
for every $|t|\leq \wt t_0$ and each $p\in K_P$. Define now $A_p:=\{s\in[0,t_0]:\left|f(s,p)-f(0,p)|>\delta\right|\}$. From \eqref{eqn:ur4} we have
$$\delta\cdot\mu_L\left([0,t]\cap A_p\right)\leq\int_{[0,t]\cap A_p}\left|f(s,p)-f(s,p)\right|\dd s\leq\int_0^t\left|f(s,p)-f(s,p)\right|\dd s\leq|t|\cdot\eps\cdot\delta,$$
hence $\mu_L\left([0,t]\cap A_p\right)\leq\eps\cdot|t|$. Consequently, for $|t|<\wt t_0$, we have
\begin{align*}
&\int_0^t\left|G(f(s,p),s,p)-G(f(0,p),s,p)\right|\dd s\leq 2\cdot\sup_{K\times[0,t_0]\times K_P}|G|\int_{[0,t]\cap A_p}1\dd s+\int_{[0,t]\setminus A_p}\eps \dd s\\
&\leq 2\cdot\sup_{K\times[0,t_0]\times K_P}|G|\cdot\mu_L\left([0,t]\cap A_p\right)+\eps|t|=\eps|t|\left( 2\cdot\sup_{K\times[0,t_0]\times K_P}|G|+1\right).
\end{align*}
Since $\eps>0$ was arbitrary, this proves \eqref{eqn:ur1}.

To prove \eqref{eqn:ur3} we proceed similarly. Again we restrict our attention to $K\times[0,t_0]\times K_P$ and fix $\eps>0$. Let $\delta>0$ be such that, for $x,y\in K$, $s\in[0,t_0]$ and $p,p^{'}\in K_P$
$$\left|G(x,s,p)-G(y,s,p^{'})\right|<\eps$$
if $|x-y|<\delta$ and $|p-p^{'}|<\delta$. 

From the uniform regularity of $f(s,p)$, there exists a number $\wt\delta>0$ such that 
$$\int_0^{t_0}\left|f(s,p)-f(s,p^{'})\right|\dd s\leq\eps\cdot \delta$$
if $p,p^{'}\in K_p$ are such that $|p-p^{'}|<\wt\delta$. From that we deduce that the set $B_{pp^{'}}:=\{s\in[0,t_0]:\left|f(s,p)-f(s,p^{'})\right|>\delta\}$ has measure smaller than $\eps$ if $|p-p^{'}|<\wt \delta$. Indeed, we can estimate
$$\delta\cdot\mu_L\left(B_{pp^{'}}\right)\leq\int_{B_{pp^{'}}}\left|f(s,p)-f(s,p^{'})\right|\dd s\leq\int_0^{t_0}\left|f(s,p)-f(s,p^{'})\right|\dd s\leq \delta\cdot\eps.$$

Therefore for $|p-p^{'}|\leq\min\{\delta,\wt\delta\}$ we can estimate
\begin{align*}
&\int_0^{t_0}\left|G(f(s,p),s,p)-G(f(s,p^{'}),s,p^{'})\right|\dd s\leq\int_{B_{pp^{'}}}\left|G(f(s,p),s,p)-G(f(s,p^{'}),s,p^{'})\right|\dd s\\
&\phantom{=}+\int_{[0,t_0]\setminus B_{pp^{'}}}\left|G(f(s,p),s,p)-G(f(s,p^{'}),s,p^{'})\right|\dd s\\
&\leq 2\cdot\sup_{K\times[0,t_0]\times K_P}|G|\cdot\mu_L\left(B_{pp^{'}}\right)+\int_0^{t_0}\eps\dd s\leq\eps\left(2\cdot\sup_{K\times[0,t_0]\times K_P}|G|+|t_0|\right).
\end{align*}
This proves \eqref{eqn:ur3}. \end{proof}

\subsection{Measurable ODEs}\label{ssec:ode}
Consider an ordinary differential equation associated with a map
$g:\R^n\times\R\ra\R^n$,
\begin{equation}\label{eqn:ode}
\dot x(t)=g(x(t),t).
\end{equation}
By a \emph{(Carath{\'e}odory) solution}\index{Carath{\'e}odory solution}\index{measurable solution of ODE|see{Carath{\'e}odory solution}} of \eqref{eqn:ode} on an
interval $I=[t_0,t_1]$ we shall mean an AC map $t\mapsto x(t)$
which satisfies \eqref{eqn:ode} a.e. For the solutions in the
above sense one can develop the standard theory of existence,
uniqueness, and parameter dependence, as done in
\cite{bressan}. Let us recall the most important results of this
theory. Assume the following:
\begin{equation}\label{ass:A}\tag{A}
\text{$t\mapsto g(x,t)$ is measurable for
every $x$,} \text{ and } \text{$x\mapsto g(x,t)$ is continuous for
every $t$;}
\end{equation}
\begin{equation}\label{ass:B}\tag{B}
\text{$g(x,t)$ is locally bounded and locally Lipschitz w.r.t. $x$;}
\end{equation}
that is, for every compact set $K\subset\R^n\times\R$ there exist
constants $C_K$ and $L_K$ such that $|g(x,t)|\leq C_K$ and
$|g(x,t)-g(y,t)|\leq L_K|x-y|$ for every $(x,t),(y,t)\in K$.

\begin{theorem}[existence and uniqueness of solutions]\label{thm:exist}\index{Carath{\'e}odory solution!existence and uniqueness}
Assuming that \eqref{ass:A} and \eqref{ass:B} hold, for every
$x_0\in\R^n$ there exists a unique solution $x(t,x_0)$ of
\eqref{eqn:ode} with the initial condition $x(t_0)=x_0$, defined
on some interval $[t_0,t_0+\eps]$. If $g$ is globally bounded and
globally Lipschitz (so that the constants $C_K$ and $L_K$ in
\eqref{ass:B} can be chosen universally for all $K$'s), then the
solution is also defined globally. Moreover, if $x(t,x_0)$ is defined on the interval $[t_0,t_1]$ then so are the solutions $x(t,x_0^{'})$ for $x_0^{'}$ close enough to $x_0$. 
\end{theorem}
\begin{proof}[Sketch of the proof]
The proof uses the standard Picard's Method. One constructs a contracting map 
$$A_{x_0}:x(t)\longmapsto x_0+\int_{t_0}^tg(x(\tau),\tau)\dd\tau$$
and uses it to define inductively a sequence of functions $x^0(t,x_0)=x_0$, $x^{n+1}(t,x_0)=A_{x_0}(x^n(t,x_0))$ which converges uniformly in $t$ to the solution $x(t,x_0)$. The length of the interval $[t_0,t_1]$ on which the solution is well-defined depends on the Lipschitz bound of $g(x,t)$. The details can be found in \cite[Thm. 2.1.1]{bressan}. 
\end{proof}

\subsection{Parameter dependence} 
Consider now differential equation \eqref{eqn:ode} with an additional parameter dependence
\begin{equation}\label{eqn:ode1}
\dot{x}(t)=g(x(t),t,s),
\end{equation}
where $g:\R^n\times\R\times\R\lra\R^n$. Assume the following:
\begin{equation}\label{ass:A1}\tag{$A^{'}$}
\text{$t\mapsto g(x,t,s)\,,\ s\mapsto g(x,t,s)$
are measurable,} \text{ and } \text{$x\mapsto g(x,t,s)$ is
continuous;}
\end{equation}
\begin{equation}\label{ass:B1}\tag{$B^{'}$}
\text{$g(x,t,s)$ is locally bounded and locally Lipschitz w.r.t. $x$;}
\end{equation}
that is, for every compact set $K\subset\R^n\times\R\times\R$
there exist constants $C_K$, $L_K$ such that $|g(x,t,s)|\leq C_K$
and $|g(x,t,s)-g(y,t,s)|\leq L_K|x-y|$ for every
$(x,t,s),(y,t,s)\in K$.

\begin{theorem}[parameter dependence]\label{thm:param}\index{Carath{\'e}odory solution!regularity}
Assume that \eqref{ass:A1} and \eqref{ass:B1} hold, and denote by
$x(t,x_0,s)$ the solution of \eqref{eqn:ode1} for a fixed
parameter $s$ and the initial condition $x(t_0,x_0,s)=x_0$ (we
know that such solutions locally exist by Theorem
\ref{thm:exist}). Then the dependence $x_0\mapsto x(t,x_0,s)$ is
continuous, whereas, for any bounded measurable map $s\mapsto
x_0(s)$, the map $s\mapsto x(t,x_0(s),s)$ is also bounded and
measurable for every $t$.
\end{theorem}
\begin{proof}[Sketch of the proof]
As before one constructs a sequence $x^n(t,x_0,s)$ defined by means of the contracting map
$$A_{x_0,s}:x(t)\longmapsto x_0+\int_{t_0}^t g(x(\tau),\tau,s)\dd\tau.$$
The sequence converges to the solution $x(t,x_0,s)$ uniformly w.r.t. $t$ and $x_0$, which implies continuity of the solution w.r.t. the initial value. If $x_0(s)$ is measurable w.r.t. $s$, so is the sequence $x^n(t,x_0(s),s)$. The limit $x(t,x_0(s),s)$ is measurable as a point-wise limit of measurable functions. Moreover, since $A_{x_0,s}$ is a contraction, $\|x(t,x_0(s),s)\|$ is bounded by a constant times $\|x_0(s)\|$. Details can be found in \cite{bressan}.\end{proof}

\noindent Assuming higher regularity of  $g(x,t,s)$, one can
prove a stronger result.

\begin{theorem}[differentiability w.r.t. the initial value]\label{thm:param_dif}\index{Carath{\'e}odory solution!regularity}

Assume that the function $g(x,t,s)$ satisfies \eqref{ass:A1} and
\eqref{ass:B1}, it is differentiable w.r.t. $x$, and the derivative
$\frac{\pa g}{\pa x}(x,t,s)$ satisfies \eqref{ass:A1} and is
locally bounded. Then the solution $x(t,x_0,s)$ of
\eqref{eqn:ode1} is differentiable w.r.t. the initial
condition $x_0$. Moreover, the derivative $\frac{\pa x}{\pa
x_0}(t,x_0,s)$ is continuous in $x_0$, AC in $t$, and measurable
in $s$.
\end{theorem}
\begin{proof}[Sketch of the proof]
We proceed again according to the standard method paying more attention to measurability. Consider a variation of \eqref{eqn:ode1}
\begin{align*}
\dot x(t)&=g(x(t),t,s),\\
\dot y(t)&=\frac{\pa g}{\pa x}(x(t),t,s)
\end{align*}
with the initial conditions $x(t_0)=x_0$ and $y(t_0)=\id$. The above equations satisfy the assumptions of Theorem \ref{thm:param}, hence the solution $x(t,x_0,s)$ and $y(t,x_0,s)$ is a uniform (in $t$ and $x_0$) limit of the Picard's sequence $x^n(t,x_0,s)$ and $y^n(t,x_0,s)$. We observe that $\frac{\pa x^n}{\pa x_0}(t,x_0,s)=y^n(t,x_0,s)$, hence also $\frac{\pa x}{\pa x_0}(t,x_0,s)=y(t,x_0,s)$. The derivative $y(t,x_0,s)$ satisfies the regularity conditions by Theorem \ref{thm:param}. 

Again a detailed proof (the only difference is the absence of
the parameter $s$) can be found in \cite[Thm.
2.3.2]{bressan}. 
\end{proof}

\subsection{Weak differential equations}
We will now concentrate on a very particular class of differential
equations, which will be an object of our primary interest in this
article.

Consider two functions $a(t,s)$ and $b(t,s)$ defined on a
rectangle $K=[t_0,t_1]\times[0,1]\ni(t,s)$. We will treat $a$ and
$b$ as $\R$-valued, yet all results remain valid for $\R^n$-valued
maps. Assume that $a$ and $b$ are bounded and measurable w.r.t. both variables
separately. Let $c(t,s)$ be a fixed continuous
function on $K$. We will say that the pair $(a,b)$ is a \emph{weak
solution} (\emph{W-solution})\index{weak solution}\index{W-solution|see{weak solution}} of the differential equation
\begin{equation}\label{eqn:basic}
\pa_t b(t,s)=\pa_sa(t,a)+c(t,s)b(t,s)a(t,s)
\end{equation}
if for every function $\varphi\in C^\infty_0(K)$ the following
equality holds:
\begin{equation}\label{eqn:basic_weak}
-\iint_K\Big[b(t,s)\pa_t\varphi(t,s)-
a(t,s)\pa_s\varphi(t,s)+c(t,s)b(t,s)a(t,s)\varphi(t,s)\Big]\dd
t\dd s=0.
\end{equation}
Observe that, since we assumed only measurability of $a$ and $b$,
the boundary values on $\pa K$ are, in general, not well-defined.
To guarantee this, we have to impose extra conditions on $a$ and $b$.

\begin{theorem}\label{thm:w_wt}
Let $a$ and $b$ be a bounded W-solution of \eqref{eqn:basic}.
Assume in addition that $a$ and $b$ satisfy the following regularity conditions:
\begin{equation}\label{eqn:traces}
\int_{t_0}^{t_1}\int_0^\eps|a(t,s)-a(t,0)|\frac 1\eps\dd s\dd
t\underset{\eps\to 0}{\to}0\,,\
\int_{t_0}^{t_1}\int_{1-\eps}^1|a(t,s)-a(t,1)|\frac 1\eps\dd s\dd t\underset{\eps\to 0}{\to}0\,,
\end{equation}
\begin{equation}\int_{0}^{1}\int_{t_0}^{t_0+\eps}|b(t,s)-b(t_0,s)|\frac 1\eps\dd
t\dd s\underset{\eps\to 0}{\to}0\,,\
\int_{0}^{1}\int_{t_1-\eps}^{t_1}|b(t,s)-b(t_1,s)|\frac 1\eps\dd
t\dd s\underset{\eps\to 0}{\to}0.\label{eqn:traces1}
\end{equation}
Then the traces $a(t,0)$, $a(t,1)$, $b(t_0,s)$, $b(t_1,s)$ are well defined and $(a,b)$ is a W-solution of \eqref{eqn:basic} with a well-defined trace; that is, for every $\psi\in C^\infty(K)$ we have
\begin{equation}\label{eqn:basic_wt}
\begin{split}
-\iint_K\Big[b(t,s)\pa_t\psi(t,s)- a(t,s)\pa_s\psi(t,s)+c(t,s)b(t,s)a(t,s)\psi(t,s)\Big]\dd t\dd s\\
=\int_0^1b(t_1,s)\psi(t_1,s)-b(t_0,s)\psi(t_0,s)\dd
s+\int_{t_0}^{t_1}a(t,0)\psi(t,0)-a(t,1)\psi(t,1)\dd t.
\end{split}
\end{equation}
\end{theorem}
\begin{proof}
Fix an element $\psi\in C^\infty(K)$ and choose $\eps>0$. The idea of the proof is standard: we will approximate $\psi$ by another function $\varphi\in C^\infty_0(K)$ and, using \eqref{eqn:basic_weak} for $\varphi$ and the regularity conditions, show that \eqref{eqn:basic_wt} holds with $\eps$-accuracy. 

Define a rectangle $K_\eps:=[t_0+\eps,t_1-\eps]\times[\eps,1-\eps]\subset K$, and choose a smooth ''hat function'' $\chi_{[a,b]}:[a,b]\ra\R$ which satisfies the following conditions:
\begin{align*}
&\chi_{[a,b]}(a)=0=\chi_{[a,b]}(b)=0,& \chi_{[a,b]}(t)=1\text{\ for $t\in[a+\eps,b-\eps]$}& &\text{and }&& \|D\chi_{[a,b]}\|<\frac 2\eps .
\end{align*}
Now define $\chi(t,s):=\chi_{[t_0,t_1]}(t)\cdot\chi_{[0,1]}(s)$. Obviously, $\chi\equiv 1$ on $K_\eps$, $\chi\in C^\infty_0(K)$ and  $\|D\chi_{[a,b]}\|<\frac 4\eps$. Moreover, $\pa_t\chi(t,s)=0$ for $t\in[t_0+\eps,t_1-\eps]$, and $\pa_s\chi(t,s)=0$ for $s\in[\eps,1-\eps]$.

Now define $\varphi:=\chi\cdot\psi\in C^\infty_0(K)$ and $\wt\psi:=(1-\chi)\psi$. Decomposing $\psi=\varphi+\wt\psi$, we get 
\begin{align*}
&-\iint_K\Big[b\pa_t\psi- a\pa_s\psi+cba\psi\Big]\dd t\dd s\\
&=-\iint_K\Big[b\pa_t\varphi- a\pa_s\varphi+cba\varphi\Big]\dd t\dd s-\iint_K\Big[b\pa_t\wt\psi- a\pa_s\wt\psi+cba\wt\psi\Big]\dd t\dd s
\overset{\eqref{eqn:basic_weak}}{=} \\
&=0-\iint_Kb\pa_t\wt\psi\dd t\dd s +\iint_K a\pa_s\wt\psi\dd t\dd s-\iint_K cba\wt\psi\dd t\dd s=:I_1+I_2+I_3.
\end{align*}
We will now concentrate on the tree last summands. Observe that 
\begin{align*}
|I_3|=\left|\iint_K cba\wt\psi\dd t\dd s\right|=\left|\iint_{K-K_\eps} cba\wt\psi\dd t\dd s\right|\leq\mu(K\setminus K_\eps)\|cba\wt\psi\|\leq \eps\cdot C_3,
\end{align*}
where $C_3$ is a constant depending on $\|a\|$, $\|b\|$, $\|\psi\|$, and $\|c\|$.

Now 
\begin{align*}
I_2&=\iint_Ka\pa_s\wt\psi\dd t\dd s=\iint_{K\setminus K_\eps}a\pa_s\wt\psi\dd t\dd s\\
&=\iint_{K\setminus K_\eps}a\psi\pa_s(1-\chi)\dd t\dd s+\iint_{K\setminus K_\eps}a(1-\chi)\pa_s\psi\dd t\dd s.
\end{align*}
The last summand can be estimated by $\eps\cdot C_2$ in the same way as $I_3$ (with $C_2$ depending additionally on $\|D\psi\|$).
Now, since $\pa_s\chi(t,s)=0$ for $s\in[\eps,1-\eps]$,
\begin{align*}
\iint_{K\setminus K_\eps}a\psi\pa_s(1-\chi)\dd t\dd s=-\int_{t_0}^{t_1}\int_0^\eps a\psi\pa_s\chi\dd s\dd t-\int_{t_0}^{t_1}\int_{1-\eps}^1 a\psi\pa_s\chi\dd s\dd t=:I_4+I_5.
\end{align*}
We can write $I_4$ as 
\begin{align*}
I_4=&\int_{t_0}^{t_1}\left[a(t,s)\psi(t,s)-a(t,0)\psi(t,0)\right]\pa_s\chi\dd s\dd t+\\
&+\int_{t_0}^{t_1}a(t,0)\psi(t,0)\left[\int_0^\eps \pa_s\chi\dd s\right] \dd t=I_6+\int_{t_0}^{t_1}a(t,0)\psi(t,0)\dd t. 
\end{align*}
Now 
\begin{align*}
I_6=&\int_{t_0}^{t_1}\int_0^\eps\left(a(t,s)-a(t,0)\right)\psi(t,0)\pa_s\chi\dd s\dd t+\\
 &+\int_{t_0}^{t_1}\int_0^\eps a(t,s)\left(\psi(t,s)-\psi(t,0)\right)\pa_s\chi\dd s\dd t=:I_7+I_8.
\end{align*}
Clearly,
$$|I_7|\leq\int_{t_0}^{t_1}\int_0^\eps|a(t,s)-a(t,0)|\cdot\|\psi\|\cdot\|D\chi\| \dd s\dd t\leq\int_{t_0}^{t_1}\int_0^\eps|a(t,s)-a(t,0)|\cdot\|\psi\|\frac 4\eps \dd s\dd t,$$
so by assumptions it converges to 0 as $\eps \to 0$. Finally, using $|\psi(t,s)-\psi(t,0)|\leq s\cdot\|D\psi\|$, we get
$$|I_8|\leq\int_{t_0}^{t_1}\int_0^\eps\|a\|s\|D\psi\|\cdot\|D\chi\|\dd s\dd t\leq\int_{t_0}^{t_1}\int_0^\eps\|a\|\eps\|D\psi\|\frac 4\eps\dd s\dd t\leq \eps\cdot C_8.$$
As a consequence, we get 
$$I_4\to \int_{t_0}^{t_1}a(t,0)\psi(t,0)\dd t\quad\text{as $\eps\to 0$}.$$
Analogous estimations can be done for $I_5$. As a result we get that
$$\left|I_2-\int_{t_0}^{t_1}\left[a(t,0)\psi(t,0)-a(t,1)\psi(t,1)\right]\dd t\right|\underset{\eps\to 0}{\to} 0.$$

\noindent We can repeat the above considerations for $I_1$ to prove that
$$\left|I_1+\int_{0}^{1}\left[b(t_0,s)\psi(t_0,s)-b(t_1,s)\psi(t_1,s)\right]\dd s\right|\underset{\eps\to 0}{\to} 0.$$

The estimations for $I_1$, $I_2$ and $I_3$ show that, for a fixed $\psi$, the equality \eqref{eqn:basic_wt} is satisfied with an accuracy converging to 0 as $\eps\to 0$. \end{proof}

\begin{definition}
If a W-solution $(a,b)$ of \eqref{eqn:basic} satisfies regularity conditions \eqref{eqn:traces} and \eqref{eqn:traces1} (hence \eqref{eqn:basic_wt} holds for every  $\psi\in C^\infty(K)$), we will say that the  pair $(a,b)$  is a \emph{ WT-solution} of \eqref{eqn:basic}.
\end{definition}

For WT-solutions we can formulate a uniqueness result. 

\begin{lemma}[uniqueness of WT-solutions]\label{lem:wt_uniq}\index{WT-solution!uniqueness}
Let $a:K\ra\R$ be a bounded measurable map (w.r.t. both variables
separately), and let $b_0:[0,1]\ra \R$ be any bounded measurable
map. Then there exists at most one bounded measurable map
$b:K\ra\R$ such that $(a,b)$ is a WT-solution of \eqref{eqn:basic},
and $b(t_0,s)=b_0(s)$. Moreover, the trace $b(t_1,s)$ is
 determined uniquely.

\end{lemma}
\begin{proof}
Assume that $b(t,s)$ and $\wt b(t,s)$ are two such solutions for a
fixed $a(t,s)$. The difference $\del b(t,s):=b(t,s)-\wt b(t,s)$ is
a bounded measurable map which is a WT-solution of the linear
equation
\begin{equation}\label{eqn:1}
\pa_t\del b(t,s)=c(t,s)\del b(t,s)a(t,s)
\end{equation}
such that $\del b(t_0,s)=0$. Let us define
$B(\tau,s):=0+\int_{t_0}^\tau c(t,s)\del b(t,s)a(t,s)\dd t.$
Clearly, $B(t,s)$ is ACB w.r.t. $t$ and measurable w.r.t. $s$.
Moreover, we have $\pa_t B(t,s)=c(t,s)\del
b(t,s)a(t,s)$ in the sense of Carath{\'e}odory and, since $B(t,s)$ is continuous w.r.t. $t$, also WT.
Consequently, $(B-\del b)$ satisfies $\pa_t (B-\del
b)\overset{\text{WT}}{=}0$ and, since $B(t_0,s)=\del b(t_0,s)=0$,
we have
\begin{equation}\label{eqn:2}
-\iint_K(B-\del b)(t,s)\pa_t\psi(t,s) \dd t\dd s= \int_0^1 (B-\del
b)(t_1,s)\psi(t_1,s)\dd s
\end{equation}
for every $\psi\in C^\infty(K)$. Taking $\psi(t,s)=\phi(s)$, where
$\phi\in C^\infty(I)$, we get that $\int_0^1(B-\del
b)(t_1,s)\phi(s)\dd s=0$, thus $B(t_1,s)=\del b(t_1,s)$ a.e. In
the light of this observation \eqref{eqn:2} reads as
$$\iint_K(B-\del b)(t,s)\pa_t\psi(t,s) \dd t\dd s=0,$$
for every $\psi\in C^\infty(K)$. Since $\pa_t\psi$ can be an
arbitrary smooth function, we conclude that $B(t,s)=\del b(t,s)$
a.e. Consequently, $\del b$ is a Carath{\'e}odory solution of
\eqref{eqn:1}. Now observe that
\begin{align*}
|\del b(\tau,s)|=|\int_{t_0}^\tau\pa_t\del b(t,s)\dd
t|\leq\int_{t_0}^\tau|\pa_t\del b(t,s)|\dd t=\int_{t_0}^\tau
|c\cdot a||\del b(t,s)|\dd t,
\end{align*}
which, in view of the integral Gronwall Inequality (see, for instance, \cite[App. B]{Evans}), implies $\del b=0$.
\end{proof}


\section{Homotopy of admissible paths}\label{sec:htp}

\subsection{Definition of $E$-homotopy} In this section we will define and investigate the notion of homotopy of $E$-paths,
which is crucial for this paper. This definition will be given in
two steps, as has been done while introducing admissible paths.
First, we will describe the smooth case and later generalize the
concept to measurable $E$-paths, more suitable in control theory.

\begin{definition} Let $A_0,A_1:\sT\R|_I\lra E$ be two skew-algebroid morphisms (equivalently, two smooth admissible paths $a_0,a_1:I\lra E$).
An \emph{$E$-homotopy} between $A_0$ and $A_1$ is an algebroid
morphism
$$H:\sT\R|_I\times\sT\R|_{[0,1]}\lra E$$
such that $A_0(\cdot)=H(\cdot,\theta_0)$ and
$A_1(\cdot)=H(\cdot,\theta_1)$, where $\theta_0\in\sT_0\R$ and
$\theta_1\in\sT_1\R$ are the null vectors.
\end{definition}

Note that this definition agrees with the notion of homotopy of Lie algebroid morphisms as introduced by
Kubarski \cite{Kubarski}. Equivalently, we may understand an $E$-homotopy as a pair of maps $a,b:I\times[0,1]\lra
E$, over the same base map $x:I\times[0,1]\lra M$ with $a_0(\cdot)=a(\cdot,0)$ and $a_1(\cdot)=a(\cdot,1)$,
such that
\begin{subequations}
\begin{align}
&&\label{eqn:a_adm} t\mapsto a(t,s) &\quad\text{is admissible for every $s\in[0,1]$}\,, \\
&&\label{eqn:b_adm} s\mapsto b(t,s) &\quad\text{is admissible for
every $t\in I$}\,,
\end{align}
\end{subequations}
and, moreover, the maps $a$ and $b$ satisfy a system of
differential equations given in local coordinates $(x^a,y^i)$ in
$E$ by
\begin{equation}\label{eqn:htp_smooth}
\pa_tb^i(t,s)-\pa_sa^i(t,s)=c^i_{jk}(x(t,s))b^j(t,s)a^k(t,s).
\end{equation}

Here the functions $c^i_{jk}$ are the skew-algebroid structure
functions.  Clearly, the maps $a$ and $b$ are the $H$-images of
the canonical sections $(\pa_t,\theta_s)\in\sT_t\R\times\sT_s\R$
and $(\theta_t,\pa_s)\in\sT_t\R\times\sT_s\R$, respectively, while
admissibility conditions \eqref{eqn:a_adm} and \eqref{eqn:b_adm}
and \eqref{eqn:htp_smooth} can be easily derived by means
of \eqref{eqn:alg_morph}.

The $E$-paths $b_0(s):=b(t_0,s)$ and $b_1(s):=b(t_1,s)$ will be
called \emph{initial-point} and \emph{final-point $E$-homotopies},
respectively. We will say that an $E$-homotopy has \emph{fixed
end-points} if $b_0\equiv \theta_{x(t_0)}$ and $b_1\equiv
\theta_{x(t_1)}$.

Extending the definition  of an $E$-homotopy to
measurable $E$-paths requires a little attention. Since $E$-paths also
make sense in a measurable category, we should allow $a$ and
$b$ to be bounded measurable maps (w.r.t. both variables
separately) over an ACB map $x$ and let \eqref{eqn:a_adm} and
\eqref{eqn:b_adm} hold in a measurable sense. Condition
\eqref{eqn:htp_smooth} should be replaced by a condition that
$(a,b)$ is a WT-solution of \eqref{eqn:htp_smooth}; i.e.
\begin{eqnarray}\nonumber
&\iint_{I\times[0,1]}\Big[b^i(t,s)\pa_t\psi_i(t,s)- a^i(t,s)\pa_s\psi_i(t,s)+c^i_{jk}(x(t,s))b^j(t,s)a^k(t,s)\psi_i(t,s)\Big]\dd t\dd s\\
&=\int_{[0,1]}\Big[b^i(0,s)\psi_i(0,s)-b^i(1,s)\psi_i(1,s)\Big]\dd s\label{eqn:htp_weak}\\
&\ -\int_I\Big[a^i(t,0)\psi_i(t,0)-a^i(t,1)\psi_i(t,1)\Big]\dd t\nonumber
\end{eqnarray}
holds for every family of functions $\psi_i\in
C^\infty(I\times[0,1];\R)$. Note that considering only W-solutions
of \eqref{eqn:htp_smooth} would not be enough, since otherwise the
boundary paths $a_0(t)$, $a_1(t)$, $b_0(s)$ and $b_1(s)$ would not
be well defined. The notion of the initial-point and the
final-point $E$-homotopies, as well as the $E$-homotopy with fixed
end-points, also remains valid in this new setting.

Finally, observe that any two measurable maps
$a,b:I\times[0,1]\lra E$ over the same AC base map
$x:I\times[0,1]\lra M$ define a measurable bundle map
$H:\sT\R|_I\times\sT\R|_{[0,1]}\lra E$ (that is, a measurable map
linear on fibers), where $H(\pa_t,\theta_s)=a(t,s)$ and
$H(\theta_t,\pa_s)=b(t,s)$.

\subsection{Uniqueness of $E$-homotopies}
As a direct consequence of the definition of an
$E$-homotopy and Lemma \ref{lem:wt_uniq} we get the following
result.

\begin{lemma}[uniqueness of $E$-homotopies]\index{algebroid homotopy!uniqueness}\label{lem:htp_unique} Let $a:I\times[0,1]\lra E$ be a bounded
measurable map covering $x:I\times[0,1]\lra M$ such that $t\mapsto
a(t,s)$ is admissible for every $s$. Then there exists at most
one bounded measurable map $b:I\times[0,1]\lra E$ covering $x$
such that $(a,b)$ is an $E$-homotopy with a given initial-point
$E$-homotopy $b(t_0,s)=b_0(s)$.
\end{lemma}

\subsection{The $E$-homotopy via Stokes theorem} 

We shall now give another, more geometrical, description of an
$E$-homotopy by means of a Stokes-like formula. First, we
will introduce the notion of an integral of an $E$-$k$-form, i.e.,
an element $\omega\in\Sec(\Lambda^kE^\ast)$, over a bundle
morphism $\Phi:\sT N\lra E$. We define
$$\int_{\Phi(N)}\omega:=\int_N\Phi^\ast\omega\,,$$
where the last integral is the standard integration of the
differential $k$-form $\Phi^\ast\omega$ on the manifold $N$. Now,
if $N$ is a manifold with boundary $\pa N$, we define
$$\int_{\pa\Phi(N)}\omega:=\int_{\pa N}\Phi^\ast\omega.$$
Observe that in case $\Phi:\sT N\lra\sT M$ is the tangent lift of
a diffeomorphism $\varphi:N\lra M$, the above definitions coincide
with the standard concept of differential form integration. The
morphism $\Phi$ need not be differentiable. Since, given local
coordinates $N\supset V\overset \psi\lra
V^{'}\subset\R^n\ni(y^1,\hdots,y^n)$ on $N$,
$$\int_V\Phi^*\omega=\int_{V^{'}\subset\R^n}\omega\left(\Phi(\pa_{y^1}),\hdots,\Phi(\pa_{y^n})\right)\dd y^1\cdots\dd y^n,$$
we shall require only that $\Phi$ maps smooth sections of $\sT N$
into bounded measurable sections of $E$.

Now assume that $\Phi:\sT\R|_I\times\sT\R|_{[0,1]}\lra E$ over
$\varphi:\R|_I\times\R|_{[0,1]}\lra M$ is a bundle map defined by
means of measurable maps $a(t,s)$ and $b(t,s)$ as in the
definition of an $E$-homotopy. Assume, moreover, that conditions
\eqref{eqn:a_adm} and \eqref{eqn:b_adm} are satisfied. Take any
$E$-1-form $\alpha\in\Sec(E^\ast)$; in local coordinates,
$\alpha\sim(x^a,\alpha_i(x))$. Now
\begin{align*}
\int_\Phi\dd_E\alpha&=\iint_{I\times[0,1]}\dd_E\alpha\left(\Phi(\pa_t),\Phi(\pa_s)\right)\dd
t\dd s=
\iint_{I\times[0,1]}\dd_E\alpha\left(a(t,s),b(t,s)\right)\dd t\dd s=\\
&=\iint_{I\times[0,1]}\left(\rho^a_i(x)a^i\frac{\pa\alpha_j}{\pa
x^a}b_j-\rho^a_i(x)b^i\frac{\pa\alpha_j}{\pa
x^a}a_j-\alpha_ic^i_{jk}(x)a^jb^k\right)\dd t\dd s.
\end{align*}
Having in mind that $\rho^a_i(x(t,s))a^i=\pa_tx^a(t,s)$ and
$\rho^a_i(x(t,s))b^i=\pa_sx^a(t,s)$, and defining $\wt
\alpha_i(t,s):=\alpha_i(x(t,s))$, we get
$$\int_\Phi\dd_E\alpha=\iint_{I\times[0,1]}\left(b^j\pa_t\wt\alpha_j-a^j\pa_s\wt\alpha_j-\wt\alpha_ic^i_{jk}a^jb^k\right)\dd t\dd s.$$
Similarly,
\begin{eqnarray*}\int_{\pa\Phi}\alpha &=&\int_I\left(\wt\alpha_i(t,0)a^i(t,0)-\wt\alpha_i(t,1)a^i(t,1)\right)\dd t\\&&-
\int_{[0,1]}\left(\wt\alpha_i(0,s)b^i(0,s)-\wt\alpha_i(1,s)a^i(1,s)\right)\dd s.
\end{eqnarray*}
As we see, \eqref{eqn:htp_weak} holds for all $\wt\alpha_i$
if and only if
\begin{equation}\label{eqn:stokes}
\int_\Phi\dd_E\alpha=\int_{\pa\Phi}\alpha,
\end{equation}
which can be understood as a generalized Stokes formula\index{Stokes theorem}.

\begin{remark}
In fact, \eqref{eqn:htp_weak} is more general than
\eqref{eqn:stokes} since $\wt\alpha_i$ being the pull-back of
$\alpha$ via the map $\Phi$ cannot be an arbitrary function of $t$
and $s$. We can, however, easily overcome this drawback by using
the graph of $\Phi$ in $\sT\R|_I\times\sT\R|_{[0,1]}\times E$ and
$(\sT\R|_I\times\sT\R|_{[0,1]}\times E)$-1-forms instead of $\Phi$
and $E$-1-forms.
\end{remark}\medskip

\begin{remark}
The condition \eqref{eqn:a_adm} for $a$ (and analogously
\eqref{eqn:b_adm} for $b$) can be expressed in the Stokes-like way
as well. Consider, namely, the map
$\Phi_s(\cdot):=\Phi(\cdot,\theta_s):\sT\R|_I\lra E$. The
admissibility of $a$ reads as
$$\int_{\Phi_s}\dd_Ef=\int_{\pa\Phi_s}f$$
for every $f\in C^\infty(M)$ and $s\in[0,1]$.
\end{remark}\medskip

\subsection{Lie groupoid reduction}
To  formulate an important theorem which gives a nice
interpretation of the concept of an $E$-homotopy, we need to recall
briefly the notion of a Lie groupoid and the associated Lie
algebroid. The reader unfamiliar with these concepts can consult
\cite{weinstein_silva,mackenzie} or think of a Lie group -- Lie algebra reduction as
an example.

Consider a Lie groupoid $\GG$ with the target and the source maps
$\alpha,\beta: \GG\lra M$ and the identity section
$i:M\hookrightarrow \GG$. The groupoid $\GG$ acts on itself by
left multiplication. This action preserves leaves of $\beta$
and interchanges leaves of the foliation
$\GG^\alpha=\{\alpha^{-1}(x):x\in M\}$; hence it induces an action
on the distribution tangent to that foliation,
$$\sT^\alpha \GG:=\ker \sT\alpha\subset \sT \GG.$$
In the space of sections, $\Sec(\sT^\alpha \GG)$, we may
distinguish $\chi_L(\GG)$, the class of vector fields invariant
w.r.t. this action. Clearly, $\chi_L(\GG)$ is closed w.r.t. the Lie bracket on $\GG$ and hence forms a Lie algebra.
Moreover, since the value of a left invariant field $X$ at
$g\in\GG$ is $X(g)=(L_g)_\ast X(i_{\beta(g)})$, every element of
$\chi_L(\GG)$ is uniquely determined by its value along the
identity section $i:M\hookrightarrow G$. Thus we have an
identification (via the maps $(L_{g^{-1}})_\ast$)
$$\chi_L(\GG)\approx\Sec(E),\quad \textrm{where}\ E=\sT^\alpha \GG|_{i(M)}.$$
The subbundle $E$ is a vector bundle with the base $M$. Its space
of sections $\Sec(E)$ inherits the Lie bracket from $\chi_L(\GG)$.
This bracket, together with $\rho$, the restriction of
$\sT\beta:\sT \GG\lra \sT M$ to $E$, gives a Lie algebroid
structure on $E$. We will denote it by $\AG$. We have thus
described the structure of the \emph{Lie algebroid associated with
the Lie groupoid $\GG$}. The left action of $\GG$ allows us to
identify $\sT_g^\alpha \GG:=\sT^\alpha \GG\cap \sT_g\GG$ with
$E_{\beta(g)}$ via the map $(L_{g^{-1}})_\ast$ for every $g\in
\GG$. The Lie groupoid -- Lie algebroid reduction is hence a
smooth bundle map
\begin{equation}\label{eqn:reduction}
\xymatrix{ {\sT ^\alpha\GG}\ar[d]_{\sT\beta}\ar[rr]^{R}  &&
E\ar[d]_{\rho} \\ \sT M\ar[rr]^{\id_{\sT M}} && \sT M }
\end{equation}
on each fiber given by the isomorphism
$\sT^\alpha_g\GG\overset{(L_{g^{-1}})_\ast}\lra E_{\beta(g)}$.

A basic example of a Lie groupoid is a Lie group $G$ (with
$M=\{\ast\}$ being a single point and trivial $\alpha$, $\beta$)
with its group structure. The procedure described above is the
standard Lie group -- Lie algebra reduction. Another well-known
example is the product groupoid $\GG=M\times M$ with maps
$\alpha(x,y)=x$, $\beta(x,y)=y$ and the multiplication
$(x,y)(y,z)=(x,z)$. In this case, the associated Lie algebroid is
simply the tangent bundle $\sT M$ and left-invariant fields on
$\GG$ are canonically identified with vector fields on $M$.

An example generalizing the above two is a groupoid
$\GG_P=P\times P/G$ associated with a principal bundle $G\ra
P\overset\pi\ra M$. The source and target maps are simply
$\alpha[(p,q)]=\pi(p)$ and $\beta[(p,q)]=\pi(q)$, and the
multiplication is $[(p,q)][(q,r)]=[(p,r)]$. Each leaf of
$\GG_P^\alpha$ is canonically isomorphic with $P$. The associated
bundle $E$ can be identified with the quotient $\sT P/G$, and the
Lie algebroid structure on $E$ is simply the Lie algebra of
$G$-invariant vector fields on $P$. This construction, called the
\emph{Atiyah algebroid} associated with $P$, will be described
later in more details (cf. Section \ref{sec:exmples}).

\subsection{Integration of $E$-homotopies}
Now we show that admissible paths and algebroid homotopies on an integrable algebroid $\AG$ are true paths and true homotopies on an integrating groupoid $\GG$ reduced to $\AG$ by means of the reduction map \eqref{eqn:reduction}.

\begin{theorem}[integration]\label{thm:int_htp}\index{admissible path!integration}\index{algebroid homotopy!integration}
Let $A(\GG)\ra M$ be a Lie algebroid of a Lie groupoid
$\GG$. Fix $x_0,y_0\in M$ and an element $g_0\in\alpha^{-1}(y_0)\cap\beta^{-1}(x_0)$.

There is a 1-1 correspondence between:
\begin{itemize}
\item bounded measurable
admissible paths $a:[t_0,t_1]\lra \AG$ over an ACB path $x:[t_0,t_1]\lra M$ such that $x(t_0)=x_0$, and
\item ACB paths $g:[t_0,t_1]\lra \GG_{y_0}$ such that $g(t_0)=g_0$ and  $x(t)=\beta(g(t))$. 
\end{itemize} The correspondence is given by means of the reduction map \eqref{eqn:reduction}; i.e., $a(t)=\mathcal{R}(\pa_t g(t))=\T R_{g(t)^{-1}}(\pa_t g(t))$. 

Similarly, there is a 1-1 correspondence between:
\begin{itemize}
\item bounded
measurable algebroid homotopies $a,b:[t_0,t_1]\times[0,1]\lra \AG$ over an ACB map $x:[t_0,t_1]\times[0,1]\lra M$ such that $x(t_0,0)=x_0$, and
\item ACB homotopies $h:[t_0,t_1]\times [0,1]\lra \GG_{y_0}$ (i.e., $h$ is ACB w.r.t. both variables) such that $h(t_0,0)=g_0$ and  $x(t,s)=\beta(h(t,s))$. 
\end{itemize} 
Again, the correspondence is given by means of the reduction map \eqref{eqn:reduction}; i.e., $a(t,s)=\mathcal{R}(\pa_th(t,s))=\T R_{h(t,s)^{-1}}(\pa_t h(t,s))$ and $b(t,s)=\mathcal{R}(\pa_sh(t,s))=\T R_{h(t,s)^{-1}}(\pa_s h(t,s))$.  
\end{theorem}

\begin{proof}
We will work in local coordinates $(z^i)$ on $\GG_{y_0}$, $(x^a)$
on $M$, and linear coordinates $(x^a,y^i)$ on
$\AG$. We have induced coordinates $(z^i,\dot z^j)$ on $\T\GG_{y_0}$ and
$(x^a,\dot x^b)$ on $\T M$. 

For $g\in\GG_{y_0}$, $\T R_{g^{-1}}$ maps $\T_g\GG_{y_0}=\T^\alpha_g\GG$ isomorphically into $\AG_{\beta(g)}$. In coordinates, $\T R_{g^{-1}}:(z^i,\dot z^j)\mapsto(x^a,y^i)$ can be expressed as 
\begin{align*}
&x^a=\beta^a(z),\\
&y^i=F^i_j(z)\dot z^j,
\end{align*}
where $\beta^a(z)$ and  $F^i_j(z)$ are smooth and $F^i_j(z)$ is invertible. By $f^j_i(z)$ we will denote the inverse matrix of $F^i_j(z)$. The structure functions of the algebroid $\AG$ in these coordinates satisfy
\begin{align*}
&\rho^a_i(\beta(z))F^i_j(z)\dot z^j=\frac{\pa \beta^a(z)}{\pa z^j}\dot z^j,\\
&c^i_{jk}(\beta(z))F^j_m(z)F^k_n(z)\dot z^m\dot
z^n=\left(\frac{\pa F^i_n(z)}{\pa z_m}-\frac{\pa F^i_m(z)}{\pa
z_n}\right)\dot z^m\dot z^n,
\end{align*}
since $\rho$ is the reduced $\T\beta$, and the $\AG$-bracket is the reduced Lie bracket on $\GG$. From the above we get
\begin{align*}
&\rho^a_i(\beta(z))=\frac{\pa \beta^a(z)}{\pa z^j}f^j_i(z),\\
&c^i_{jk}(\beta(z))=\left(\frac{\pa F^i_n(z)}{\pa z_m}-\frac{\pa
F^i_m(z)}{\pa z_n}\right)f^m_j(z)f^n_k(z).
\end{align*}

To prove the first part of the assertion, observe that, if $g:[t_0,t_1]\lra\GG_{y_0}$ is an ACB path over an ACB path $x:[t_0,t_1]\lra M$, then the derivative $\pa_t g(t)\in\T_{g(t)}\GG_{y_0}$ is a bounded measurable
path, and so is  $a(t)=\mathcal{R}(\pa_t g(t))=\T R_{g(t)^{-1}}(\pa_t g(t))$, since $\mathcal{R}$ is smooth. Clearly,  $\dot x(t)=\T\beta(\pa_t g(t))=\rho(a(t))$ (cf. diagram \eqref{eqn:reduction}), so $a(t)$ is
a bounded measurable $\AG$-path.

Conversely, consider a bounded measurable admissible path $a:[t_0,t_1]\lra \AG$
over an ACB path $x:[t_0,t_1]\lra M$. For every $t\in[t_0,t_1]$ and all $g$ satisfying $\beta(g)=x(t)$ we may lift $a(t)\in \AG_{x(t)}$ to a vector $A(t,g):=\T R_g(a(t))\in\T_g\GG_{y_0}$. We would like to define $g(t)$ as a solution of the differential equation in $\GG_{y_0}$
$$\pa_t g(t)=A(g,t)$$
with the initial condition $g(t_0)=g_0$. Then, clearly, $\T R_{g(t)^{-1}}(\pa_t g(t))=\T R_{g(t)^{-1}}A(g(t),t)=a(t)$ as in the assertion. The problem is that, since $A(g,t)$ is defined only on a subset of $\GG_{y_0}$ it is not clear that the solution exists, nor that it is unique.  To overcome this difficulty consider a differential equation on $\GG_{y_0}$ given in local coordinates by
\begin{equation}\label{eqn:int_a}
\dot z^i=f^i_j(z)a^j(t).
\end{equation}
It satisfies the assumptions of Theorem \ref{thm:exist}
for measurable ODEs, so it has an ACB solution $z(t)$, unique up
to the choice of the initial point. In particular, let $z(t)$ be the solution  with $z(t_0)=g_0\in\GG_{y_0}$. The base trajectory $\wt x(t)=\beta(z(t))$ satisfies 
\begin{align*}
\pa_t\wt x^a(t)&=\frac{\pa \beta^a(z(t))}{\pa z^j}\dot z^j(t)=\frac{\pa \beta^a(z(t))}{\pa z^j}f^j_i(z(t))a^i(t)=\rho^a_i(\wt x(t))a^i(t),\\
\wt x(t_0)&=\beta(z(t_0))=\beta(g_0)=x_0.
\end{align*}
On the other hand, by admissibility of $a(t)$, we have $\pa_t x^a(t)=\rho^a_i(x(t))a^i(t)$ and $x(t_0)=x_0$,
hence; clearly, $\wt x(t)=x(t)$. This, in turn, implies that $\dot z^i(t)=f^i_j(t) a^j(t)=A^i(z(t),t)$, i.e., $g(t)=z(t)$ as above is well defined and unique. 
 
Now consider a homotopy $h:[t_0,t_1]\times[0,1]=:K\lra \GG_{y_0}$ over $x:[t_0,t_1]\times[0,1]\lra M$, which is ACB w.r.t. both variables. In local coordinates it is given by $z^i(s,t)$. Repeating the argument from the previous part, we can prove that the maps $t\mapsto a(t,s):=\T R_{h(t,s)^{-1}}\left(\pa_th(t,s)\right)$ and $s\mapsto b(t,s):=\T R_{h(t,s)^{-1}}(\pa_sh(t,s))$ are bounded measurable admissible paths over $t\mapsto x(t,s)$ and $s\mapsto x(t,s)$, respectively. In local coordinates,
\begin{align*}
&a^i(t,s)=F^i_j(z(t,s))A^j(t,s),\\
&b^i(t,s)=F^i_j(z(t,s))B^j(t,s),
\end{align*}
where we denoted $A^i(t,s):=\pa_tz^i(t,s)$ and
$B^i(t,s)=\pa_sz^i(t,s)$.

Since $h$ is a homotopy, we have
$$\iint_K z^i(t,s)\pa_t\pa_s\phi_i(t,s)\dd t\dd s=\iint_Kz^i(t,s)\pa_s\pa_t\phi_i(t,s)\dd t\dd s,$$ 
 for every $\phi_i\in C^\infty(K)$. Integrating the above equality several times by parts, we get that $A^i(t,s)$ and $B^i(t,s)$ satisfy the differential equation
$$\pa_sA^i(t,s)\overset{\text{WT}}=\pa_tB^i(t,s).$$

Now calculating the WT-derivatives of $a^i(t,s)$ and $b^i(t,s)$ we get
\begin{eqnarray*} &\pa_t b^i(t,s)-\pa_s a^i(t,s)\overset {\text{WT}}= \pa_t\left(F^i_j(z(t,s))B^j(t,s)\right)-
\pa_s\left(F^i_j(z(t,s))A^j(t,s)\right)\\
&\overset {\text{WT}}=\left(\frac{\pa F^i_n}{\pa
z^m}(z(t,s))-\frac{\pa F^i_n}{\pa z^m}(z(t,s))
\right)B^n(t,s)A^m(t,s)\\&+F^i_j(z(t,s))\left(\pa_tB^j(t,s)-\pa_sA^j(t,s)\right)\\
&=c^i_{jk}(\beta(z(t,s)))b^j(t,s)a^k(t,s)+0
=c^i_{jk}(x(t,s))b^j(t,s) a^k(t,s).
\end{eqnarray*}
We see that $(a,b)$ is an $\AG$-homotopy.

Conversely, let $a,b:[t_0,t_1]\times[0,1]\lra \AG$ over $x:[t_0,t_1]\times[0,1]\lra
M$ be an algebroid homotopy. By the first part of the assertion we can uniquely integrate the admissible path $s\mapsto b_0(t_0,s)$ to an ACB path $g_0(s)\in\GG_{y_0}$ with $g_0(0)=g_0$. Next we can uniquely integrate each admissible path $t\mapsto a(t,s)$ to an ACB path $g(t,s)\in\GG_{y_0}$ such that $g(t_0,s)=g_0(s)$. In local coordinates $g(t,s)$
is a solution of the differential equation (cf. the previous part of this proof)
\begin{align*} 
\pa_t z^i(t,s)&=f^i_j(z(t,s))a^j(t,s),\\
 z^i(t_0,s)&=z_0^i(s),
\end{align*}
where $z_0(s)=g_0(s)$ is ACB. By Theorem \ref{thm:param}, $g(t,s)$ is ACB w.r.t. both variables. Now $g(t,s)$ is an ACB homotopy in $\GG_{y_0}$ hence, as has already been proved, it reduces to an algebroid homotopy $\wt a,\wt b:[t_0,t_1]\lra\AG$. By construction, $\wt a(t,s)=a(t,s)$ and $\wt b(t_0,s)=b(t_0,s)$. We see that $(a,b)$ and $(a,\wt b)$ are two WT-solutions of \eqref{eqn:htp_smooth} with the same initial-point $\AG$-homotopy $b(t_0,s)$. By Lemma \ref{lem:htp_unique} $b(t,s)= \wt b(t,s)$.
\end{proof}

\begin{corollary}\label{cor:htp_P}
Theorem \ref{thm:int_htp} establishes the equivalence between
$A(\GG)$-paths/homotopies and standard paths/homotopies in a
single $\alpha$-fibre in the groupoid $\GG$. For the
groupoid $\GG_P=P\times P/G$ and the associated Atiyah algebroid
$\T P/G$, these fibres are canonically isomorphic to $P$, so
the $E$-homotopies are just standard homotopies in $P$ reduced to
$\T P/G$. The two are equivalent up to the choice of the initial point.
\end{corollary}

\begin{remark}{\rm
Theorem \ref{thm:int_htp} is closely related to the ideas of
Crainic and Fernandes. In \cite{crainic_fernandes} they use the notion of an
$E$-homotopy to study the problem of integrability of Lie
algebroids. The crucial observation is the correspondence between
paths in the groupoid $\GG$ and admissible paths in the associated
algebroid $A(\GG)$, and between homotopies in $\GG$ and
$A(\GG)$-homotopies. If $\GG$ is simply connected, then there is a
one-to-one correspondence between points in $\GG$ and homotopy
classes of paths in $\GG$. Since the notion of an $E$-homotopy allows
us to speak of homotopies and homotopy classes without referring
to the integral object $\GG$, one can use the structure of an
algebroid on $E$ alone to reconstruct the integral object, which is
the space of $E$-paths divided by the $E$-homotopy equivalence.
This space, together with the multiplication constructed in terms
of $E$-path composition, has the structure of a simply connected
topological groupoid. In some cases (the obstructions have been
described in \cite{crainic_fernandes}) it has a natural smooth structure and
yields the Lie groupoid integrating our Lie algebroid. Note,
however, that Crainic and Fernandes use the smooth version of
Theorem \ref{thm:int_htp}, whereas we work in the measurable
category.
}\end{remark}\medskip

\subsection{$E$-homotopy classes} On the set of measurable $E$-paths we can introduce the following equivalence relation.

\begin{definition}\label{def:E_htp_class}
We will say that two measurable $E$-paths $a_0,a_1:I\lra E$ are
\emph{equivalent} (and write $a_0\sim a_1$) if there exists an
$E$-homotopy \underline{with fixed end-points} joining $a_0$ and
$a_1$. The equivalence class of an element $a$ will be denoted by
$[a]$ (or sometimes $[a(t)]_{t\in I}$) and called the
\emph{$E$-homotopy class}.
\end{definition}
\begin{remark}{\rm \label{rem:E_path_equiv} Clearly, in light of Theorem \ref{thm:int_htp}, for an integrable
algebroid $E=A(\GG)$ two $E$-paths are equivalent if their
integral paths in $\GG$ are homotopic in the standard sense.
So far, the above definition does not allow us to compare
$E$-paths defined on different time intervals. This can be done by
adding a natural condition that the composition with a null path not change the equivalence class:  $[a\circ
\theta_{x(t_1)}]=[a]=[\theta_{x(t_0)}\circ a]$.
}\end{remark}\medskip

Observe that, since admissible paths can be composed, so can
$E$-homotopies. Let, namely, $a,b:I\times[0,1]\lra E$ over $x$ and
$\ol a,\ol{b}:J\times[0,1]\lra E$ over $\ol{x}$ (where
$I=[t_0,t_1]$ and $J=[t_1,t_2]$) be two $E$-homotopies. Assume
that the final-point $E$-homotopy of the first and the
initial-point $E$-homotopy of the second coincide; i.e.
$b(t_1,s)=\ol b(t_1,s)$ a.e. (hence $x(t_1,s)=\ol x(t_1,s)$, so
$a(\cdot,s)$ and $\ol{a}(\cdot,s)$ are composable for every
$s\in[0,1]$). The maps $\wt a,\wt b:I\cup J\times[0,1]\lra E$
defined as 
$$\wt a(t,s)=\begin{cases}
a(t,s)& \text{for $t\leq t_1$},\\
\ol a(t,s)& \text{for $t>t_1$}
\end{cases}$$
and
$$\wt b(t,s)=\begin{cases}
b(t,s)& \text{for $t\leq t_1$}\\,
\ol b(t,s)& \text{for $t>t_1$}
\end{cases}$$
clearly form an $E$-homotopy joining $a(\cdot,0)\circ \ol
a(\cdot,0)$ and $a(\cdot,1)\circ \ol a(\cdot,1)$. The
initial-point $E$-homotopy is $b(t_0,\cdot)$, while the
final-point $E$-homotopy is $\ol b(t_2,\cdot)$.

Due to the above construction the multiplication on $E$-homotopy
classes given by
$$[a][\ol a]:=[a\circ \ol a]\quad \text{when $a$ and $\ol a$ are composable,}$$
is well defined.

\subsection{$E$-homotopies as families of $E$-paths} Now we shall discuss some properties of $E$-homotopies.
The following lemma emphazes the role of AL algebroids. Roughly
speaking, it turns out that for AL algebroids one-parameter
families of $E$-paths are $E$-homotopies.

\begin{lemma}[generating $E$-homotopies]\label{lem:gen_E_htp}\index{almost Lie algebroid!characterisation}
Let $E$ be an AL algebroid, and let $a:I\times[0,1]\lra E$ be a
one-parameter family of bounded measurable $E$-paths (that is,
$t\mapsto a(t,s)$ is admissible for every $s$) covering
$x:I\times[0,1]\lra M$. Assume that $a(t,s)$ is ACB w.r.t. $s$; that is, $\pa_sa(t,s)$ is defined a.e. and is bounded
measurable w.r.t. both variables. Let $b_0(s)$ be an
arbitrary bounded measurable $E$-path covering $x(t_0,s)$.

Then there exists an unique $E$-homotopy $a,b:I\times[0,1]\lra E$
such that $b(t_0,s)=b_0(s)$. Moreover, $b(t,s)$ is ACB w.r.t. $t$
(that is, $\pa_tb(t,s)$ is defined a.e. and is bounded measurable
w.r.t. both variables).
\end{lemma}

\begin{proof}
By the definition of an $E$-homotopy, $b(t,s)$ should be a map covering
$x(t,s)$ such that $s\mapsto b(t,s)$ is admissible and 
\eqref{eqn:htp_weak} is WT-satisfied. Observe that,
since $\pa_s a(t,s)$ is well defined a.e., the system of equations
\begin{equation}\label{eqn:gen_htp}
\pa_tb^i(t,s)=\pa_s a^i(t,s)+c^i_{jk}(x(t,s))b^j(t,s)a^k(t,s)
\end{equation}
for $b(t,s)$ satisfies the assumptions of Theorem \ref{thm:param}. Consequently,
 it has a unique Carath{\'e}odory solution $b^i(t,s)$ for a given
initial condition $b^i(t_0,s)= b_0^i(s)$. The solution $b(t,s)$ is
ACB w.r.t. $t$ and, since the parameter-$s$-dependence of both right-hand side
of \eqref{eqn:gen_htp} and the initial condition is bounded
measurable, so is the $s$-dependence of the solution $b(t,s)$.
Consequently, the right-hand side of \eqref{eqn:gen_htp} is bounded and
measurable w.r.t. both $t$ and $s$, and hence so is $\pa_tb(t,s)$---
the left-hand side of \eqref{eqn:gen_htp}. Clearly, $a$ and
$b$ are regular enough to satisfy the assumptions of Theorem
\ref{thm:w_wt}, so the integral condition \eqref{eqn:htp_weak}
holds.

To prove that thus constructed $(a,b)$ is indeed an
$E$-homotopy, it is enough to show that $s\mapsto b(t,s)$ is
admissible for every fixed $t$. Consider a map
$$\chi^a(t,s):=\pa_s x^a(t,s)-\rho^a_i\big(x(t,s)\big)b^i(t,s).$$
We shall show that $\chi^a= 0$ a.e. Observe that, since $a(t,s)$
is a family of admissible paths, we have
\begin{equation}\label{eqn:gen_htp_2}
\pa_tx^a(t,s)=\rho^a_k\left(x(t,s)\right)a^k(t,s) \text{ a.e.}
\end{equation}
The right-hand side of this equation is differentiable with
respect to $s$, and hence so is the left-hand side, and
$$\pa_s\pa_tx^a(t,s)=\frac{\pa\rho^a_k}{\pa x^b}\left(x(t,s)\right)\pa_s x^b(t,s)a^k(t,s)+\rho^a_k\left(x(t,s)\right)\pa_s a^k(t,s).$$

Consequently, as $\pa_t\pa_s
x^a(t,s)\overset{\text{WT}}{=}\pa_s\pa_tx^a(t,s)$ (since $x(t,s)$ is a
true homotopy in $M$),
\begin{eqnarray*}\pa_t \chi^a&=&\pa_t\pa_s
x^a-\pa_t\left(\rho^a_k b^k\right)\overset {\text{WT}}=\pa_s\pa_t
x^a-\pa_t\left(\rho^a_k b^k\right)\\&=&\frac{\pa \rho^a_k}{\pa x^b}\pa_s
x^b a^k+\rho^a_k\pa_s a^k-\frac{\pa \rho^a_i}{\pa x^b}\pa_t x^b
b^i- \rho^a_i\pa_t b^i\,,
\end{eqnarray*}
which, in view of \eqref{eqn:gen_htp_2}
and \eqref{eqn:gen_htp}, equals
\begin{align*}
&\frac{\pa\rho^a_k}{\pa x^b}\pa_s x^b a^k+\rho^a_k\pa_s
a^k-\frac{\pa \rho^a_i}{\pa x^b}\rho^b_k a^k b^i-
\rho^a_i\left(\pa_s a^i+c^i_{jk}b^ja^k\right)=\\
&=\left(\frac{\pa}{\pa
x^b}\rho^a_k\right)a^k\chi^b+\left[\left(\frac{\pa}{\pa
x^b}\rho^a_k\right)\rho^b_j-\left(\frac{\pa}{\pa
x^b}\rho^a_j\right)\rho^b_k-\rho^a_ic^i_{jk}\right]b^ja^k.
\end{align*}
Since $E$ is an AL algebroid, the last term vanishes and we have
$$\pa_t\chi^a\overset{\text{WT}}=\left(\frac{\pa}{\pa x^b}\rho^a_k\right)a^k\chi^b.$$
Thus $\chi^a$ is a WT-solution of a linear differential equation
with measurable r.h.s. and the initial condition $\chi^a(t_0,s)=0$
(since $b_0(s)$ is admissible). Repeating the argument from the proof of Lemma
\ref{lem:wt_uniq} we conclude that $\chi^a=0$ a.e.
\end{proof}

It turned out that, when a skew-algebroid is almost Lie,
$E$-homotopies are the true homotopies in the space of $E$-paths
(i.e. one-parameter families of $E$-paths). This has already been
observed in \cite[Thm. 3]{GG_var_calc} in a slightly different form.

\subsection{Further properties of $E$-homotopies}

\begin{lemma}\label{lem:reparam1}\index{admissible path!reparametrisation}
Let $a:[t_0,t_1]\lra E$ be a bounded measurable $E$-path over $x(t)$, and let $h:[0,1]\ra[t_0,t_1]$ be an invertible $C^1$-function. Define
\begin{equation}\label{eqn:reparam}
\begin{split}
a(t,s)&:=\frac{h(s)-t_0}{t_1-t_0} a\left(t_0+\frac{t-t_0}{t_1-t_0}(h(s)-t_0)\right),\\
b(t,s)&:=\frac{t-t_0}{t_1-t_0}\dot h(s) a\left(t_0+\frac{t-t_0}{t_1-t_0}(h(s)-t_0)\right).
\end{split}
\end{equation}
Then the pair $(a,b)$ is an $E$-homotopy over $x(t,s):=x\left(t_0+\frac{t-t_0}{t_1-t_0}(h(s)-t_0)\right)$.
\end{lemma}

\begin{proof}
For notation simplicity assume that $[t_0,t_1]=[0,1]$. Then
\begin{align*}
a(t,s)&=h(s)a(t h(s))\quad \text{and}\\
b(t,s)&=t\dot h(s) a(t h(s)).
\end{align*}
First, note that $t\mapsto a(t,s)$ and $s\mapsto b(t,s)$
are admissible. Indeed, from $\dot{x}(t)=\rho\left(a(t)\right)$ we deduce that 
$$\pa_tx(t, s)=h(s)\dot{x}(th(s))=h(s)\rho\left(a(th(s))\right)=\rho\left(
a(t,s)\right).$$
Similarly, $\pa_sx(t,s)=\rho\left(b(t,s)\right)$. 
Now we will check that
$$\pa_t[t\dot h(s)a^i(th(s))]\overset {\text{W}}=\pa_s[h(s)a^i(th(s))]+c^i_{jk}(x(t,s))[t\dot h(s)a^j(th(s))]\cdot[h(s)sa^k(th(s))].$$
By the skew-symmetry of $c^i_{jk}$, the last term vanishes, so we
have to check if
$$\pa_t[t\dot h(s)a^i(th(s))]\overset {\text{W}}=\pa_s[h(s)a^i(th(s))].$$
The latter is certainly true, as both sides are equal 
$\dot h(s)a^i(th(s))+th(s)\dot h(s) A^i(th(s))$, where $A^i(t)$ is the distributive derivative of $a^i(t)$. 

To finish the proof we shall show that
$a(t,s)$ and $ b(t,s)$ satisfy the regularity conditions \eqref{eqn:traces} and \eqref{eqn:traces1}.
 
Assume that $h(0)=t_0=0$ and $h(1)=t_1=1$ (the case $h(1)=t_0=0$ and $h(0)=t_1=1$ is completely analogous). Now
\begin{align*}\int_0^1\int_0^\eps| a(t,s)- a(t,0)|\frac 1\eps\dd s\dd t=\frac 1\eps\int_0^\eps\int_0^1|h(s) a(th(s))| \dd t \dd s
\leq \frac 1\eps\int_0^\eps h(s)\|a\|\dd s.
\end{align*}
The later converges to 0 as $\eps\to 0$, since $h(s)\overset{s\to 0}\lra h(0)=0$.
Next, 
\begin{align*}\int_0^1\int_{1-\eps}^1| a(t,s)- a(t,1)|\frac 1\eps \dd s\dd t=\frac 1\eps\int_{1-\eps}^1\int_0^1|h(s) a(th(s))-h(1)a(th(1))|\dd t\dd s\\
\leq\frac 1\eps\int_{1-\eps}^1 h(s)\int_0^1|a(th(s))-a(th(1))|\dd t\dd s+\frac 1\eps\int_{1-\eps}^1|h(s)-h(1)|\int_0^1 |a(t h(s))|\dd t \dd s.
\end{align*}
The second factor converges to 0 as $\eps\to 0$ because $h(s)$ is continuous at $s=1$ and $a$ is bounded. By Lemma \ref{lem:reg} the measurable function $g(s):=\int_0^1|a(th(s))-a(th(1))|\dd t$ is regular at $s=0$ and, moreover, $g(s)=0$. We conclude that
$$\frac 1\eps\int_{1-\eps}^1h(s)g(s)\dd s\to 0 \quad\text{as}\quad \eps\to 0.$$
Consequently, $\int_0^1\int_{1-\eps}^1| a(t,s)- a(t,1)|\frac 1\eps \dd s\dd t\overset{\eps\to 0}\lra 0$ and conditions \eqref{eqn:traces} are fulfilled.

Now check \eqref{eqn:traces1}. The first of the two conditions is a matter of a simple estimation:
\begin{align*}
&\int_0^1\int_0^\eps|b(t,s)-b(0,s)|\frac 1\eps \dd t\dd s=\frac 1\eps\int_0^\eps\int_0^1|t\dot h(s) a(th(s))|\dd s\dd t\leq \frac 1\eps\int_0^\eps t\|\dot h\|\cdot\|a\|\dd t\overset{\eps\to 0}\lra 0.
\intertext{For the second we can estimate:}
&\int_0^1\int_{1-\eps}^1|b(t,s)-b(1,s)|\frac 1\eps\dd t\dd s=\frac 1\eps\int_{1-\eps}^1\int_0^1|t\dot h(s)a(t h(s))-\dot h(s)a(h(s))|\dd s\dd t\\
&\leq\frac 1\eps\int_{1-\eps}^1 t\|\dot h\|\int_0^1|a(t h(s))-a(h(s))|\dd s\dd t+\frac 1\eps \int_{1-\eps}^1|1-t|\cdot\|\dot h\|\int_0^1|a(h(s))|\dd s\dd t. 
\end{align*}
The last factor clearly converges to 0 as $\eps\to 0$. Using Lemma \ref{lem:reg} we show that the measurable function $k(t):=\int_0^1|a(t h(s))-a(h(s))|\dd s$ is regular at $t=1$ and, moreover, $k(1)=0$. We conclude that
$$\frac 1\eps\int_{1-\eps}^1t\|\dot h\|k(t)\dd t\overset{\eps\to 0}\lra 0,$$ 
which proves that conditions \eqref{eqn:traces1} are satisfied. By Theorem \ref{thm:w_wt} the pair $(a,b)$ is a WT-solution of \eqref{eqn:htp_smooth}, and hence $E$-homotopy. 
 \end{proof}

\noindent As a corollary we obtain the following fact.
\begin{lemma}[shrinking an $E$-path]\label{lem:van_curv}
Let $a:[0,1]\lra E$ be a measurable $E$-path over $x(t)$. Define $a(t,s):=s a(ts)$ and $b(t,s):=t
a(ts)$ for $t,s\in [0,1]$. The pair $(a,b)$ is
an $E$-homotopy over $ x(t,s)=x(ts)$. Its initial-point
$E$-homotopy is $\theta_{x(0)}$, and the final-point $E$-homotopy
is $a$.

Similarly, consider $\wt{a}(t,s):=(1-s)a(1-(1-t)(1-s))$ and
$\wt{b}(t,s):=(1-t)a(1-(1-t)(1-s))$ where $t,s\in[0,1]$. The pair $(\wt{a},\wt{b})$ is an $E$-homotopy over $\wt x(t,s):=x(1-(1-t)(1-s))$. Its
initial-point $E$-homotopy is $a$, and the final-point $E$-homotopy
is $\theta_{x(1)}$.
\end{lemma}

\begin{proof}  The assertion follows from Lemma \ref{lem:reparam1}. For $(a,b)$ we simply take $a(t)$ and $h(s)=s$. 

\noindent For $(\wt a,\wt b)$ we use Lemma \ref{lem:reparam1} with $\wt a(t):=a(1-t)$ defined on an interval $[t_0=1,t_1=0]$ and $h(s)=s$. 
\end{proof}

\noindent We can now state the following important result.
\begin{lemma}\label{lem:htp}\index{algebroid homotopy!class}
Let $a,b:[t_0,t_1]\times[0,1]\lra E$ be an $E$-homotopy covering
$x:[t_0,t_1]\times[0,1]\lra M$. Then we have the following equality of $E$-homotopy classes:
$$[a(t,0)]_{t\in[t_0,t_1]} [b(t_1,s)]_{s\in[0,1]}=[b(t_0,s)]_{s\in[0,1]} [a(t,1)]_{t\in[t_0,t_1]}.$$
\end{lemma}

\begin{proof} The first part of Lemma \ref{lem:van_curv}, applied to the curve $s\mapsto b(t_0,s)$,
gives us the existence of $E$-homotopy $c,d:[0,1]\times[0,1]\lra
E$ such that $c(t,0)=\theta_{x(t_0)}$, $c(s,1)=b(t_0,s)$,
$d(0,s)=\theta_{x(t_0)}$, and $d(1,s)=b(t_0,s)$. Similarly, using
the second part of Lemma \ref{lem:van_curv} for $s\mapsto
b(t_1,s)$, we obtain $E$-homotopy $e,f:[0,1]\times[0,1]\lra E$
such that $e(s,0)=b(t_1,s)$, $e(t,1)=\theta_{x(t_1)}$,
$f(0,s)=b(t_1,s)$, and $f(1,s)=\theta_{x(t_1)}$.

Clearly, $E$-homotopies $(c,d)$, $(a,b)$, and $(e,f)$ are
composable and their composition is an $E$-homotopy with fixed
end-points which establishes an equivalence of $E$-paths
$\theta_{x(t_0)}\circ a(\cdot,0)\circ b(t_1,\cdot)$ and
$b(t_0,\cdot)\circ a(\cdot,1)\circ\theta_{x(t_1)}$.
\end{proof}

\begin{remark} \label{rem:lem_htp}
The above lemma is very important, as it shows the relation
between $E$-homotopies with and without fixed end-points. If
$a,b:I\times[0,1]\lra E$ is an $E$-homotopy joining $a_0$ and
$a_1$, then the composition of $a_0$ with the final-point
$E$-homotopy $b(t_1,\cdot)$ is equivalent to the composition of
the initial-point $E$-homotopy $b(t_0,\cdot)$ with $a_1$. Thus, if
the initial-point $E$-homotopy $b(t_0,\cdot)$ vanishes, in order
to check whether $[a_0]=[a_1]$, it is enough to check whether
$[b(t_1,\cdot)]=0$. Thus the problem of equivalence of $a_0$ and
$a_1$ can be solved by investigating the final-point $E$-homotopy. 
\end{remark}\medskip

\noindent Finally, as a corollary from Lemmas \ref{lem:reparam1} and \ref{lem:htp} we obtain a  result about reparametrisation of
$E$-paths.

\begin{lemma}[reparametrization]\label{lem:reparam}\index{admissible path!reparametrisation}
Let $a:[t_0,t_1]\lra E$ be a measurable $E$-path, and let $h:[0,1]\ra[t_0,t_1]$ be an invertible $C^1$-function. Define $\wt a(t):=\dot h(t)a(h(t))$ for $t\in[0,1]$. Then
\begin{align}
&[\wt a(t)]_{t\in[0,1]}=[a(t)]_{t\in[t_0,t_1]} \quad\text{if $h(0)=t_0$ and $h(1)=t_1$,}\label{eqn:reparam1}\\
&[a(t)]_{t\in[t_0,t_1]}\cdot[\wt a(t)]_{t\in[0,1]}=0 \quad\text{if $h(0)=t_1$ and $h(1)=t_0$.}\label{eqn:reparam2} 
\end{align}
\end{lemma}

\begin{proof} Consider a homotopy \eqref{eqn:reparam} from Lemma \ref{lem:reparam1}. If $h(0)=t_0$ and $h(1)=t_1$ we have $a(t,0)=0$, a(t,1)=a(t)$, b(t_0,s)=0$, and $b(t_1,s)=\dot h(s)a(h(s))=\wt a(s)$. By Lemma \ref{lem:htp}, 
$$[\theta_{x(t_0)}]\cdot[\wt a(s)]_{s\in[0,1]}=[\theta_{x(t_0)}]\cdot[a(t)]_{t\in [t_0,t_1]}.$$

Analogously, for $h(0)=t_1$ and $h(1)=t_0$ we have  $a(t,0)=a(t)$, a(t,1)=0$, b(t_0,s)=0$, and $b(t_1,s)=\dot h(s)a(h(s))=\wt a(s)$. From Lemma \ref{lem:htp} we deduce that 
$$[a(t)]_{t\in[t_0,t_1]}\cdot[\wt a(s)]_{s\in[0,1]}=[\theta_{x(t_0)}]\cdot[\theta_{x(t_0)}]=0.$$
\end{proof}

\section{Optimal control problem on AL algebroids}\label{sec:ocp}

From this section on we will restrict ourselves to AL algebroids
only. This choice is justified by the properties of $E$-homotopies
on AL algebroids discussed in Lemma  \ref{lem:gen_E_htp}.

\subsection{Control systems on algebroids} A classical control system is a differential equation on a manifold $M$
depending on a control parameter $u$ taking values in some
topological space $U$. This can be realized by a continuous map
$f:M\times U\lra \sT M$ such that $f(\cdot, u)$ is a section of
$\sT M$ for each $u\in U$. We propose the following generalization
of this notion.

\begin{definition}\label{def:con_sys}
A \emph{control system}\index{control system on AL algebroid} on an AL algebroid $E$ is a continuous map
\begin{equation}\label{eqn:def_con_sys}
f:M\times U\lra E
\end{equation}
such that, for every $u\in U$, the map $f(\cdot,u):M\lra E$ is a
$C^1$-section of $E$.
We will assume that $U$ is a subset of some Euclidean space $\R^r$. Moreover, we demand that the maps $f:M\times U\lra E$ and $\T_xf:\T M\times U\lra\T E$ are continuous. In local coordinates, if $f\sim (f^i(x,u), x^a)$, this means that $f^i(x,u)$ is continuous w.r.t. $x$ and $u$, differentiable w.r.t. $x$, and that $\frac{\pa f^i}{\pa x_a}(x,u)$ is continuous w.r.t. $x$ and $u$. 
\end{definition}

Observe that, for the tangent algebroid $E=\T M\ra M$, the above definition coincides with the classical one.
On the other hand, one easily sees (cf. Theorem \ref{thm:int_htp})
that a right-invariant control system on a Lie groupoid $\GG$
reduces to a system of the above form on the associated Lie
algebroid $A(\GG)$. For example, a right-invariant control system
on a gauge groupoid $\GG_P=P\times P/G$ of a principal bundle $G\ra
P\ra M$ is determined by its values on a single leaf of
$\GG^\alpha_P$ canonically isomorphic to $P$. Consequently, it is
equivalent to a $G$-invariant control system on $P$ and reduces to
a control system of the form \eqref{eqn:def_con_sys} on the Atiyah
algebroid $E=\sT P/G$. In particular, for a right-invariant system
on a Lie group $G$, Definition \ref{def:con_sys} coincides
with the reduced control system on its Lie algebra $\g$ as
described in \cite[Ch. 12]{jurdjevic}.

Now, for a given function $u:I\lra U$ (\emph{control}\index{control}), the map
\eqref{eqn:def_con_sys} defines a first-order ODE on $M$,
\begin{equation}\label{eqn:con_sys}
\dot{x}(t)=\rho\left(f(x(t),u(t))\right).
\end{equation}
We will restrict our attention only to functions $u$ of a certain
class (called \emph{admissible controls}\index{admissible controls|main}). In this paper these
are controls which are bounded and measurable, but one can
think of smaller classes: piecewise continuous or piecewise
constant functions. The set of all admissible controls will be
denoted by $\Uadm$.

Clearly, if $u(\cdot)$ is admissible, the map
$g(x,t)=\rho\left(f(x,u(t))\right)$ is differentiable w.r.t. $x$
and measurable w.r.t. $t$, so the assumptions of Theorem
\ref{thm:exist} hold. Consequently, we have the results of local
existence and uniqueness for the solutions of \eqref{eqn:con_sys}.
Observe that if $x(\cdot)$ is a solution of \eqref{eqn:con_sys}
for $u(\cdot)\in\Uadm$, then the path $f(x(\cdot),u(\cdot)):I\lra
E$ is a measurable $E$-path over $x(\cdot)$. This path will be
called a \emph{trajectory} of \eqref{eqn:con_sys}, whereas for the
pair $(x(\cdot), u(\cdot))$ we will use the term \emph{controlled
pair}.

\subsection{Optimal control problem on algebroids}
We introduce now a \emph{cost function}\index{cost function} $L:M\times U\lra \R$. We will assume the same regularity conditions for $L$ as in the case of
$f$, namely, that $L$ is a continuous function on $M\times U$,
which is of class $C^1$ w.r.t. the first variable and that the derivative $\dd_x L:\T M\times U\lra\R$ is continuous. If now
$(x(t), u(t))$, with $t\in I$, is a controlled pair for
\eqref{eqn:con_sys}, we define the \emph{total cost}\index{total cost|main} of this pair
to be $\int_{t_0}^{t_1}L\big(x(t),u(t)\big)\dd t.$ Note that, since
$L$ is continuous, $u(t)$ is bounded measurable, and the interval
$[t_0,t_1]$ is compact, the above integral is finite whenever the
solution $x(t)$ exists. Now we can define
optimal control problem for the data introduced above. This
definition may seem unnatural at first sight, yet we will
motivate it in a moment.

\begin{definition}\label{def:ocp} For a control system \eqref{eqn:def_con_sys} and a cost function $L$ we define
an \emph{OCP} as follows:
\begin{equation}\tag{P}\label{eqn:P}
\begin{split}
\text{minimize} \int_{t_0}^{t_1}L\big(x(t),u(t)\big)\dd t
\text{ over all controlled pairs $(x,u)$ of \eqref{eqn:con_sys} such that}\\
\text{ the $E$-homotopy class of the trajectory $f(x(t),u(t))$
equals $[\sigma]$\,,}
\end{split}
\end{equation}
where $\sigma$ is a fixed $E$-path. The interval $[t_0,t_1]$ is to
be determined as well.
\end{definition}

\medskip
The second condition in \eqref{eqn:P} plays the role of the
fixed-end-point boundary condition in standard control theory.
Indeed, let us compare two OCPs in the case $E=\sT M$:
\begin{equation}\label{eqn:P1}\tag{$\textrm{P}_1$}
\begin{split}
\text{minimize } \int_{t_0}^{t_1}L\big(x(t),u(t)\big)\dd t \text{ over all controlled pairs $(x,u)$}\\
\text{satisfying\ } x(t_0)=x_0,\ x(t_1)=x_1,
\end{split}
\end{equation}
and
\begin{equation}\label{eqn:P2}\tag{$\textrm{P}_2$}
\begin{split}
\text{minimize } \int_{t_0}^{t_1}L\big(x(t),u(t)\big)\dd t \text{ over all controlled pairs $(x,u)$}\\
\text{for which the homotopy class of $x(t)$ equals $[\sigma]$},
\end{split}
\end{equation}
where $\sigma$ is a fixed path in $M$ joining $x_0$ and $x_1$.

\medskip
First, note that every solution of \eqref{eqn:P1} gives a
solution of \eqref{eqn:P2} for some $[\sigma]\in\Pi_1(M,x_0,x_1)$.
On the other hand, a controlled pair which has the minimal total
cost among the solutions of \eqref{eqn:P2} for all
$[\sigma]\in\Pi_1(M,x_0,x_1)$ is a solution of \eqref{eqn:P1}.

Clearly, if $M$ is simply connected, then every homotopy class of
a path in $M$ is uniquely determined by its end-points, so
\eqref{eqn:P1} and \eqref{eqn:P2} are equivalent. If this is not
the case, we can always lift the control system, the path
$\sigma$, and the cost function to $\wt M$,the universal cover
of $M$. Now \eqref{eqn:P2} on $M$ is equivalent to \eqref{eqn:P1}
on $\wt M$ for the end-points $\wt \sigma(t_0)$ and $\wt
\sigma(t_1)$ (here $\wt \sigma$ is any lift of the path $\sigma$).
Finally, observe that PMP gives only necessary conditions for
\emph{local} optimality; hence the controlled pairs pointed
by the PMP will be candidates for solution of both \eqref{eqn:P1} and
\eqref{eqn:P2}. Concluding, for standard control systems on
manifolds, OCPs \eqref{eqn:P1} and \eqref{eqn:P2} are almost
equivalent and, what is more, they are not distinguished by the PMP.
Now we can finally motivate Definition \ref{def:ocp}.

\begin{remark}{\rm \label{rem:ocp_red}
Since, by Theorem \ref{thm:int_htp}, every homotopy on a Lie
groupoid $\GG$ reduces to an $E$-homotopy on a Lie algebroid
$E=A(\GG)$, it is clear that a left-invariant OCP of the form
\eqref{eqn:P2} on a Lie groupoid $\GG$ reduces to an OCP
\eqref{eqn:P} on $E=A(\GG)$. In light of this observation, the
OCP from Definition \ref{def:ocp} may be regarded as a general
reduction scheme for left-invariant systems on Lie groupoids. The
sense of using the $E$-homotopy instead of the fixed-end-point
condition is to express all the data in terms of reduced objects,
without referring to the integral ones. What is more, the OCP
\eqref{eqn:P} can be also regarded as a problem on its own (not
just reduced OCP), since not every AL (and even Lie) algebroid
comes from some reduction.
}\end{remark}

\subsection{OCP in terms of the product algebroid $\bm{E\times\sT\R}$}
For a control system on an AL algebroid $E$, similar to the classical
situation of the tangent algebroid $E=\T M\ra M$, there is an elegant formulation of the OCP \eqref{eqn:P}  in terms of the product algebroid
$E\times\T\R$. The idea is to incorporate the cost function into
the control system \eqref{eqn:con_sys}.

Denote by $\bm{A}$ the product algebroid structure on
$\bm\tau=(\tau_E,\tau_{\T\R}): E\times\T\R\lra M\times\R$ (see Section \ref{sec:ala}). We
will consequently use bold letters to emphasise objects associated
with $\bm{A}$, whereas objects associated with the $\T\R$-component of
$\bm A$ will be distinguished by underlining. For example,
$\bm{x}=(x,\ul x)\in M\times\R$ and $\bm{a}=(a,\ul a)\in
E\times\T\R=\bm A$.

Introduce now a new variable $\ul x\in\R$ and, for a given admissible control
$u\in\Uadm$, consider the following extension of the differential
equation \eqref{eqn:con_sys}:
\begin{equation}\label{eqn:con_sys_A}
\left\{ \begin{aligned}
\dot{x}(t)&=\rho\left(f(x(t),u(t))\right),\\
\dot{\ul x}(t)&=L\left(x(t),u(t)\right).
\end{aligned}\right.
\end{equation}
Clearly, $\ul x(t_1)-\ul x(t_0)=\int_{t_0}^{t_1}L(x(t),u(t))\dd t$
is the total cost of the controlled pair $(x(t),u(t))$ of
\eqref{eqn:con_sys}. Equation \eqref{eqn:con_sys_A} is a
differential equation associated with the following control system
on $\bm A$:
\begin{equation}
\label{def:con_sys_A} \bm{f}=(\wt f,\ul f):(M\times\R)\times U\lra
E\times\T\R=\bm{A},
\end{equation}
where $\ul f\left((x,\ul x),u\right):=\left(\ul
x,L(x,u)\right)\in\R\times\R\approx\T\R$ and $\wt{f}\left((x,\ul
x),u\right):=f(x,u)\in E$. For a given $u\in\Uadm$ the base
trajectory $\bm x(t)=(x(t),\ul x(t))$ of \eqref{eqn:con_sys_A}
contains information on both the base trajectory $x(t)$ of
\eqref{eqn:con_sys} (for the same control $u$) and the total cost
of the controlled pair $(x(t),u(t))$. Observe that the trajectory
$\bm f(\bm x(t),u(t))$ of \eqref{eqn:con_sys_A} projects onto the
trajectory $f(x(t),u(t))$ of \eqref{eqn:con_sys} under the
canonical algebroid projection $p_E:\bm A=E\times\sT\R\lra E$.
Now the OCP \eqref{eqn:P} can be reformulated in terms of control
system \eqref{def:con_sys_A} as follows:
\begin{align}\label{eqn:P_A}\tag{\textbf{P}}
\begin{split}
&\text{minimise } \ul x(t_1)
\text{ over controlled pairs }\ (\bm x(t),u(t))=\left((x(t),\ul x(t)),u(t)\right)\\
& \text{ of \eqref{eqn:con_sys_A} satisfying the following:}\\
&1.\text{ the $E$-projection $f(x(t),u(t))$ of the trajectory $\bm f(\bm x(t),u(t))$}\\
 &\text{belongs to a fixed $E$-homotopy class $[\sigma]$;}\\
&2.\text{ }\ul x(t_0)=0.
\end{split}
\end{align}

\subsection{The $E$-homotopy associated with a control system}

As has been observed in \cite{crainic_fernandes}, algebroid homotopies can be generated
by time-dependent algebroid sections. Since the control system
\eqref{eqn:con_sys} is a family of $E$-sections $f(\cdot,u)$,
fixing an admissible control $u(t)$ gives a time-dependent section
$f(\cdot,u(t))$.  The associated $E$-homotopy can be well
understood in terms of Lemma \ref{lem:gen_E_htp}.

Solving \eqref{eqn:con_sys} for a one-parameter family of
initial conditions $x(t_0,s)=x_0(s)$ produces a one-parameter
family of base paths $x(t,s)$. It follows from Theorem \ref{thm:exist} that, if the solution $x(t,0)$ is defined on $I=[t_0,t_1]$, then so is $x(t,s)$
at least for $x_0(s)$'s close enough to $x_0(0)$. With $x(t,s)$
we can associate a one-parameter family of trajectories
$$ a(t,s):=f\big(x(t,s),u(t)\big).$$
One easily sees that for \eqref{eqn:con_sys} the
assumptions of Theorem \ref{thm:param_dif} are satisfied.
Consequently, the base trajectories $x(t,x_0)$ are continuous
differentiable w.r.t. the initial condition $x_0$ (and ACB in
$t$). As $x(t,s)=x(t,x_0(s))$ if $s\mapsto x_0(s)$ is an ACB map,
we deduce that $x(t,s)$ is ACB w.r.t. the second variable; that
is, $\pa_s x(t,s)$ is a well-defined measurable function of both
variables. Consequently, the derivative $\pa_sa^i(t,s)=\frac{\pa
f^i}{\pa x^a}(x(t,s),u(t))\pa_sx^a(t,s)$ satisfies the assumptions
of Lemma \ref{lem:gen_E_htp}. Thus, the conclusions of Lemma
\ref{lem:gen_E_htp} hold; namely, for a given bounded measurable
$E$-path $b_0(s)$ covering $x_0(s)$, there exists a measurable map
(AC w.r.t. the first variable) $b:I\times[0,1]\lra E$ with
$b(t_0,s)=b_0(s)$ such that $(a,b)$ is an $E$-homotopy. The
$t$-evolution of $b(t,s)$ is given by \eqref{eqn:gen_htp}. Observe that, since
$\pa_sx(t,s)=\rho\left(b(t,s)\right)$, we have
$\pa_sa^i(t,s)=\frac{\pa f^i}{\pa
x^a}\left(x(t,s),u(t)\right)\rho^a_k\left(x(t,s)\right)b^k(t,s)$.
Consequently, $b(t,s)\sim\left(x^a(t,s),b^i(t,s)\right)$ is a
solution of the following differential equation
\begin{equation}\label{eqn:par_trans}
\left\{
\begin{aligned} \pa_tb^i(t,s)&=\frac{\pa
f^i}{\pa x^a}\Big(x(t,s),u(t)\Big)
\rho^a_k\big(x(t,s)\big)b^k(t,s)\\&+c^i_{jk}\big(x(t,s)\big)b^j(t,s)f^k\big(x(t,s),u(t)\big),\\
\pa_tx^a(t,s)&=\rho^a_i\big(x(t,s)\big)f^i\big(x(t,s),u(t)\big),
\end{aligned}\right.\end{equation}
with the initial conditions $b^i(t_0,s)=b_0^i(s)$ and
$x^a(t_0,s)=x^a_0(s)$.

 The above differential equation is well understood  in terms of the tools introduced in Section \ref{sec:ala}.
 For every $u\in U$, the section $f_u(\cdot):=f(\cdot,u):M\lra E$ gives rise to a linear vector field $\dd_\sT(f_u)$ on $E$.
 Evaluating it on $u(t)$ gives a time-dependent family of vector fields $\dd_\sT(f_{u(t)})$. Equation \eqref{eqn:par_trans}
 is simply the evolution along this family, $\pa_tb(t,s)=\dd_\sT\left(f_{u(t)}\right)(b(t,s))$.
On the other hand, with a time-dependent family of section
$f_u(t)$ we may associate the family of linear functions
$h_t(x,\xi):=\< f(x,u(t)),\xi>_{\tau}$ on $E^\ast$, and the
corresponding family of Hamiltonian vector fields $\X_{h_t}$. In
local coordinates,
$$\X_{h_t}(x,\xi)=\rho^b_j(x)f^j(x,u(t))\pa_{x^b}+\left(c^k_{ij}(x)f^i(x,u(t))\xi_k-
\rho^a_j(x)\frac{\pa f^i}{\pa
x^a}(x,u(t))\xi_i\right)\pa_{\xi_j}.$$ As we have seen in Section
\ref{sec:ala} (equations \eqref{eqn:ham_vf}--\eqref{eqn:hvf_tgl}),
the fields $\dd_\sT(f_{u(t)})$ and $\X_{h_t}$ give the same base
evolution (given by \eqref{eqn:con_sys}), and are related by
$\<\dd_\sT(f_{u(t)}),\X_{h_t}>_{\sT\tau}=0$.

\begin{definition}\label{def:op_par_tr}
The flows of the fields  $\dd_\sT(f_{u(t)})$ and $\X_{h_t}$ (for a
given $u\in\Uadm$) will be called \emph{operators of parallel
transport}\index{parallel transport} (in $E$ and $E^\ast$ respectively) along the solution $x(t)$ of the system \eqref{eqn:def_con_sys}. We will denote them
with $B_{tt_0}$ and $B^\ast_{tt_0}$, respectively. Analogously we
define operators  $\bm{B}_{tt_0}$ and ${\bm{B^\ast}}_{tt_0}$ for
the control system \eqref{def:con_sys_A}. Note that, by
construction, $B_{tt^{'}}\circ B_{t^{'}t_0}=B_{tt_0}$ and
$B^\ast_{tt^{'}}\circ B^\ast_{t^{'}t_0}=B^\ast_{tt_0}$.
\end{definition}

\begin{remark} \label{rem:B_htp}
Let us see that, by construction, the map
$b(t,s)=B_{tt_0}\left(b_0(s)\right)$ together with $a(t,s)$ forms
an $E$-homotopy. Moreover, $B_{tt_0}(b_0)$ is continuous w.r.t. $b_0$, $t$, and $t_0$. Indeed, $b(t)=B_{tt_0}(b_0)$ is
the solution of \eqref{eqn:par_trans} for $s=0$. The right-hand side is
measurable in $t$ and locally Lipschitz (linear) in $b$, so, by
Theorem \ref{thm:param}, $b(t)$ is AC w.r.t. $t$ and
continuous w.r.t. the initial condition $b_0$.
\end{remark}


\begin{remark} \label{rem:B_paring}
Note also that the operators $B$ and $B^\ast$ have the property of
preserving the parring $\<\cdot,\cdot>_{\tau}$; that is, for every
$a\in E_{x(t^{'})}$ and $\xi\in E^\ast_{x(t^{'})}$ over the same
base point ${x(t^{'})}\in M$,
$$\< B_{tt^{'}}(a),B^\ast_{tt^{'}}(\xi)>_{\tau}=\< a,\xi>_{\tau}\quad \text{ for every $t\in I$}.$$
Indeed, since by definition the pairing
$\<\cdot,\cdot>_{\sT\tau}:\sT E\times_{\sT M}\sT E^\ast\lra\R$ is
the tangent map of $\<\cdot,\cdot>_\tau:E\times_M E^\ast\lra\R$,
we have
\begin{align*}
&\pa_t\<B_{tt^{'}}(a),B^\ast_{tt^{'}}(\xi)>_{\tau}=\<\pa_tB_{tt^{'}}(a),\pa_tB^\ast_{tt^{'}}(\xi)>_{\sT\tau}\\&=
\<\dd_\sT(f_{u(t)})\left(B_{tt^{'}}(a)\right),\X_{h_t}\left(B_{tt^{'}}^\ast(\xi)\right)>_{\sT\tau}=0.
\end{align*}
\end{remark}\medskip

Finally, observe that the evolution of $\bm\xi(t)=\bm
{B^\ast}_{tt_0}(\bm\xi_0)$ for the control system
\eqref{def:con_sys_A} is trivial on the $\sT^\ast\R$-component.
Indeed, the associated linear Hamiltonian
$$\bm H_t(\bm x,\bm\xi)=\<\bm f(\bm x,u(t)),\bm\xi>_{\bm\tau}=\<f(x,u(t)),\xi>_\tau+\ul\xi L(x,u(t))$$
does not depend on the $\R$-component of $\bm x=(x,\ul x)\in
M\times\R$; hence $\bm\xi(t)=\left(\xi(t),\ul\xi(t)\right)\in\bm
{A^\ast}= E^\ast\times\T^\ast\R$ (in local coordinates,
$(\xi,\ul\xi)\sim \left((x^a,\xi_i),(\ul x,\ul\xi)\right)$)
evolves due to equations
\begin{equation}
\label{eqn:par_trans_*A} \left\{ \begin{aligned}
\pa_t\xi_k(t)=&-\rho^a_k\left(x\right)\left(\frac{\pa f^i}{\pa
x^a}\left(x,u(t)\right)\xi_i(t)+\frac{\pa L}{\pa
x^a}\left(x,u(t)\right)\ul\xi(t)\right)\\ &+c^i_{jk}\left(x\right)f^j\left(x,u(t)\right)\xi_i(t),\\
\pa_t{\ul\xi}(t)=&0,\\
\dot{x}(t)=&\rho\left(f(x(t),u(t))\right),\\
\dot{\ul x}(t)=&L\left(x(t),u(t)\right).
\end{aligned}\right.\end{equation}
In other words, $\ul\xi(t)\equiv\ul\xi_0$ is a constant and
$\xi(t)\in E^\ast$ evolves due to a time-dependent family of
Hamiltonian vector fields $\X_{H_t}$ on $E^\ast$, where
$H_t(x,\xi)=\<f(x,u(t)),\xi>_\tau+\ul\xi_0L(x,u(t))$.

\section{The Maximum Principle}\label{sec:pmp}

In the previous section we have introduced the OCP \eqref{eqn:P}
in the AL algebroid setting. The main difference in comparison
with the classical formulation is, apart from using
algebroid-valued velocities, the fixed-homotopy boundary
condition. This formulation is well motivated by the following
reasoning. First, in the standard case for $E=\sT M$, the
fixed-homotopy condition \eqref{eqn:P2} and fixed end-point
condition \eqref{eqn:P1} are almost equivalent (equivalent after
passing to the universal covering of the control system) and are
not distinguished by the PMP. Second, reducing a fixed-homotopy
condition on a Lie groupoid $\GG$ leads to a fixed-$E$-homotopy
condition on the associated Lie algebroid $E=A(\GG)$.

\begin{theorem}\label{thm:pmp}\index{Pontryagin maximum principle|main}
Let $(x(t),u(t))$, with $t\in[t_0,t_1]$, be a controlled pair of
\eqref{eqn:con_sys} solving the optimal control problem
\eqref{eqn:P}. Then there exists a curve $\xi:[t_0,t_1]\lra
E^\ast$ covering $x(t)$ and a constant $\ul\xi_0\leq 0$ such that the following holds:
\begin{itemize}
    \item the curve $\xi(t)$ is a trajectory of the time-dependent family of Hamiltonian vector fields
    $\X_{H_t}$ associated with Hamiltonians $H_t(x,\xi):=H(x,\xi,u(t))$, where
    $$H(x,\xi,u)=\< f\left(x,u\right), \xi>_\tau+\ul\xi_0 L\left(x,u\right);$$
    \item the control $u$ satisfies the ``maximum principle''
    $$H(x(t),\xi(t),u(t))=\sup_{v\in U}H(x(t),\xi(t),v)$$
and $H(x(t),\xi(t),u(t))=0$ at every regular point $t$ of $u$;
    \item if $\ul\xi_0=0$, the covector $\xi(t)$ is nowhere-vanishing.
\end{itemize}
\end{theorem}

The above result clearly reduces to the standard PMP for the case of the tangent algebroid $E=\T M\lra M$. It also covers the known results for system with symmetry on Lie groups and, more generally, principal bundles. A more detailed discussion and examples will be given in Section \ref{sec:exmples}. 

Theorem \ref{thm:pmp} has an equivalent formulations in terms of the product algebroid $\bm A=E\times\T\R$.

\begin{theorem}\label{thm:pmp_A} Let $(\bm x(t),u(t))$, with $t\in[t_0,t_1]$, be a controlled pair of\index{Pontryagin maximum principle}
\eqref{eqn:con_sys_A} solving the optimal control problem
\eqref{eqn:P_A}. There exists a nowhere-vanishing curve $\bm
\xi=(\xi,\ul\xi):[t_0,t_1]\lra\bm A^\ast= E^\ast\times\sT^\ast\R$
covering $\bm x(t)$, with $\ul\xi(t_1)\leq 0$, such that the following hold:
\begin{itemize}
    \item the curve $\bm\xi(t)$ is a trajectory of the time-dependent family of Hamiltonian vector fields
    $\X_{\bm H_t}$, for
    $\bm H_t(\bm x,\bm\xi):=\bm H(\bm x,\bm\xi,u(t))$, where
    $$\bm H(\bm x,\bm\xi,u)=\<\bm f\left(\bm x,u\right), \bm \xi>_{\bm\tau};$$
    \item the control $u$ satisfies the ``maximum principle''
    $$\bm H(\bm x(t),\bm \xi(t),u(t))=\sup_{v\in U}\bm H(\bm x(t),\bm\xi(t),v)=0$$
at every regular point $t$ of $u$.
\end{itemize}
\end{theorem}

The above formulation will be more useful in the proof. The equivalence of Theorems \ref{thm:pmp} and \ref{thm:pmp_A} is obvious in light of our previous considerations. Indeed, the covector $\bm\xi(t)$ can
be decomposed as $\left(\xi(t),\ul\xi(t)\right)$, and its evolution
along $\X_{\bm H_t}$ is given by \eqref{eqn:par_trans_*A}.
Consequently, as we have observed at the very end of Section
\ref{sec:ocp}, covector $\ul\xi(t)=\ul\xi_0$ is constant and the
evolution of $\xi(t)$ is given by $\X_{H_t}$.  Since
$H(x,\xi,u)+\ul{\xi_0}L(x,u)=\bm H(\bm x,\bm\xi,u)$ for $\bm\xi=(\xi,\ul\xi_0)$, and
$\bm x=(x,\ul x)$, the corresponding statements in Theorems
\ref{thm:pmp} and \ref{thm:pmp_A} are equivalent.

\begin{remark}{\rm \label{rem:pmp_versions}
There are many different versions of the PMP -- for autonomous and
non-autonomous systems, with free or fixed end-points, with free
or fixed time intervals $[t_0,t_1]$, etc. Our results cover the
case of an autonomous OCP with a fixed end-point boundary condition
and an unrestricted time interval. We expect an AL-algebroid version
of the PMP for all other cases, yet in this paper we decided to
restrict our attention only to this single case in order to avoid
further technical complications by defining homotopy classes
relative to an immersed subalgebroid. This study will be
undertaken in a forthcoming work.
}\end{remark}

\section{OCP on Lie groups and principal bundles}\label{sec:exmples}

\subsection{The Atiyah algebroid}
Let $G\ra P\overset\pi\ra M$ be a principal $G$-bundle. Denote by
$p\mapsto pg$ the action of $G$ on $P$ and by $R_g$ the induced
action on the tangent bundle $\sT P$. The $G$-invariant vector
fields on $P$ give rise to a certain Lie algebroid structure,
called the \emph{Atiyah algebroid}. Observe that every
$G$-invariant vector field on $P$ can be canonically identified
with a section of the quotient bundle $E:=\sT P/G\ra M$, and since
$G$-invariant vector fields on $P$ with the Lie bracket of vector
fields $[\cdot,\cdot]_{\sT P}$ form a Lie algebra, we have the
induced Lie algebra structure $(\Sec(E),[\cdot,\cdot]_E)$ on the
space of sections of $E$. The bracket $[\cdot,\cdot]_E$ satisfies
the Leibniz rule \eqref{eqn:lieb_rule} for the natural anchor map
$\rho:E\ra\sT M$ given by $[X]\mapsto \sT\pi(X)$; hence
$(E,[\cdot,\cdot]_E,\rho)$ is a Lie algebroid. Now we shall
investigate this structure in detail.

The action $R_g$ restricts to the space $VP\subset \sT P$ of
vertical vectors (i.e. vectors tangent to the fibers of $\pi$).
Since $VP$ is spanned by the fundamental vector fields of the
$G$-action on $P$, we have a canonical isomorphism $VP\simeq
P\times\g$, where $\g$ is the Lie algebra of $G$. The action $R_g$
in this identification reads as $R_g(p,a)=(pg,Ad_{g^{-1}}a)$, so that
$K:=VP/G\simeq P\times_G\g$, where the right action of $G$ on $\g$
is $g\mapsto Ad_{g^{-1}}$. What is more, the bracket of two
$G$-invariant vertical vector fields on $P$ corresponds, in this
identification, to the canonical Lie bracket $[\cdot,\cdot]_\g$ on
$\g$. That is, if $x\mapsto (x, a(x))$ and $x\mapsto (x,b(x))$ are
two $G$-invariant sections of $P\times\g\simeq VP$, then
$$[(x,a(x)),(x,b(x))]_{\sT P}=(x,[a(x),b(x)]_\g)\in P\times\g.$$
This shows that the bundle $K$ is a Lie algebroid with the trivial
anchor and the Lie algebra structure in fibers isomorphic to $\g$.
Hence, we get the following sequence of Lie algebroid morphisms,
called the \emph{Atiyah sequence}:
$$ 0\to K:=P\times_G{\g}\to E\overset{\rho}{\ra}TM\to 0\,.
$$
Introduce a local trivialization $\phi_s:V\times G\ra P|_V$
obtained from a local section $s:V\ra P$ by the formula
$\phi_s(x,g)=s(x)g$. Clearly, we may identify $E|_V$ with $\sT
P|_{s(V)}$; thus $E|_V\simeq \sT M|_V\times\g$. For two local
sections $x\mapsto (X(x),a(x))$ and $x\mapsto (Y(x),b(x))$ in this
trivialization the Lie bracket reads as
\begin{equation}\label{eqn:bracket1}
\left[(X,a),(Y,b)\right]_E=\left([X,Y]_{\sT
M},[a,b]_\g+X(b)-Y(a)\right)
\end{equation}
and the anchor $\rho\left((X,a)\right)=X$. In fact, this is just
the product of the Lie algebroids $\sT M|_V\lra V$ and $\g$ (cf.
the last paragraph of Section \ref{sec:ala}). We have a similar
description globally if the principal bundle $P$ is trivial. Note,
however, that, as we will work with principal bundles over a
neighborhood of a path in $M$, we can always assume that our
principal bundle is trivial.

In some applications one has to work with a principal connection
on $P$. It corresponds to a $G$-invariant horizontal distribution
in $\sT P$ and is represented by a splitting $E=\sT M\oplus_M K$
given by a bundle embedding $\nabla:\sT M\ra E$ such that
$\rho\circ\nabla=id_{\sT M}$. If the bundle $K$ is trivial, then
we get another trivialization $E\simeq_\nabla\sT M\times\g$,
associated with the connection $\nabla$, in which the Lie bracket
on sections reads as
\begin{equation}\label{eqn:bracket}
\left[(X,a),(Y,b)\right]_{E}=\left([X,Y]_{\sT M},
F_\nabla(X,Y)+[a,b]_\g+X(b)-Y(a)\right),
\end{equation}
where $F_\nabla$ is the \emph{curvature} of the connection
$\nabla$, i.e.
$$F_\nabla(X,Y)=[\nabla(X),\nabla(Y)]_E-\nabla\left([X,Y]_{\sT M}\right).$$

\subsection{PMP on an Atiyah algebroid}
The already proven results on the Lie groupoid -- Lie algebroid
reduction of a control system and homotopy (cf. Theorem
\ref{thm:int_htp}, Corollary \ref{cor:htp_P} and Section
\ref{sec:ocp}) allow us to formulate the following result which can
be understood as a general reduction scheme of the PMP on a
principal $G$-bundle. Consider a principal $G$-bundle $G\ra
P\overset \pi\ra M$.

\begin{theorem}\label{thm:OCP_P}
Let $F:P\times U\lra \sT P$ be a $G$-invariant control system on
$P$ and let $L:P\times U\lra\R$ be a $G$-invariant cost function.

Then the OCP \eqref{eqn:P2} on $P$ for a path joining two fixed
points is equivalent to the OCP \eqref{eqn:P} for the system
$f:M\times U \lra E=\sT P/G$ where $f(\pi(p),u):=[F(p,u)]$, with
the cost function $l:M\times U\lra\R$ defined by
$l(\pi(p),u):=L(p,u)$ and the $E$-homotopy class $[\sigma]$ being
the reduction  of the corresponding fixed homotopy class in $P$.
\end{theorem}

\medskip
For the OCP on $E$ described above we can use
Theorem \ref{thm:pmp}. This theorem does not differ much from the
standard PMP (the definition of the pseudo-Hamiltonian function,
the maximum condition, non-vanishing of the covector, etc.) apart
from using the different Poisson structure to describe the
Hamiltonian evolution of the covector $\xi(t)\in E^\ast$. We shall
now find this structure for the descriptions of the Atiyah
algebroid by means of the principal connection introduced in the
previous paragraph.

Suppose for simplicity that the vertical subbundle $K=VP/G$ is
trivial (e.g. $P$ is trivial), $K=M\times\g$, and consider a
splitting $E\simeq_\nabla\sT M\oplus_MK$ induced by a principal
connection $\nabla:\sT M\ra E$, so that we get an identification
$E\simeq_\nabla \sT M\times\g$. Let $E^\ast\simeq_\nabla\sT^\ast
M\times\g^\ast$ be the corresponding identification of the dual
bundle. The proof of the following theorem is straightforward by
explicit coordinate calculations (cf. Section \ref{sec:ala}).
\begin{theorem}\label{thm:ham_E1}
The Poisson tensor $\Pi_E$ associated with the Lie algebroid
structure on $E$ in the identification
$E^\ast\simeq_\nabla\sT^\ast M\times\g^\ast$ reads as
$$\Pi_E(p_x,\zeta)=\left(\Pi_{\sT^\ast M}(p_x)+V_{\<\zeta,F_{\nabla}(x)(\cdot,\cdot)>}\right)\times \Pi_{\g^\ast}(\zeta),$$
where $p_x\in T^\ast_x M$, $\xi\in\g^\ast$, $\Pi_{\sT^\ast M}$ and
$\Pi_{\g^\ast}$ are the standard Poisson tensors, and
$V_{\<\zeta,F_{\nabla}(x)(\cdot,\cdot)>}$ is the $2$-form
$\<\zeta,F_{\nabla}(x)(\cdot,\cdot)>$ associated with the
curvature $F_\nabla(x):\bigwedge^2\sT_x M\ra \g$ understood as a
vertical tensor on $\sT^\ast M$.

The Hamiltonian vector field defined by means of $\Pi_E$ and a
Hamiltonian $h:\sT^\ast M\oplus\g\lra\R$ reads as
$$\X^E_h(p_x,\zeta)=\left(\X^{\sT^\ast M}_{h(\cdot,\zeta)}(p_x)+V_{\<\zeta,F_{\nabla}(x)(\frac{\pa h}{\pa p}(p_x,\zeta),\cdot)>},
\X^{\g^\ast}_{h(p_x,\cdot)}(\zeta)\right).$$ In local coordinates,
$p\sim(x^a,p_b)$ and $\zeta\sim(\zeta_\alpha)$,
\begin{eqnarray*}\X^E_h(x,p,\zeta)&=\frac{\pa h}{\pa p_a}(x,p,\zeta)\pa_{x^a}+\left(
\zeta_\alpha F^\alpha_{ab}(x)\frac{\pa h}{\pa
p_a}(x,p,\zeta)-\frac{\pa h}{\pa x^b}(x,p,\zeta)\right)\pa_{p_b}\\&+
\zeta_\gamma C^\gamma_{\alpha\beta}\frac{\pa
h}{\pa\zeta_\alpha}\pa_{\zeta_\beta}\,,
\end{eqnarray*}
where $F^\alpha_{ab}(x)$
are the coefficients of the curvature $F_\nabla$ and
$C^\alpha_{\beta\gamma}$ are the structure constants of $\g$.
\end{theorem}

Note that for the case of a Lie group ($P=G$, $M=\{\ast\}$) we recover the results of Jurdjevic \cite[Chap.12,
Thms 5,6]{jurdjevic}.

\subsection{The falling cat problem}
Now we will reconsider the well-known results of Montgomery
\cite{montgomery_isohol} (see also \cite[ch. 7.1]{bloch} and \cite{cendra_holm_marsden}) on
the isoholonomic problem by means of the PMP in the Atiyah
algebroid setting.

Let $G\ra P\ra M$ be a principal $G$-bundle, let $\cH\subset \sT P$ be
a $G$-invariant horizontal distribution, and let $\mu(\cdot,\cdot)$ be
a $G$-invariant sub-Riemannian metric on $\cH$ ($\mu(\cdot,\cdot)$
can be understood as a base metric lifted to $\cH$ by the
horizontal lift). The problem is now to find a horizontal curve
$q(t)$ with $t\in[0,1]$ joining two fixed points $q_0$, $q_1$ in
$P$ and minimizing the total energy
$$\frac 12\int_0^1\mu(\dot q(t),\dot q(t))\dd t.$$

Clearly, due to the $G$-invariance of the problem, after changing
the fixed-end-point condition into a fixed-homotopy condition (as
discussed in detail in Section \ref{sec:ocp}), the above problem
is equivalent to an OCP of the form \eqref{eqn:P} on the Atiyah
algebroid $E=\T P/G$.

With the invariant distribution $\cH$, understood as a principal
connection, we can associate a map $\nabla:\sT M\ra E$ inducing a
splitting $E\simeq_\nabla\sT M\times\g$. Our control system will
be $f:\sT M\lra \sT M\times\g$ given by $f(X)=(X,0)$ (this assures
that the trajectory is horizontal), the cost function $L:\sT
M\lra\R$ reads as $L(X)=\frac 12 \mu(X,X)$, and the fixed
$E$-homotopy class is simply a reduction of a classical homotopy
class in $P$.

Note two differences with the formulation of the OCP
\eqref{eqn:P}. Firstly, our control and cost functions have
arguments in $\sT M$ instead of in $M\times U$. Of course, this makes
no big difference, since locally $\sT M\approx M\times\R^n$.
Secondly, our time interval is fixed. This, in turn, results in
substituting the condition $H(x(t),\xi(t),u(t))=0$ by
$H(x(t),\xi(t),u(t))=\mathrm{const}$ in the assertion of Theorem
\ref{thm:pmp}.

Now we can apply Theorem \ref{thm:pmp} with the Hamiltonian
evolution described in Theorem \ref{thm:ham_E1}. The covector
$\xi\in E^\ast$ can be decomposed as $\xi=(p,\zeta)\in\sT^\ast
M\times\g^\ast$, and the corresponding Hamiltonian is
$$H(p,\zeta,X)=\<X,p>+\frac 12\lambda_0\mu(X,X)=:h(p,X),$$
with $\lambda_0\leq 0$. The maximum principle reads as
$p(t)=-\lambda_0\mu(X(t),\cdot)$; hence on the optimal trajectory,
$H(p(t),\zeta(t),X(t))=-\frac 12\lambda_0\mu(X(t),X(t))$ (which is
constant in $t$). The evolution of $p(t)$ and $\zeta(t)$ is given
by
\begin{align*}
&\dot\zeta(t)=0\,,\\
&\dot p(t)=\X^{\sT^\ast M}_{h(p,X)}+\<\xi(t),F_\nabla(X,\cdot)>\,;
\end{align*}
hence $\zeta(t)=\mathrm{const}$. The second equation is equivalent to
$$\lambda_0\nabla_X^\mu X=\<\zeta,F_\nabla(X,\cdot)>^{\#\mu},$$
where $\nabla^\mu$ denotes the Levi-Civita covariant derivative on
$(M,\mu)$ and $A^{\#\mu}$ is the vector dual to $A$ by means of
$\mu$. Indeed, the equation $\dot p(t)=\X^{\sT^\ast M}_{h(p,X)}$,
together with $p(t)=-\lambda_0\mu(X(t),\cdot)$, is the PMP for a
geodesic problem on $(M,\mu)$. Passing to the dual vector
$p(t)^{\#\mu}=\lambda_0X(t)$ we should obtain the geodesic
equation multiplied by the factor $\lambda_0$. The equation
$\dot\zeta(t)=0$ means that the curve $\zeta(t)\in\g^\ast$ is
covariantly constant, hence
$$\nabla_X\zeta=0.$$
We have thus obtained the Wong equations as in \cite{montgomery_isohol}.

The abnormal case $\lambda_0=0$ implies $p(t)=0$ and
$\<\zeta,F_\nabla(X,\cdot)>=0$. This allows us to exclude abnormal
solutions in certain situations. For example, if $P$ is a bundle
of circles over a two-dimensional base and the connection is
non-integrable (i.e., $F_\nabla$ is non-vanishing), we have
$\<\zeta,F_\nabla(X,\cdot)>=0$ if and only if $X=0$ (hence the
solution is trivial) or $\zeta=0$, which can be excluded by the
non-vanishing of the covector in the PMP.

\subsection{Lagrangian reduction, Hammel equations, Euler-Poincare equations}
As a special case of the OCP from Theorem \ref{thm:OCP_P},
consider a $G$-invariant variational problem on $P$ with a
Lagrangian $L:\sT P\lra\R$. By a standard argument this problem is
equivalent to an unconstrained $G$-equivariant OCP on $P$ which
can be reduced to an OCP on the Atiyah algebroid $E=\sT P/G$. The
control function is $f=\id_E:E\lra E$, and the reduced Lagrangian
$l:E\lra \R$ is the cost function. As there are no constraints on
the set of velocities, the abnormal case can be excluded. Now,
using the PMP on $E$ in the setting of Theorem \ref{thm:ham_E1},
one will obtain the reduced Euler-Lagrange equations of \cite[Sect.
5]{cendra_holm_marsden}. A similar calculation (for a principal connection
induced by a local trivialization) will give the Hammel equations
\cite[Sect. 5]{cendra_holm_marsden}. As a special case, when $E=\g$, one gets the
Euler-Poincare equations. The exact calculations are easy and are
left to the reader.

\subsection{Two-point time OCP on $\so$}
Consider now a rigid body in $\R^3$ which can rotate with constant
angular velocity along two fixed axes in the body. At every moment
the position of the body is described by an element $q\in SO(3)$.
The rotation axes can be represented by elements of the Lie
algebra $l_+,l_-\in\so$. The rotation along the axis $l_\pm$ is
described by the equation
$$\pa_t q=ql_\pm.$$
It would be suitable to write $l_+=a+b$ and $l_-=a-b$. The above
equation can be regarded as a control system on the Lie group
$SO(3)$ with the control function $F(q,u)=qf(u)$, where
$f(u)=a+ub$ and the set of controls is simply $U=\{-1,1\}$. We
would like to find a control $u(t)$ which moves the body from a
position $q_0\in SO(3)$ to $q_1\in SO(3)$ (or such that the
trajectory belongs to a fixed homotopy class in $SO(3)$) in the
shortest possible time.

It is obvious that the above OCP on the Lie
group reduces to the OCP on the Lie algebra
$\so$ with the control function $f$ and the cost function $L\equiv
1$. Fix a basis $(e_1,e_2,e_3)$ on $\so$, and denote by
$c^\alpha_{\beta\gamma}$ the structure constants of the Lie
algebra in this basis. Let $u(t)$, for $t\in[t_0,t_1]$, be a
solution of the above OCP. It follows from
Theorem \ref{thm:pmp} that there exist a number $\lambda_0\leq 0$
and a curve $\zeta(t)\in \so^*$ such that
$$H(\zeta(t),u(t))=\< \zeta(t),a+u(t)b>+\lambda_0=\max_{v=\pm 1}\< \zeta(t),a+vb>+\lambda_0.$$
This implies that $u(t)=\operatorname{sgn}\big(\<
\zeta(t),b>\big)$. Moreover, the evolution of $\zeta(t)$ is given
by the equation
$$\pa_t \zeta_\beta(t)=c^\gamma_{\alpha\beta}(a^\alpha+u(t)b^\alpha)\zeta_\gamma(t).$$
We have obtained the same equation as in (\cite[Sec.
19.4]{agrachev}). We refer the reader to this book for the detailed discussion on
solutions.

\subsection{An application to a nonholonomic system}

In \cite{GGU_geom_mech} and \cite{GG_var_calc} a framework of geometric mechanics on
general algebroids was presented. Roughly speaking, the structure
of an algebroid on a bundle $\tau:E\ra M$ allows one to develop
Lagrangian formalism for a given Lagrangian function $L:E\ra\R$.
Moreover, if $E$ is an AL algebroid, then the associated
Euler-Lagrange equations have a variational interpretation: a
curve $\gamma:[t_0,t_1]\ra E$ satisfies the Euler-Lagrange
equations if and only if it is an extremal of the action
$\mathcal{J}(\gamma):=\int_{t_0}^{t_1}L(\gamma(t))\dd t$
restricted to those $\gamma$'s which are admissible and belong to
a fixed $E$-homotopy class \cite{GG_var_calc}.  Hence, the trajectories of
the Lagrange system should be derivable from our version of the PMP
for the unconstrained control system on $E$ with the cost function
$L$.

In \cite{grabowski_nonholonomic} it has been shown that if $D\subset E$ is a
subbundle and $L$ is of mechanical type (that is, $L(a)=\frac 12
\mu(a,a)-V(\tau(a))$, where $\mu$ is a metric on $E$ and $V$ is an
arbitrary function on the base), then nonholonomically constrained
Euler-Lagrange equations associated with $D$ can be obtained as
unconstrained Euler-Lagrange equations on the skew-algebroid
$\left(D,\rho_E|_D,[\cdot,\cdot]_D:=\PP_D[\cdot,\cdot]_E\right)$,
where $\PP_D:E\ra D$ denotes the projection orthogonal w.r.t.
$\mu$. It follows that if $D$ with the algebroid structure
defined above is AL, then the solutions of the nonholonomically
constrained Euler-Lagrange equations are extremals of the
unconstrained OCP on $D$ with the cost function $L|_D$. On the
other hand, using our version of the PMP on the algebroid $E$ with
controls restricted to $D$ and the cost function $L$, one will
obtain nonholonomically constrained Euler-Lagrange equations
associated with $D$. Note that the algebroid bracket
$[\cdot,\cdot]_D$ need not satisfy Jacobi identity even if
$[\cdot,\cdot]_E$ does. Concluding, the PMP on general (not
necessarily Lie) AL algebroids can be used in the theory of
nonholonomic systems. To our knowledge this point of view is
completely novel.

To give a concrete example we will use PMP to study the Chaplygin
sleigh. It is an example of a nonholonomic system on the Lie
algebra $\mathfrak{se}(2)$ which describes a rigid body sliding on
a plane. The body is supported in three points, two of which
slide freely without friction, while the third point is a knife
edge. This imposes the constraint of no motion orthogonal to this
edge (see \cite{chaplygin,neimark}).

The configuration space before reduction is the Lie group
$G=SE(2)$ of the Euclidean motions of the two-dimensional plane
$\R^2$. Elements of the Lie algebra $\mathfrak{se}(2)$ are of the
form
$$\hat{\xi}=
\begin{pmatrix}
0&\omega&v_1\\
-\omega&0&v_2\\
0&0&0
\end{pmatrix}=v_1E_1+v_2E_2+\omega E_3,
$$
where $[E_3,E_1]=E_2$, $[E_2, E_3]=E_1$, and  $[E_1, E_2]=0$.

The system is described by the purely kinetic Lagrangian function
$L:\mathfrak{se}(2)\ra\R$, which reads as
$$L(v_1, v_2, \omega)=\frac{1}{2}\left[ (J+m(a^2+b^2))\omega^2 + mv_1^2+m v_2^2-2bm\omega v_1-2am\omega v_2\right].$$
Here $m$ and $J$ denote the mass and the moment of inertia of
the sleigh relative to the contact point, while $(a, b)$
represents the position of the centre of mass w.r.t. the
body frame, determined by placing the origin at the contact point
and the first coordinate axis in the direction of the knife axis.
Additionally, the  system is subjected to the nonholonomic
constraint determined by the linear subspace
$$ D=\{(v_1, v_2, \omega)\in \mathfrak{se}(2)\; |\; v_2=0\}\subset\mathfrak{se}(2).$$
Instead of $\{E_1, E_2, E_3\}$ we take another basis of
$\mathfrak{se}(2)$:
$$
e_1=E_3,\quad e_2=E_1,\quad e_3= -ma E_3-mab E_1+(J+ma^2) E_2,$$ adapted
to the decomposition $D\oplus D^\perp$; $D=\hbox{span }\{ e_1,
e_2\}$ and $D^\perp=\hbox{span }\{ e_3\}$. The induced
skew-algebroid structure on $D$ is given by
$$
[e_1, e_2]_{D}=\frac{ma}{J+ma^2} e_1+\frac{mab}{J+ma^2}e_2.$$
 Therefore, the structural constants are ${\mathcal C}^1_{12}=\frac{ma}{J+ma^2}$ and
${\mathcal C}^2_{12}=\frac{mab}{J+ma^2}$. The algebroid $D$ is
almost Lie (in fact, in this simple case it is a Lie algebra).
Next, we will use theorem \ref{thm:pmp} to derive the nonholonomic
equations of motion. Set $U=\R^2\ni(y^1,y^2)$ and the control
function to be a map $f:U\ra D$ given by
$$ f(y^1,y^2)=y^1e_1+y^2e_2\in D.$$
The Lagrangian restricted to $D$ defines the cost function
$L:U\ra\R$,
$$
L(y^1, y^2)=\frac{1}{2}\left[ (J+m(a^2+b^2))(y^1)^2 +
m(y^2)^2-2bmy^1y^2\right].$$ For a curve
$\xi(t)=\xi_1(t)e^\ast_1+\xi_2(t)e^\ast_2\in D^\ast$ the maximum
principle reads
\begin{eqnarray} \label{eqn:max_nh} H(\xi(t),y(t))&=&\xi_1(t)y^1+\xi_2(t)y^2+\ul\xi_0 \cdot L(y^1,y^2)\\&=&
\max_{(v^1,v^2)\in\R^2}(\xi_1(t)v^1+\xi_2(t)v^2+\ul\xi_0 \cdot
L(v^1,v^2))\,.\nonumber\end{eqnarray} 
If $\ul\xi_0=0$, then maximality would
give $\xi(t)=0$, which is impossible. Hence, we may assume that
$\ul\xi_0=-1$. Now from \eqref{eqn:max_nh} we will get
\begin{equation}\label{eqn:ev_nh}
\begin{split}
\xi_1(t)&=\left(J+m(a^2+b^2)\right)y^1-bmy^2\,,\\
\xi_2(t)&=my^2-bmy^2\,.
\end{split}
\end{equation}
Finally, the Hamiltonian evolution \eqref{eqn:par_trans_*A} is
simply
\begin{align*}
\dot \xi_1&=\mathcal{C}^1_{21}y^2\xi_1+\mathcal{C}^2_{21}y^2\xi_2=-\frac{ma}{J+ma^2}y^2(\xi_1+b\xi_2),\\
\dot
\xi_1&=\mathcal{C}^1_{21}y^1\xi_1+\mathcal{C}^2_{12}y^1\xi_2=\frac{ma}{J+ma^2}y^1(\xi_1+b\xi_2).
\end{align*}
In view of \eqref{eqn:ev_nh} and the above equations we conclude
that the equations of motion are
\begin{align*}
(J+m(a^2+b^2))\dot{y}^1 -bm \dot{y}^2&= -ma y^1y^2,\\
m\dot{y}^2-bm\dot{y}^1&=ma(y^1)^2,
\end{align*}
which completely agrees with \cite{grabowski_nonholonomic}.

\section{Needle variations and \emph{E}-homotopy}\label{sec:needle}

\subsection{The variation of controls and trajectories}

In order to prove Theorem \ref{thm:pmp_A} we shall somehow compare the cost on the optimal trajectory $\bm f(\bm x(t),u(t))$ with costs of nearby trajectories. As our assumptions input on the
set of controls $U$ are very mild, we cannot use the natural
concept of a continuous deformation, as in the standard calculus
of variations ($U$ can be for instance discrete). Instead, we introduce
the notion of \emph{needle variations} after \cite{pontryagin}.
For a given admissible control $u:[t_0,t_1]\lra U$ this variation
will be, roughly speaking, the family of controls $ u_s(t)$
obtained by substituting $u(t)$ by given elements $v_i\in U$ on a
small intervals $I_i=(\tau_i-s\del t_i,\tau_i]\subset[t_0,t_1]$. Our main result in this section is Theorem \ref{thm:1st_main}, where we study the $\bm A$-homotopy classes of the trajectories of the system \eqref{eqn:con_sys_A} obtained for controls $u_s(t)$. 
We finish this section with the definition of $\bm K_\tau^u$---the set of infinitesimal variations of the trajectory $\bm f(\bm x(t), u(t))$.

\subsection{Needle variation of controls and trajectories}

Throughout this section we will work with a fixed admissible control $u:[t_0,t_1]\lra U$ and fixed trajectory $\bm a(t):=\bm f(\bm x(t),u(t))$. 

Choose points $t_0<\tau_1\leq\tau_2\leq\hdots\leq\tau_k\leq\tau<t_1$,
being regular points of $u$. Next, choose non-negative numbers
$\delta t_1,\hdots,\delta t_k$ and an arbitrary real number
$\delta t$. Finally, take (not necessarily different) elements
$v_1,\hdots,v_k\in U$. The whole set of data $(\tau_i,
v_i,\tau,\delta t_i, \delta t)_{i=1,\hdots,k}$ will be denoted by
$\A$ and called a \emph{symbol}. Its role will be to encode the variation of the control $u(t)$. Intuitively, points $\tau_i$ emphasise moments in which we substitute $u(t)$ by a constant control $v_i$ on an interval $I_i=(\tau_i-s\del t_i,\tau_i]$ of length $s\del t_i$, while $s\del t$ is responsible for shortening or lengthening the time for which $u(t)$ is defined. The precise
definition is quite technical, because one should take care to make the intervals $I_i$  pair-wise disjoint.

Take
$$l_i=\begin{cases}
\del t-(\del t_i+\hdots+\del t_k) &\text{when $\tau_i=\tau$;}\\
\phantom{x.}-(\del t_i+\hdots+\del t_k) &\text{when $\tau_i=\tau_k<\tau$;}\\
 \phantom{x.}-(\del t_i+\hdots+\del t_j) &\text{when
$\tau_i=\tau_{i+1}=\hdots=\tau_j<\tau_{j+1}$,}
\end{cases} $$
and define $s$-dependent intervals  $I_i:=(\tau_i+s l_i,\tau_i+s(l_i+\del t_i)]$. As we see, $I_i$ is left-open and right-closed and it has
length $s\cdot\delta t_i$. If $\tau_i<\tau_{i+1}$, or $i=k$ and
$\tau_k<\tau$, the end-point of $I_i$ lies at $\tau_i$. If
$\tau_i=\tau_{i+1}$, then the end-point of $I_i$ coincides with the
initial-point of $I_{i+1}$. If $\tau_k=\tau$, we set the
end-point of $I_k$ at $\tau+s\delta t$. Clearly, for $s$ small
enough, the intervals $I_i$ lie inside $[t_0,t_1]$ and are
pairwise disjoint.


\begin{definition}\label{def:need_var} For a symbol $\A=(\tau_i, v_i,\tau,\delta t_i, \delta t)_{i=1,\hdots,k}$ we introduce a $s$-dependent
family of admissible controls defined on intervals
$[t_0,t_1+s\delta t]$:
\begin{equation}
\label{eqn:need_var} u^\A_s(t)=
\begin{cases}
v_i & \text{for $t\in I_i$},\\
u(t) & \text{for $t\in [t_0,\tau+s\del t]\setminus \bigcup_i
I_i$}\\
u(t-s\del t) &\text{for $t\in(\tau+s\del t,t_1+s\del t]$}.
\end{cases} \end{equation}
We will call $u^\A_s$ a (\emph{needle})
\emph{variation of the control} $u$ \emph{associated with the symbol
$\A$}. 
\end{definition}

Using $u^\A_s(t)$ and an AC path $s\mapsto \bm x_0(s)\in M\times\R$ where $\bm x_0(0)=\bm x_0$ we can define the variation of $\bm a(t)$.

\begin{definition}
The family of trajectories
$$\bm{a}(t,s):=\bm{f}\left(\B{x}(t,s),u^\A_s(t)\right)$$
of the system \eqref{eqn:con_sys_A}, with the initial conditions $\bm x(t_0,s)=\bm
x_0(s)$, where $t\in[t_0,t_1+s\del t]$, will be called a
\emph{variation of the trajectory} $\bm a(t)=\bm f(\bm
x(t),u(t))$ \emph{associated with the symbol $\A$ and the initial base-point variation $\bm x_0(s)$}.
\end{definition}

\begin{remark}\label{rem:var_dt=0}
Observe that, when $\delta t_i=0$, the interval $I_i$ is empty.
It follows that adding a triple $(\tau_i, v_i, \del t_i=0)$ to the symbol $\A$ does not change the variation 
$u^\A_s$ and, consequently, the associated variations $\bm a(t,s)$.
\end{remark}

\subsection{The associated $\bm A$-homotopy}

Our goal now is to compare the $\bm A$-homotopy classes of the trajectory $\bm a(t)$ and its variation $\bm a(t,s)$ introduced above. We need this because the OCP \eqref{eqn:P} is defined in term of algebroid homotopy classes. Having in mind Lemma \ref{lem:gen_E_htp} and the construction of a $\bm A$-homotopy associated with a control system \eqref{eqn:con_sys_A} given in Section \ref{sec:ocp}, we may expect that the family of trajectories $\bm a(t,s)$ forms an $\bm A$-homotopy (for some initial-point homotopy $\bm b_0(s)$). Consequently, the description of $\bm A$-homotopy classes of $\bm a(t,s)$ should be possible by meas of Lemma \ref{lem:htp}. This is indeed the case, yet some technical work is needed in order to reparametrise $\bm a(t,s)$ in a suitable way.

\begin{theorem}\label{thm:1st_main}  Let $s\mapsto \bm b_0(s)$ be a  bounded measurable $\bm A$-path covering $s\mapsto \bm x_0(s)$, where $\bm x_0(0)=\bm x_0$. Consider a variation $\bm a(t,s)=\bm f(\bm x(t,s),u^\A_s(t))$ of the trajectory $\bm a(t)=\bm f(\bm x(t),u(t))$ associated with a symbol $\A=(\tau_i, v_i,\tau,\delta t_i, \delta t)_{i=1,\hdots,k}$ and initial base-point variation $\bm x_0(s)$. 

Then there exists a number $\theta>0$ and an $\bm A$-path $s\mapsto \bm d^{\A}(s)$ defined for $0\leq s\leq\theta$ such that
\begin{equation}\label{eqn:htp_1st_lem}
[\bm b_0(s)]_{s\in[0,\eps]}[\bm{a}(t,\eps)]_{t\in[t_0,t_1+\eps\del t]}=[\bm
{a}(t)]_{t\in[t_0,t_1]}[\bm{d}^{\A}(s)]_{s\in[0,\eps]},\end{equation} 
for every $\eps\leq \theta$. 

Moreover, if $(\tau_i,v_i,\tau)$ in $\A$ are fixed, we can choose $\theta>0$ universal for all $(\del t_i,\del t)$ belonging to a fixed compact set.

Finally, if $\bm b_0(s)$ is regular at $s=0$, then  $\bm d^{\A}(s)$, regarded as a function of $s$, $\del t_i$ and $\del t$, is uniformly regular w.r.t. $\del t_i$ and $\del t$ at $s=0$. What is more,
\begin{equation}
\label{eqn:d_at_0} 
\begin{split}\bm{d}^{\A}(0)=&\bm B_{t_1\tau}[\bm{f}(\bm{x}(\tau),u(\tau))]\del t+\bm{B}_{t_1 t_0}(\bm
b_0(0))\\&+\sum_{i=1}^k
\bm{B}_{t_1\tau_i}\Big[\bm{f}(\bm{x}(\tau_i),v_i)-\bm{f}(\bm{x}(\tau_i),u(\tau_i))\Big]\del
t_i\in\bm{A}_{\bm{x}(t_1)}.
\end{split}
\end{equation} 
\end{theorem}

\begin{proof}
The proof is technically complicated, yet conceptually not very difficult. The idea is to decompose
$\bm a(t,s)$ into several parts, which, after a
suitable reparametrisation, form an $\bm A$-homotopy. As one may have expected, these parts correspond to ''switches'' in the needle variation associated with the symbol $\A$. Our argument will be therefore inductive w.r.t. $k$---the number of ''switches'' in $\A$. Formula \eqref{eqn:htp_1st_lem} will be obtained from the repetitive usage of Lemma \ref{lem:htp} for the partial homotopies, and  \eqref{eqn:d_at_0} will follow from the concrete form of these homotopies. The preservation of the uniform regularity will be obtained using the technical results introduced in Section \ref{sec:meas}. 

Finally, let us explain the role of the number $\theta$. We know from Theorem \ref{thm:exist} that if a solution of the ODE for a fixed initial condition $\bm x_0$ is defined on an interval $[t_0,t_1]$, then so are the solutions for initial conditions close enough to $\bm x_0$. Since the base variation $\bm x(t,s)$ associated with $u^\A_s(t)$ is obtained as a composition of the solutions of \eqref{eqn:con_sys} with perturbations on intervals of length $s\del t_i$ and $s\del t$, it is clear that, if numbers $\del t_i$ and $\del t$ are bounded and $\tau_i$, $\tau$ and $v_i$ fixed, for a given $\bm x_0(s)$, we can chose $\theta>0$ such that the trajectory $\bm x(t,s)$ will stay close enough to $\bm x(t)$ to be well-defined for all $0\leq s\leq\theta$ and all $[t_0,t_1]$. 

In our inductive reasoning it will be more convenient to assume that all the data depends on an additional parameter $p\in P$ (i.e., we have $\bm x_0(s,p)$ instead of $\bm x_0(s)$, $\bm a(t,s,p)$ instead of $\bm a(t,s)$, etc.). In the assertion we demand that \eqref{eqn:htp_1st_lem} and \eqref{eqn:d_at_0} hold for each fixed $p\in P$. Moreover, for fixed $(\tau_i,v_i,\tau)$ we want $\bm d^{\A,p}(s)$ to be uniformly regular w.r.t. $p$, $\del t_i$, and $\del t$ at $s=0$ if $\bm b_0(s,p)$ is uniformly regular w.r.t. $p$ at $s=0$.

In what follows we will need two technical lemmas. 
\begin{lemma}\label{lem:sub_lem1} Let $t\mapsto\bm a(t,s,\wt p)= \bm f(\bm x(t,s,\wt p),v(t))$, with $t\in[\wt t_0,\wt t_1]$, be a family of bounded measurable admissible paths over $\bm x(t,s,\wt p)$ parameterised by $\wt p\in \wt P$. Let $s\mapsto \wt{\bm b}_0(s,\wt p)$ be a family of  bounded measurable $\bm A$-paths over $\bm x(t_0,s,\wt p)$. There exists a number $\theta>0$ and a family of bounded measurable $\bm A$-paths $s\mapsto\bm d_1^{\wt p}(s)$ defined for $0\leq s\leq\theta$ such that
\begin{equation}\label{eqn:needle1}
[\wt{\bm b}_0(s,\wt p)]_{s\in[0,\eps]}[\bm{f}(\bm x(t,\eps,\wt p),v(t))]_{t\in[\wt t_0,\wt t_1]}=[\bm
{f}(\bm x(t,0,p),v(t))]_{t\in[\wt t_0,\wt t_1]}[\bm{d}_1^{\wt p}(s)]_{s\in[0,\eps]},\end{equation} 
for all $\eps\leq \theta$. 
 
Explicitly, $\bm d_1^{\wt p}(s)=\bm B^v_{t\wt{t_0}}\left[\wt{\bm b}_0(s,\wt p)\right]$, where $\bm B^v_{t\wt{t_0}}$ is a parallel transport operator associated with the control $v(t)$. Moreover, if $\wt{\bm b}_0(s,\wt p)$ is uniformly regular w.r.t. $\wt p$ at $s=0$, then so is $\bm d_1^{\wt p}(s)$.
\end{lemma}
\begin{proof}[Proof of the lemma]
Fix $\wt p\in\wt P$ and consider $\bm b(t,s,\wt p):=\bm B^v_{t\wt t_0}\left[\wt{\bm b}_0(s,\wt p)\right]$. It follows from the definition of the operator of parallel
transport $\bm B^v_{t\wt{t_0}}$ that the pair $(\bm a(t,s,\wt p),\bm b(t,s,\wt p))$ is and $\bm A$-homotopy over $\bm x(t,s,\wt p)$ (see Remark \ref{rem:B_htp}).  Now \eqref{eqn:needle1} follows directly from Lemma \ref{lem:htp}, since $\bm b(\wt t_1,s,p)=\bm B^v_{\wt t_1\wt t_0}\left[\wt{\bm b}_0(s,\wt p)\right]=\bm d_1^{\wt p}(s)$. 

Finally, since by Remark \ref{rem:B_htp} the map $\bm B^v_{\wt t_1\wt t_0}(\cdot)$ is  continuous  for every fixed $\wt t_1$ and $\wt t_0$, in light of Lemma \ref{lem:ur_comp}, it preserves the uniform regularity of $\wt{\bm b}_0(s,\wt p)$ .
\end{proof}

\noindent The second lemma is the following one.

\begin{lemma}\label{lem:sub_lem2} Let $t\mapsto\bm a(t,s,\wt p)= \bm f(\bm x(t,s,\wt p),v(t))$ be a family of bounded measurable admissible
paths over $\bm x(t,s,\wt p)$ parametrised by $\wt p\in P$. Let $s\mapsto \wt{\bm b}_0(s,\wt p)$ be a family of bounded measurable $\bm A$-paths over $\bm x(\wt t_0+sc,s,\wt p)$. Then there exists a number $\theta>0$ and a family of bounded measurable $\bm A$-paths $s\mapsto\bm d_2^{\wt p,c,d}(s)$ defined for $0\leq s\leq\theta$ such that
\begin{equation}\label{eqn:needle3}
[\wt{\bm b}_0(s,\wt p)]_{s\in[0,\eps]}[\bm{f}(\bm x(\wt t_0+t,\eps,\wt p),v(t))]_{t\in[\eps c,s\eps d]}=[\bm{d}_2^{\wt p,c,d}(s)]_{s\in[0,\eps]},\end{equation} 
for every $\eps\leq \theta$. 
 
 Moreover, if $\wt{\bm b}_0(s,\wt p)$ is uniformly regular  w.r.t. $\wt p$ at $s=0$ and $\wt t_0$ is a regular point of the control $v(t)$,  then $\bm d_2^{\wt p,c,d}(s)$ is uniformly regular  w.r.t. $\wt p$, $c$, and $d$ at $s=0$. Finally, 
\begin{equation}\label{eqn:needle3a}
\bm{d}_2^{\wt p,c,d}(0)=\wt{\bm b}_0(s,\wt p)+(d-c)\bm f(\bm x(\wt t_0,0,\wt p),v(\wt t_0)).\end{equation}
\end{lemma}

\begin{proof}[Proof of the lemma]
For notation simplicity let forget about the $\wt p$-dependence. Define
$$\wh {\bm b}_0(s):=\bm B_{\wt t_0(\wt t_0+sc)}\left[\wt{\bm b}_0(s)-c\bm f(\bm x(\wt t_0+cs),v(\wt t_0+cs))\right].$$
Clearly, $\wh{\bm b}_0(s)$ is an admissible paths over $\bm x(\wt t_0,s)$. Now define a pair of maps
\begin{align*}
\bm a(t,s)&=s\bm f\left(\bm x(\wt t_0+ts,s),v(\wt t_0+ts)\right),\\
\bm b(t,s)&=\bm B^v_{(\wt t_0+ts)\wt t_0}\left[\wt{\bm b}_0(s)\right]+t\bm f\left(\bm x(\wt t_0+ts,s),v(\wt t_0+st)\right),
\end{align*}
where $t\in[c,d]$ and $s\in[0,\theta]$. We shall prove that this pair is an $\bm A$-homotopy.

If this is the case, then clearly \eqref{eqn:needle3} follows form Lemma \ref{lem:htp} since the initial-point $\bm A$-homotopy is 
\begin{align*}
\bm b(c,s)=&\bm B^v_{(\wt t_0+cs)\wt t_0}\bm B^v_{\wt t_0(\wt t_0+sc)}\left[\wt{\bm b}_0(s)-c\bm f\left(\bm x(\wt t_0+sc),v(\wt t_0+sc)\right)\right]+\\
\phantom{=}&+c\bm f\left(\bm x(\wt t_0+sc),v(\wt t_0+sc)\right)=\wt {\bm b}_0(s),
\intertext{the final-point $\bm A$-homotopy is}
\bm d_2^{c,d}(s):=&\bm b(d,s)=\bm B^v_{(\wt t_0+ds)\wt t_0}\left[\wt{\bm b}_0(s)-c\bm f\left(\bm x(\wt t_0+sc),v(\wt t_0+sc)\right)\right]+\\
\phantom{=}&+d\bm f\left(\bm x(\wt t_0+sd),v(\wt t_0+sd)\right), 
\end{align*}
and, by Lemma \ref{lem:reparam}, $\left[\bm a(t,s)\right]_{t\in[c,d]}=\left[\bm f(\wt t_0+t,s),v(\wt t_0+t)\right]_{t\in[cs,ds]}$.
Evaluating the formula for $\bm d_2^{c,d}(s)$ at $s=0$ we get \eqref{eqn:needle3a}. 

Finally, $\bm d_2^{c,d}(s)$ is uniformly regular if $\wt{\bm b}_0(s)$ is and $\bm t_0$ is a regular point of $v(t)$. Indeed, we can use the results from Appendix \ref{sec:meas}. The point is to observe that $\bm d_2^{c,d}(s)$ is obtained from measurable maps $v(\wt t_0+s)$ and $\bm b_0(s)$ regular at $s=0$ by operations described in Propositions \ref{prop:ur_cont}--\ref{prop:ur_multiplication} and Lemmas \ref{lem:ur_rescal}--\ref{lem:ur_comp} which preserve the uniform regularity. One has also to use the fact that $\bm B_{t_1 t_0}^v(\bm b_0)$, $f(x,u)$, and $x(t,s)$ are continuous maps (cf. Remark \ref{rem:B_htp}). 

Now it remains to check that $\bm a(t,s)$ and $\bm b(t,s)$ are indeed an $\bm A$-homotopy. Let us calculate the WT-derivatives:
\begin{align*}
\pa_s\bm a^i(t,s)\overset{\text{WT}}=&\bm f^i\left(\bm x(\wt t_0+ts,s),v(\wt t_0+ts)\right)+ts\pa_{\ol t}\bm f^i\left(\bm x(\ol t,s),v(\ol t)\right)|_{\ol t=\wt t_0+ts}+\\
&+s\pa_{s}\bm f^i\left(\bm x(\ol t,s),v(\ol t)\right)|_{\ol t=\wt t_0+ts},
\intertext{and}
\pa_t\bm b^i(t,s)\overset{\text{WT}}=&\pa_{\ol t}\bm B^v_{\ol t\wt t_0}\left[\wh{\bm b}_0(s)\right]^i\Big|_{\ol t=\wt t_0+ts}+\bm f^i\left(\bm x(\wt t_0+ts,s),v(\wt t_0+ts)\right)+\\
&+ts\pa_{\ol t}\bm f^i\left(\bm x(\ol t,s),v(\ol t)\right)|_{\ol t=\wt t_0+ts}.
\end{align*}
Now, since $\bm f\left(\bm x(\ol t,s),v(\ol t)\right)$ and 
$\bm B^v_{\ol 
t\wt t_0}\left[\wh{\bm b}_0(s)\right]$ is an $\bm A$-homotopy (cf. Lemma \ref{lem:sub_lem1}), we have
$$\pa_{\ol t}\bm B^v_{\ol t\wt t_0}\left[\wh{\bm b}_0(s)\right]^i-\pa_{s}\bm f^i\left(\bm x(\ol t,s),v(\ol t)\right)\overset{\text{WT}}=c^i_{jk}(\bm x(\ol t,s))\bm B^v_{\ol t\wt t_0}\left[\wh{\bm b}_0(s)\right]^j\bm f^k\left(\bm x(\ol t,s),v(\ol t)\right).$$
Consequently, 
\begin{align*}
\pa_t \bm b^i(t,s)-\pa_s\bm a(t,s)\overset{\text{WT}}=&s\left[\pa_{\ol t}\bm B^v_{\ol t\wt t_0}\left[\wh{\bm b}_0(s)\right]^i-\pa_{s}\bm f^i\left(\bm x(\ol t,s),v(\ol t)\right)\right]\Big|_{\ol t=\wt t_0+ts}\\\overset{\text{WT}}=&
sc^i_{jk}(\bm x(\ol t,s))\bm B^v_{\ol t\wt t_0}\left[\wh{\bm b}_0(s)\right]^j\bm f^k\left(\bm x(\ol t,s),v(\ol t)\right)\big|_{\ol t=\wt t_0+ts}\\
=&c^i_{jk}(x(\wt t_0+ts,s))\bm b^j(t,s)\bm a^k(t,s).
\end{align*} \end{proof}

No we return to the inductive proof of Theorem \ref{thm:1st_main}. We will prove first that the assertion is true for $t_1=\tau$. Our argument will be inductive w.r.t. $k$---the number of switches in the symbol $\A=(\tau_i, v_i,\tau,\delta t_i, \delta t)_{i=1,\hdots,k}$.

\underline{Step 1, $k=0$.} We start with $k=0$. This means that $u^\A_s(t)=u(t)$. For \underline{$\del t=0$} we simply have $\bm a(t,s)=\bm f(\bm x(t,s,p),u(t))$ with $t\in[t_0,\tau]$, where $\bm x(t_0,s,p)=\bm x_0(s,p)$. Now we can use the Lemma \ref{lem:sub_lem1} taking $\wt t_0=t_0$, $\wt t_1=\tau$, $\wt p=p$, $\wt {\bm b}_0(s,\wt p)=\bm b_0(s,p)$, and $v(t)=u(t)$ to get the assertion. 

If \underline{$\del t\neq 0$}, things are a little more complicated. We have $\bm a(t,s,p)=\bm f(\bm x(t,s,p),u(t))$ where $t\in[t_0,\tau+s\del t]$ and $\bm x(t_0,s,p)=\bm x_0(s,p)$. We can decompose
\begin{equation}\label{eqn:needle2}\begin{split}
&\left[\bm f(\bm x(t,\eps,p),u(t))\right]_{t\in[t_0,\tau+\eps\del t]}\\&=\left[\bm f(\bm x(t,\eps,p),u(t))\right]_{t\in[t_0,\tau]}\cdot\left[\bm f(\bm x(t,\eps,p),u(t))\right]_{t\in[\tau,\tau+\eps\del t]}.\end{split}
\end{equation}
Now using the assertion for $\del t=0$ we get
\begin{equation}\label{eqn:needle4}
[\bm b_0(s, p)]_{s\in[0,\eps]}[\bm{f}(\bm x(t,\eps, p),u(t))]_{t\in[ t_0,\tau]}=[\bm
{f}(\bm x(t,0,p),u(t))]_{t\in[t_0,\tau]}[\bm{d}_1^{ p}(s)]_{s\in[0,\eps]},\end{equation} 
where $\bm d^p_1(s)$ is uniformly regular w.r.t. $p$ at $s=0$, and $\bm d^p_1(0)=\bm B_{\tau t_0}[\bm b_0(0,p)]$. 
Next, using Lemma \ref{lem:sub_lem2} for $\wt t_0=\tau$, $c=0$, $d=\del t$, $\wt p=p$, $\wt {\bm b}_0(s,\wt p)=\bm d_1^p(s)$, and $v(t)=u(t)$, we get
\begin{equation}\label{eqn:needle5}
[\bm d_1^p(s)]_{s\in[0,\eps]}[\bm{f}(\bm x(\tau+t,\eps,\wt p),v(\tau+t))]_{t\in[0,\eps\del t]}=[\bm{d}_2^{ p,\del t}(s)]_{s\in[0,\eps]},\end{equation} 
where $\bm d_2^{p,\del t}(s)$ is uniformly regular w.r.t. $p$, and $\del t$ at $s=0$ and $$\bm d_2^{p,\del t}(0)=\bm d_1^p(0)+\del t \bm f(\bm x(\tau,0,p),u(\tau)).$$ 
Multiplying \eqref{eqn:needle4} by $[\bm{f}(\bm x(\tau+t,\eps,\wt p),v(\tau+t))]_{t\in[0,\eps\del t]}$, using \eqref{eqn:needle2} and \eqref{eqn:needle5}, and taking $\bm d^{\A,p}(s):=\bm d_2^{p,\del t}(s)$, we get the assertion.

\ul{Step 2.} Assume that the assertion holds for all $l<k$. Consider a symbol $\A=(\tau_i,v_i,\tau,\del t_i, \del t)_{i=1,\hdots, k}$. We will distinguish the following two situations:

\noindent\ul{Situation 2.A.} Not all $\tau_i$ are equal. In particular, $$t_0<\tau_1\leq\hdots\leq \tau_l<\tau_{l+1}\leq\hdots\leq\tau_k\leq\tau$$ for some $l<k$. We can now use the inductive assumption for a symbol $\A_1=(\tau_i,v_i,\tau=\tau_l,\del t_i, \del t=0)_{i=1,\hdots, l}$ to get 
\begin{equation}\label{eqn:needle6}\begin{split}
&[\bm b_0(s,p)]_{s\in[0,\eps]}[\bm{f}(\bm x(t,\eps,p),u^{\A_1}_\eps(t))]_{t\in[t_0,\tau_l]}\\&
=[\bm{f}(\bm x(t,0,\eps),u(t))]_{t\in[t_0,\tau_l]}[\bm{d}^{\A_1,p}(s)]_{s\in[0,\eps]},\end{split}\end{equation} 
where $\bm d^{\A_1,p}(s)$ is uniformly regular w.r.t. $p$, $\del t_1,\hdots,\del t_l$ at $s=0$, and
$$\bm d^{\A_1,p}(0)=\bm{B}_{\tau_l t_0}(\bm
b_0(0,p))+\sum_{i=1}^l
\bm{B}_{\tau_l\tau_i}\Big[\bm{f}(\bm{x}(\tau_i),v_i)-\bm{f}(\bm{x}(\tau_i),u(\tau_i))\Big]\del
t_i.$$
Using the inductive assumption for $\A_2=(\tau_i,v_i,\tau,\del t_i, \del t)_{i=l+1,\hdots, k}$ with $t_0=\tau_l$ and $\bm b_0(s, p,\del t_1,\hdots,\del t_l)=\bm d_1^{\A_1,p}(s)$, we get  
\begin{equation}\label{eqn:needle7}\begin{split}
&[\bm d^{\A_1,p}(s)]_{s\in[0,\eps]}[\bm{f}(\bm x(t,\eps,p),u^{\A}_\eps(t))]_{t\in[\tau_l,\tau]}\\
&=[\bm
{f}(\bm x(t,0,\eps),u(t))]_{t\in[\tau_l,\tau]} [\bm{d}^{\A_1,\A_2,p}(s)]_{s\in[0,\eps]},\end{split}\end{equation} 
where $\bm d^{\A_1,\A_2,p}(s)$ is uniformly regular w.r.t. $p$, $\del t_i$ and $\del t$ at $s=0$, and
\begin{align*}\bm d^{\A_1,\A_2, p}(0)&=\bm{B}_{\tau\tau_l}(\bm
d^{\A_1,p}(0))+\sum_{i=l}^k
\bm{B}_{\tau\tau_i}\Big[\bm{f}(\bm{x}(\tau_i),v_i)-\bm{f}(\bm{x}(\tau_i),u(\tau_i))\Big]\del
t_i\\
&=\bm{B}_{\tau t_0}(\bm
b_0(0,p))+\sum_{i=1}^k
\bm{B}_{\tau\tau_i}\Big[\bm{f}(\bm{x}(\tau_i),v_i)-\bm{f}(\bm{x}(\tau_i),u(\tau_i))\Big]\del
t_i.\end{align*}
Multiplying \eqref{eqn:needle6} by $[\bm{f}(\bm x(t,\eps,p),u^{\A}_\eps(t))]_{t\in[\tau_l,\tau]}$, using \eqref{eqn:needle7}, and taking $\bm d^{\A,p}(s):=\bm d^{\A_1,\A_2,p}(s)$, we get the assertion.

\noindent\ul{Situation 2.B.} If all $\tau_i$ are equal then either $$\tau_1=\hdots\tau_k=\tau\quad \text{or}\quad \tau_1=\hdots\tau_k<\tau.$$

\noindent\ul{2.B.1.} In the first case using the assertion for $\A_1=(\tau_i,v_i,\del t_i,\tau, \del t-\del t_k)_{i=1,2,\hdots,k-1}$ we get 
\begin{equation}\label{eqn:przyp_B1}\begin{split}
&[\bm b_0(t,s)]_{s\in[0,\eps]}[\bm f(\bm x(t,\eps,p),u^{\A_1}_\eps(t))]_{t\in[t_0,\tau+(\del t-\del t_k)\eps]}\\&
=[\bm f(\bm x(t,0,p),u(t))]_{t\in[t_0,\tau]}[\bm d_1^{\A_1,p}(s)]_{s\in[0,\eps]},\end{split}
\end{equation} 
where $\bm d_1^{\A_1,p}(s)$ is uniformly regular w.r.t. $p$, $\del t_1,\hdots,\del t_{k-1}$, $\del t-\del t_k$ at $s=0$ and
$$\bm d_1^{\A_1,p}(0)=\bm B_{\tau t_0}(\bm b_0(0,p))+(\del t-\del t_k)\bm f(x(\tau),u(\tau))+\sum_{i=1}^{k+1}\bm B_{\tau\tau_i}\left[\bm f(\bm x(\tau_i),v_i)-\bm f(\bm x(\tau_i),u(\tau_i))\right]\del t_i.$$
Now using Lemma \ref{lem:sub_lem2} for $\wt t_0=\tau$, $c=\del t-\del t_k$, $\wt v(t)=v_k$, $\wt p=(p,\del t_1,\hdots,\del t_{k-1},\del t-\del t_k)$, and $\wt{\bm b}_0(s,\wt p)=\bm d_1^{\A_1,p}(s)$ we get
\begin{equation}\label{eqn:needle8}
\left[\bm d_1^{\A_1,p}(s)\right]_{s\in[0,\eps]}\left[\bm f(\bm x(\tau+t,\eps,p),v_k)\right]_{t\in[\eps(\del t-\del t_{k-1}),\eps\del t]}=\left[\bm d_2^{\A_1,\del t,\del t_k,p}(s)\right]_{s\in[0,\eps]},
\end{equation}
where $\bm d_2^{\A_1,\del t,\del t_k,p}(s)$ is uniformly regular w.r.t. $p$, $\del t_1,\hdots,\del t_k$, $\del t$ at $s=0$ and
$$\bm d_2^{\A_1,\del t,\del t_k,p}(0)=\bm d_1^{\A_1,p}(0)-\del t_k\bm f(\bm x(\tau,0,p),u(\tau)).$$
Again multiplying \eqref{eqn:needle7} by $\left[\bm f(\bm x(\tau+t,\eps,p),v_k)\right]_{t\in[\eps(\del t-\del t_{k-1}),\eps\del t]}$, using \eqref{eqn:needle8} and taking $\bm d^{\A,p}(s)=\bm d_2^{\A_1,\del t,\del t_k,p}(s)$ we get the assertion.

\noindent\ul{2.B.1.} If $\tau_1=\hdots=\tau_k<\tau$ we can use the result from 2.B.1 for a symbol $\A_1=(\tau_i,v_i,\del t_i,\tau=\tau_k,\del t=0)_{i=1,\hdots,k-1}$ and then use the inductive assumption for $[t_0,\tau]=[\tau_k,\tau]$ and $\A_2=(\tau,\del t)$ on $[t_0=\tau_k,\tau]$ in essentially the same way as in the case A. The inductive argument is now complete.

Finally, to obtain the assertion for $t_1$ not $\tau$, one has just to proceed as in the step 1 with $\del t=0$ and  use Lemma \ref{lem:sub_lem1} again, taking $\wt t_0=\tau$, $\wt t_1=t_1$, $\wt{\bm b}_0(s,\wt p)$ to be the final-point $\bm A$-homotopy $\bm d^{\A,p}(s)$ derived for $t_1=\tau$, and the control $v(t)=u(t)=u^{\A}_s(t+s\del t)$. \end{proof}

\subsection{The set of infinitesimal variations $\bm K^u_\tau$} 

\begin{remark}\label{rem:interpretation_K}
Observe that choosing $\bm b_0(s)\equiv \theta_{\bm x_0}$ in Theorem \ref{thm:1st_main} we obtain an admissible path  $\bm d^\A(s)$, regular at $s=0$, defined for $0\leq s<\theta$, and satisfying
\begin{equation}
\label{eqn:var_htp}[\bm{a}(t,\eps)]_{t\in[t_0,t_1+\eps\del
t]}=[\bm {a}(t)]_{t\in[t_0,t_1]}[\bm{d}^\A(s)]_{s\in[0,\eps]},
\end{equation}
for $\eps\leq\theta$.

We define the set $\bm K_\tau^u$ consisting of elements of the form $\bm d^\A(0)$, where $\A=(\tau_i, v_i,\tau,\delta t_i, \delta t)_{i=1,\hdots,k}$ are symbols with $\tau$ fixed: 
\begin{align*} 
\bm K^u_\tau&:=\left\{\bm B_{t_1\tau}[\bm{f}(\bm{x}(\tau),u(\tau))]\del
t+\sum_{i=1}^k
\bm{B}_{t_1\tau_i}\Big[\bm{f}(\bm{x}(\tau_i),v_i)\right.\\&-\left.\bm{f}(\B{x}(\tau_i),u(\tau_i))\Big]\del
t_i:(\tau_i, v_i,\tau, \delta t_i, \delta
t)_{i=1,\dots,k}\text{ is a symbol}\right\}\subset\bm A_{\bm x(t_1)}.
\end{align*}
We will call $\bm K_\tau^u$ the
\emph{set of infinitesimal variations of the trajectory $\bm f(\bm
x(t),u(t))$ associated with the regular $\tau\in (t_0,t_1)$}.

The set $\bm{K}^u_\tau$ can be interpreted as the set of all generalised 
directions in $\bm{A}_{\bm{x}(t_1)}$ in which one can move the final base-point $\bm{x}(t_1)$ by performing needle variations, associated with symbols $\A=(\tau_i, v_i,\tau,\delta t_i, \delta t)_{i=1,\hdots,k}$ with fixed $\tau$ and trivial initial base-point variations $\bm x_0(s)\equiv\bm x_0$. 
\end{remark}

The geometry of $\bm{K}^u_\tau$ will be an object of our main interests in Section \ref{sec:proof}. Now let us note the following property

\begin{lemma}\label{lem:K_convex}
The set $\bm{K}^u_\tau$ is a convex cone in
$\bm{A}_{\bm{x}(t_1)}$.
\end{lemma}
\begin{proof}  Take symbols $\A=(\tau_i, v_i,\tau, \del t_i, \del
t)_{i=1,\hdots,k}$, $\A^{'}=(\tau^{'}_i, v^{'}_i,\tau, \del
t^{'}_i, \del t^{'})_{i=1,\hdots,k^{'}}$ and numbers
$\nu,\nu^{'}\geq 0$. We have to find a symbol $\mathfrak{v}$ such
that
$$\bm{d}^{\mathfrak{b}}(0)=\nu\bm{d}^{\A}(0)+\nu^{'}\bm{d}^{\A^{'}}(0).$$
Due to Remark \ref{rem:var_dt=0}, we may change the symbol by
adding $(\tau_i,v_i,\delta t_i=0)$ without changing the variation
$u^\A_\eps$. As we see from the form of \eqref{eqn:d_at_0},
such an addition will not change $\bm{d}^\A(0)$. Consequently, we
may assume that $k=k^{'}$, $\tau_i=\tau^{'}_i$, $v_i=v^{'}_i$, and
the symbols $\A$ and $\A^{'}$ differ only by $\del t_i$ and $\del
t$. Now consider the symbol $\mathfrak{v}=(\tau_i, v_i,\tau,
\nu\del t_i+\nu^{'}\del t^{'}_i, \nu\del
t+\nu^{'}\del^{'})_{i=1,\hdots,k}$. The formula
\eqref{eqn:d_at_0} (for $\bm b_0(0,p)=0$) is linear with respect
to $\del t_i$ and $\del t$, hence
$$\bm{d}^{\mathfrak{v}}(0)=\nu\bm{d}^{\A}(0)+\nu^{'}\bm{d}^{\A^{'}}(0).$$
\end{proof}

\section{Important facts}\label{sec:lemma}

In this section we prove a technical result about $E$-homotopies -- Lemma \ref{lem:htp_of_a_r}, which will be crucial in the proof of Theorem \ref{thm:pmp}. To discuss briefly the result, given a family of smooth curves $x_{\vec r}:I\ra\R^m$, parameterized by $\vec r\in B^m(0,1)\subset\R^m$, which emerges from a single point $x_{\vec r}(0)=0$ and points into every direction $\dot x_{\vec r}(0)=\vec r$, it is quite obvious that, for every $t>0$ small enough, there exists a curve $x_{\vec {r_0}}$ from this family which reaches $0$ at time $t$. A similar result holds for families of admissible curves on a skew-algebroid $E$. Any such family which is sufficiently regular and emerges from a single point into every direction in $E$ will realise a zero homotopy class. This is the assertion of Lemma \ref{lem:htp_of_a_r}.

This result seems to be quite natural and it is indeed, if such an algebroid is (locally) integrable. In this case $E$-homotopy classes can be represented by points on a finite-dimensional manifold. However, if $E$ is not integrable, $E$-homotopy is just a relation in the space of bounded measurable curves. Therefore to prove the results we have to pass through the Banach space setting. The main idea in the proof is to semi-parametrise the $E$-homotopy classes by a finite dimensional-space and reduce the problem to a finite-dimensional topological problem. By a semi-parametrisation we mean an epimorphism from a finite-dimensional space to the space of $E$-homotopy classes. 

\subsection{Local coordinates} Since we are going to work in a Banach space setting it is convenient to introduce local coordinates on an algebroid $E$. Consider coordinates $(x^a,y^i)\in U\times\R^m\subset\R^n\times\R^m$ trivialising the bundle $\tau:E\ra M$ around a point $p\in M$. We may assume that $p$ corresponds to $0\in U$.  As usual, we will denote
the structural functions of $E$ in these coordinates by
$\rho^a_i(x)$ and $c^i_{jk}(x)$. Since these functions are smooth,
we can assume (after restricting ourselves to a compact
neighborhood  $\ol V\ni 0$ in $\R^n$) that they are bounded (by
numbers $C_\rho$ and $C_c$, respectively) and Lipschitz w.r.t. $x^a$ (with constants $L_\rho$ and $L_c$,
respectively). It will be convenient to think of $V\times\R^m$
with those functions as of a (local) AL algebroid. Observe that
every bounded measurable $E$-path with the base initial-point
$p$ is represented by a pair of paths
$(x(t),a(t))\in\R^n\times\R^m$, where $a(t)$ is bounded measurable
and $x(t)$ is an AC-solution of the ODE
$$\begin{cases}
\dot x^b(t)=\rho^b_i(x(t))a^i(t),\\
x^b(0)=0.
\end{cases}
$$
As we see, $x(t)$ is determined entirely by $a(t)$. We can thus
identify the space $\mathcal{ADM}_p(I,E)$ of bounded
measurable admissible paths originated at $p$ with  the
space $\BM(I,\R^m)$ of bounded measurable maps $a:I\ra\R^m$. We
will consequently speak of algebroid homotopy classes in
$\BM(I,\R^m)$. Note that $\BM(I,\R^m)$, equipped with the
$L^1$-norm, is a Banach space. We will denote this norm simply by $\|\cdot\|$. The same symbol will be also used for $L^1$-norm in $\R^m$.  In our considerations we will understand a product of Banach spaces $(\mathcal{B}_1, \|\cdot\|_1)$ and  $(\mathcal{B}_2, \|\cdot\|_2)$ as a space  $\mathcal{B}_1\times\mathcal{B}_2$ equipped with the norm $\|\cdot\|= \|\cdot\|_1+\|\cdot\|_2$.

\subsection{Technical lemma}

\begin{lemma}\label{lem:htp_of_a_r}
Let $a_{\vec{r}}(\cdot)\in\BM(I,R^m)$, where $\vec r\in  B^m(0,1)$, be
a family of $E$-paths uniformly regular at $t=0$ w.r.t. $\vec r$
and such that  $a_{\vec{r}}(0)=\vec{r}$. Then there exists a number $\eta>0$ with the following property. For every $0<\eps<\eta$ there exists a 
vector $\vec r_0$ such that the curve
$a_{\vec r_0}(t)$, after restricting to the interval $[0,\eps]$,
is null-$E$-homotopic:
$$\big[a_{\vec{r}_0}(t)\big]_{t\in[0,\eps]}=\big[0\big].$$
\end{lemma}

Let us briefly sketch the strategy of the proof. 
Denote by $c_{\vec r}$ a constant path $c_{\vec r}(s)=\vec r$ in $\R^m$. We will construct
a continuous and invertible (local) map of Banach spaces
$\Phi:\R^m\times \mathcal{B}\lra \BM(I,\R^m)$ (the space
$\mathcal{B}$ will be specified later) which will have an
additional property that the homotopy class of the image is
determined by the first factor only
$$\left[\Phi(\vec r,d)\right]=\left[c_{\vec r}\right].$$

In such a way we will realise our idea from the introduction to this section --- $\R^m$ will semi-parametrise all local $E$-homotopy classes of $\mathcal{ADM}_p(I,E)$. Next, using
the map $a:\vec r\mapsto a_{\vec r}$, we will construct a
continuous map of finite-dimensional spaces 
$$\R^m\supset B^m(0,1)\overset
a{\lra}\BM([0,1],\R^m)\overset{\Phi^{-1}}\lra\R^m\times
\mathcal{B}\overset{\operatorname{pr}_1}\lra\R^m.$$ 
A topological argument will prove that $0$ lies in the image of this map, hence 
$$\left[a_{\vec{r_0}}\right]=\left[c_0\right]=\left[0\right]\quad \text{for some $\vec{r_0}$}.$$

\begin{proof}
Consider an $E$-path with a constant $\R^m$-part $a(t)=\vec r$,
where $\vec r$ is a fixed element in $\R^m$, $t\in[0,1]$. The
associated base path $x(t)\in\R^m$ is the solution of
\begin{equation}\label{eqn:a_x}
\left\{ \begin{aligned}
\pa_t x^b(t)&=\rho^b_i(x(t))a^i(t)=\rho^b_i(x(t))r^i,\\
x^b(0)&=0.
\end{aligned}\right.
\end{equation}
Clearly, if $\|\vec r\|$ is small enough, the solution of this
equation exists for $t\in[0,1]$ and is contained entirely in $\ol
V\subset \R^n$. Now for $a(t)$ and $x(t)$ as above and fixed paths
$d\in\BM(I,\R^m)$, $b\in\mathcal{AC}(I,\R^m)$, consider the
following system of differential equations:
\begin{equation}\label{eqn:phi_1}
\left\{ \begin{aligned}
\pa_s a^i(t,s)&=d^i(t)+c^i_{jk}(x(t,s))a^j(t,s)b^k(t),\\
a^i(t,0)&=a^i(t)=r^i,\\
\pa_s x^b(t,s)&=\rho^b_i(x(t,s))b^i(t),\\
x^b(t,0)&=x^b(t).
\end{aligned}\right.
\end{equation}
The existence and regularity of the solutions of \eqref{eqn:phi_1} can be discussed using the theory developed in Subsection \ref{ssec:ode}. Let us concentrate first on the equation for $x(t,s)$. The right hand-side of this equation is AC in $t$ and Lipschitz in $x$, the initial value depends AC on a parameter $t$, and hence, by the standard theory of ODEs, the solution $x(t,s)$ is defined locally and is AC w.r.t. both variables. By shrinking the norm $\|b\|_{\sup}$ we may change the Lipschitz constant in the defining equation. Consequently, for $\|b\|_{\sup}$ (and $\|\vec r\|$) small enough, the solution $x(t,s)$ is defined for all $t,s\in[0,1]$ and entirely contained in $\ol V$.

Now the right hand-side of the first equation in \eqref{eqn:phi_1} is locally Lipschitz w.r.t. $a$ and bounded measurable w.r.t. the parameter $t$. The initial value $a^i(t,0)$ depends continuously on $t$, hence, by Theorem \ref{thm:param}, the solution $a(t,s)$ locally exists, is AC w.r.t. $s$, and is bounded measurable w.r.t. $t$. Again, shrinking $\|b\|_{\sup}$ makes the Lipschitz constant smaller, hence for $\|b\|_{\sup}$ small enough $a(t,s)$ is defined for all $t,s\in[0,1]$.

Now let us consider \eqref{eqn:phi_1} with $b(t)=\int_0^td(s)\dd
s$, where $d$ is chosen in such a way, that $b(0)=b(1)=0$. We have
\begin{equation}\label{eqn:phi_2}
\left\{ \begin{aligned}
\pa_s a^i(t,s)&=\pa_tb^i(t)+c^i_{jk}(x(t,s))a^j(t,s)b^k(t),\\
\pa_s x^b(t,s)&=\rho^b_i(x(t,s))b^i(t).
\end{aligned}\right.
\end{equation}
We recognise equations \eqref{eqn:htp_smooth} for $E$-homotopy.
Indeed, in such a situation $a(t,s)$ and $b(t,s)=b(t)$ form an
$E$-homotopy with fixed end-points (since $b(0)=b(1)=0$).
Consequently, the homotopy classes of $a(t,0)=c_{\vec r}(t)$ and
$a(t,1)$ are equal. Since $\|b\|_{\sup}\leq\|d\|$, for
$\|d\|$ and $\|\vec r\|$ small enough, this homotopy is defined
for all $t,s\in[0,1]$. For $a(t)$ and $d(t)$ as above we define
$$\Phi(\vec r,d)(t):=a(t,1).$$
This is a (local) map of Banach spaces
$$\Phi:\R^m\times \BM_0(I,\R^m)\supset W_0\lra\BM(I,\R^m),$$
where $\BM_0(I,\R^m)=\{d\in \BM(I,\R^m):\int_0^1d(s)\dd s=0\}$
is a Banach subspace of $\BM(I,\R^m)$ and $W_0$ is some open neighbourhood of the point $(0,0)$. We shall now prove the following:
\begin{description}
    \item[(A)]\label{cond:A} $\Phi$ maps $(\vec{r},0)$ into a constant path $c_{\vec{r}}\in \BM(I,\R^m)$.
     \item[(B)]\label{cond:B} $\Phi$ is a continuous map of Banach spaces.
     \item[(C)]\label{cond:C} The $E$-homotopy class of the curve $\Phi(\vec{r},d)$ is determined by $\vec{r}$; that is,
$$\big[\Phi(\vec{r},d)\big]=\big[\Phi(\vec{r},0)\big]\overset{(A)}{=}\big[c_{\vec{r}}\big].$$
     \item[(D)]\label{cond:D} The map $\widetilde{\Phi}(\vec{r},d):=\Phi(\vec{r},d)-(c_{\vec{r}}+d)$ is Lipschitz with constant $\frac{1}{6}$.
     \item[(E)]\label{cond:E} The map $\Phi$ posses a continuous inverse $\Phi^{-1}$ defined on some open neighbourhood $V_0\ni 0$ in $\BM(I,\R^m)$. Moreover, $\Phi^{-1}$ is Lipschitz with constant $6$.
\end{description}

Property (\textbf{C}) is clear from the construction of $\Phi$, as
$a(t,1)=\Phi(\vec r,d)(t)$ and $a(t,0)=c_{\vec r}(t)$ are
$E$-homotopic. 

Property (\textbf{A}) is obvious, since $\Phi(\vec
r,0)$ is the solution (taken  at $s=1$) of the differential equation $\pa_sa(t,s)=0$
with the  initial condition $a(t,0)=\vec r$.

Property (\textbf{B}) will follow from (\textbf{D}). Indeed, if
$\wt\Phi$ is Lipschitz, then $\Phi(\vec
r,d)=\wt\Phi(\vec r,d)+c_{\vec r}+d$ is continuous as a sum of
continuous maps.

Assuming (\textbf{D}) again, we will be able to prove
(\textbf{E}). As one might have expected, the existence and the
Lipschitz condition for $\Phi^{-1}$ will be proven essentially in
the same way as in the standard proof of the inverse function
theorem (cf. \cite{Lang}). First, we will establish a pair
of linear isomorphism between Banach spaces
\begin{align*}
\BM(I,\R^m)&\overset{\alpha}{\lra}\R^m\times \BM_0(I,\R^m)\,,\\
a(t)&\longmapsto \left(\int_0^1a(s)\dd s,\ a(t)-\int_0^1a(s)\dd
s\right)\,,
\intertext{and}
\R^m\times \BM_0(I,\R^m)&\overset{\beta}{\lra}\BM(I,\R^m),\\
(\vec{r},d)&\longmapsto c_{\vec{r}}+d.
\end{align*}
It is straightforward to verify that $\alpha$ and $\beta$ are
continuous inverses of each other and that $\alpha$ is Lipschitz
with constant $3$. The map $\Phi$ is defined on some open
neighbourhood $W_0\ni(0,0)$. Take $R$ such that $B(0,2R)\subset
W_0$. The map $\widetilde{\Phi}$ is Lipschitz with constant
$\frac{1}{6}$ and $\alpha$ is Lipschitz with constant $3$; hence
$\alpha\circ\widetilde{\Phi}$ is Lipschitz with constant
$\frac{1}{2}$ and, since it preserves the origin, it maps the ball
$B(0,2R)$ into the ball $B(0,R)$.

Fix now any $a\in \BM(I,\R^m)$ such that $\|a\|<\frac{R}{3}$. We
shall construct a unique element $(\vec r,d)\in B(0,2R)\subset W_0$
satisfying $\Phi(\vec r,d)=a$. Consider a map
$\Phi_a(\vec{r},d):=\alpha(a-\widetilde{\Phi}(\vec{r},d))$. From the Lipschitzity of $\alpha$ and $\wt \Phi$ we deduce that 
$$\|\Phi_a(\vec r,d)\|\leq 3\|a-\wt\Phi(\vec r,d)\|\leq 3\|a\|+3\|\wt\Phi(r,d)\|\leq 3\cdot\frac R3+\frac 12\|(\vec r,d)\|.$$ 
Consequently, $\Phi_a$ maps
the ball $B(0,2R)$ into $B(0,2R)$. Moreover,
\begin{align*}
\|\Phi_a(\vec r,d)-\Phi_a(\vec r^{'},d^{'})\|&=\|\alpha(a-\wt \Phi(\vec r,d))-\alpha(a-\wt\Phi(\vec r^{'},d^{'}))\|\leq\\
&\leq 3\|\wt\Phi(\vec r,d)-\wt\Phi(\vec r^{'},d^{'})\|\leq 3\cdot\frac 16\|(\vec r,d)-(\vec r^{'},d^{'})\|,
\end{align*}
 hence $\Phi_a$ is a contraction. Now, using the Banach fixed point theorem,
we deduce that $\Phi_a$ has a unique fixed point $(\vec{r},d)\in
B(0,2R)$. Consequently,
$$a-\widetilde{\Phi}(\vec{r},d)=\beta\circ\alpha\left(a-\widetilde{\Phi}(\vec{r},d)\right)=\beta\circ\Phi_a(\vec{r},d)=\beta(\vec r,d)=c_{\vec{r}}+d,$$
and hence $a=\wt\Phi(\vec{r},d)+c_{\vec{r}}+d=\Phi(\vec{r},d).$
We have proven the existence of $\Phi^{-1}$.

Take now $a, a^{'}\in \BM(I,\R^m)$, and let
$\Phi^{-1}(a)=(\vec{r},d)$,
$\Phi^{-1}(a^{'})=(\vec{r}^{'},d^{'})$. Using the Lipschitz
condition for $\alpha$ and $\widetilde{\Phi}$ once more, we get
\begin{align*}
&\|(\vec{r},d)-(\vec{r}^{'},d^{'})\|=\|\Phi_a(\vec{r},d)-\Phi_{a^{'}}(\vec{r}^{'},d^{'})\|=
\|\alpha\Big(a-\widetilde{\Phi}(\vec{r},d)\Big)-\alpha\Big(a^{'}-\widetilde{\Phi}(\vec{r}^{'},d^{'})\Big)\|\leq\\
&\leq
3\|a-a^{'}\|+3\|\widetilde{\Phi}(\vec{r},d)-\widetilde{\Phi}(\vec{r}^{'},d^{'})\|\leq
3\|a-a^{'}\|+3\cdot\frac{1}{6}\|(\vec{r},d)-(\vec{r}^{'},d^{'})\|.
\end{align*}
We finish the proof of property (\textbf{E}) concluding that
$$\|\Phi^{-1}(a)-\Phi^{-1}(a^{'})\|=\|(\vec{r},d)-(\vec{r}^{'},d^{'})\|\leq 6\|a-a^{'}\|.$$

We are now left with the proof of (\textbf{D}). This will be done
by introducing several integral estimations. In our calculations
we will, for simplicity, omit the indices (hence $c$ will stand
for $c^i_{jk}$, $a$ for $a^i$, etc.). Take pairs $(r,d)$ and
$(r^{'},d^{'})$ from $\R^m\times\BM_0(I,\R^m)$. Denote by
$x(t,s)$, $a(t,s)$, $a(t)$, $b(t)$ and  $x^{'}(t,s)$,
$a^{'}(t,s)$, $a^{'}(t)$, $b^{'}(t)$, respectively, the objects
defined as in the construction of $\Phi$ for pairs $(\vec r,d)$ and $(\vec r^{'},d^{'})$.
To begin with, observe that, since $b(t)=\int_0^td(s)\dd s$, we
have $|b(t)|\leq\int_0^1|d(s)|\dd s$; hence
$$\|b\|_{\sup}\leq\|d\|.$$
Similarly, $\|b^{'}\|_{\sup}\leq\|d^{'}\|$ and
$\|b-b^{'}\|_{\sup}\leq\|d-d^{'}\|$.

Let us now estimate the difference $|x(t,s)-x^{'}(t,s)|$. Since,
by \eqref{eqn:phi_2},
$x(t,s)=x(t)+\int_0^s\rho(x(t,\sigma))b(t)\dd\sigma$, we have
\begin{align*}
|x(t,s)-x^{'}(t,s)|&=\left|\int_0^s\Big(\rho(x(t,\sigma))b(t)-\rho(x^{'}(t,\sigma))b^{'}(t)\Big)\dd\sigma\right|\leq\\
&\leq\int_0^1\Big|\rho(x(t,s))b(t)-
\rho(x^{'}(t,s))b^{'}(t)\Big|\dd s
\leq\\
&\leq\int_0^1\Big|\Big(\rho(x(t,s))-\rho(x^{'}(t,s))\Big)b(t)\Big|\dd
s+\\&\phantom{=}+
\int_0^1\Big|\rho(x^{'}(t,s))\Big(b(t)-b^{'}(t)\Big)\Big|\dd s\leq\\
&\leq
L_\rho\sup_{(t,s)}|x(t,s)-x^{'}(t,s)|\|d\|+C_\rho\|d-d^{'}\|.
\end{align*}
It follows that
$$
\sup_{(t,s)}|x(t,s)-x^{'}(t,s)|\Big(1-L_\rho \|d\|\Big)\leq
C_\rho\|d-d^{'}\|.$$ For $\|d\|$ sufficiently small (in other
words, after a possible shrinking of $W_0$) we will obtain
\begin{equation}\label{est:x}
\sup_{(t,s)}|x(t,s)-x^{'}(t,s)|\leq 2C_\rho\|d-d^{'}\|.
\end{equation}
Now introduce
\begin{align*}
&\del a(t,s):= a(t,s)-(a(t)+sd(t)) \quad \text{and}\\
&\del a^{'}(t,s):= a^{'}(t,s)-(a^{'}(t)+sd^{'}(t)).
\end{align*}
Note that $\del a(t,1)=\wt\Phi(\vec r,d)$ and $\del a^{'}(t,1)=\wt\Phi(\vec r^{'},d^{'})$. 
From \eqref{eqn:phi_1} we deduce that
$$\del a(t,s)=\int_0^s c(x(t,\sigma))a(t,\sigma)b(t)\dd\sigma=\int_0^s c(x(t,\sigma))\Big(\del a(t,\sigma)+a(t)+\sigma d(t)\Big)b(t)\dd\sigma.$$
Further, for a fixed $s\in[0,1]$,
\begin{align*}
\|\del a(\cdot,s)\|&=\int_0^1|\del a(t,s)|\dd
t\leq\int_0^1\int_0^s\left|c(x(t,\sigma))\Big(\del
a(t,\sigma)+a(t)+
\sigma d(t)\Big)b(t)\right|\dd\sigma\dd t\leq\\
&\leq \int_0^1\int_0^1C_c\big(|\del a(t,s)|+\|\vec r\|+|d(t)|\big)\|b\|_{\sup}\dd s\dd t\leq\\
&\leq C_c\|b\|_{\sup}\left(\sup_s\|\del a(\cdot,s)\|+\|\vec
r\|+\|d\|\right)\leq\\
&\leq C_c\|d\|\left(\sup_s\|\del
a(\cdot,s)\|+\|\vec r\|+\|d\|\right)
\end{align*}
and we conclude that
$$\sup_s\|\del a(\cdot,s)\|(1-C_c\|d\|)\leq C_c(\|\vec r\|+\|d\|).$$
Hence, for $\|d\|$ small enough (after possible shrinking of
$W_0$), we get
\begin{equation}\label{est:del_a}
\sup_s\|\del a(\cdot,s)\|\leq 2C_c(\|\vec r\|+\|d\|).
\end{equation}
Finally,
\begin{align*}
\|\del a&(\cdot,s)-\del a^{'}(\cdot,s)\|=\int_0^1
|\del a(t,s)-\del a^{'}(t,s)|\dd t\leq \\
\leq&\int_0^1\int_0^s\left|c(x(t,s))\Big(\del
a(t,\sigma)+a(t)+\sigma d(t)\Big)b(t)\right.+\\&- \left.c(x^{'}(t,\sigma))\left(\del
a^{'}(t,\sigma)+a^{'}(t)+
\sigma d^{'}(t)\right)b^{'}(t)\right|\dd\sigma\dd t\leq\\
\leq&\int_0^1\int_0^1\left|c(x(t,\sigma))\del a(t,s)b(t)-c(x^{'}(t,s))\del a^{'}(t,s)b^{'}(t)\right|\dd s\dd t+\\
&+\int_0^1\int_0^1\left|c(x(t,s))a(t)b(t)-c(x^{'}(t,s))a^{'}(t)b^{'}(t)\right|\dd s\dd t+\\
&+\int_0^1\int_0^1\left|c(x(t,s))d(t)b(t)-c(x^{'}(t,s))
d^{'}(t)b^{'}(t)\right|\dd s\dd t=:I_1+I_2+I_3
\end{align*}
Now we estimate
\begin{align*}
I_1&\leq\int_0^1\int_0^1\left|c(x(t,s))-c(x^{'}(t,s))\right||\del a(t,s)||b(t)|\dd t\dd s+\\
&\phantom{=}+\int_0^1\int_0^1\left|c(x^{'}(t,s))\right|\left|\del a(t,s)-\del a^{'}(t,s)\right|\left|b(t)\right|\dd t\dd s+\\
&\phantom{=}+\int_0^1\int_0^1\left|c(x^{'}(t,s))\right|\left|\del a^{'}(t,s)\right|\left|b(t)-b^{'}(t)\right|\dd t\dd s\leq\\
&\leq L_c\sup_{(t,s)}|x(t,s)-x^{'}(t,s)|\sup_s\|\del
a(\cdot,s)\|\cdot\|b\|_{\sup}+
C_c\sup_s\|\del a(\cdot,s)-\del a^{'}(\cdot,s)\|\cdot\|b\|_{\sup}+\\
&\phantom{=}+C_c\sup_s\|\del a^{'}(\cdot,s)\|\cdot\|b-b^{'}\|_{\sup}\,.
\end{align*}
Finally, using \eqref{est:x}, \eqref{est:del_a} and $\|b\|_{\sup}\leq\|d\|$, we get
\begin{align*}
I_1&\leq L_c 2C_\rho\|d-d^{'}\|2C_c(\|\vec r\|+\|d\|)\|d\|+
C_c\sup_s\|\del a(\cdot,s)-\del a^{'}(\cdot,s)\|\cdot\|d\|
+C_c2C_c(\|\vec r^{'}\|\\&+\|d^{'}\|)\|d-d^{'}\|
=C_c\|d\|\sup_s\|\del a(\cdot,s)-\del
a^{'}(\cdot,s)\|+\|d-d^{'}\|\cdot
F_1(\|r\|,\|d\|,\|r^{'}\|,\|d^{'}\|),
\end{align*}
where $F_1$ converges to 0 when its arguments do. Similar
estimations for $I_2$ and $I_3$ will give
\begin{align*}
&I_2\leq \left(\|r-r^{'}\|+\|d-d^{'}\|\right)\cdot F_2(\|r\|,\|d\|,\|r^{'}\|,\|d^{'}\|),\\
&I_3\leq \left(\|r-r^{'}\|+\|d-d^{'}\|\right)\cdot
F_3(\|r\|,\|d\|,\|r^{'}\|,\|d^{'}\|),
\end{align*}
where $F_2$ and $F_3$ behave as $F_1$. Putting together the
partial results, we would get
$$\sup_s\|\del a(\cdot,s)-\del a^{'}(\cdot,s)\|(1-2C_c\|d\|)\leq  \left(\|r-r^{'}\|+\|d-d^{'}\|\right)\cdot F\left(\|r\|,\|d\|,\|r^{'}\|,\|d^{'}\|\right),$$
where $F$ converges to 0 when its arguments do. As
$\widetilde{\Phi}(r,d)=\del a(t,1)$ and
$\widetilde{\Phi}(r^{'},d^{'})=\del a^{'}(t,1)$, for $W_0$ small
enough, $\wt\Phi$ is Lipschitz with constant $\frac 16$. That
proves property (\textbf{D}).

Now using properties (\textbf{A})--(\textbf{E}) of $\Phi$ we will
make the final step of the proof of Lemma \ref{lem:htp_of_a_r}.
The family $a_{\vec r}(t)$ is uniformly regular w.r.t. $\vec r\in
\ol B^m(0,1)$ at $t=0$. The family $c_{\vec r}(t)$ has the same
properties, so $a_{\vec r}(t)-c_{\vec r}(t)$ is also uniformly
regular (cf. Proposition \ref{prop:ur_sum}) and
$$\int_0^t|a_{\vec r}(s)-c_{\vec r}(s)|\dd s=t|a_{\vec r}(0)-c_{\vec r}(0)|+f(t,\vec r)=t\cdot 0+f(t,\vec r),$$
where $\frac 1tf(t,\vec r)$ converges uniformly to $0$ as $t\to
0$. Hence, there exists a number $\eta>0$ such that
$$\int_0^\eps|a_{\vec r}(s)-c_{\vec r}(s)|\dd s<\frac 1{12}\eps,$$
for every $0\leq\eps\leq\eta$ and $\vec{r}\in\ol B^m(0,1)$. Reparametrising the
paths by the rule
$$\wt a_{\vec r}(t):=\eps a_{\vec r}(\eps t),$$
we will obtain another uniformly regular family of paths satisfying $\wt a_{\vec r}(0)=\eps \vec r$.

By the uniform regularity of $\wt a_{\vec r}$, the map $\wt a:\vec r\mapsto \wt a_{\vec r}$ is a continuous map from $B^m(0,1)$ to $\BM(I,\R^m)$ with $L_1$-topology. For $\eta$ small enough $\wt a$ takes values in $\Phi(W_0)$. 
Composing $\wt a$ with
$\Phi^{-1}$ we will obtain a continuous map
$$\ol B^m(0,1)\xra{(\psi,\phi)}\R^m\times \BM_0(I,\R^m).$$

Observe that 
$$\|\wt a_{\vec r}-c_{\eps\vec r}\|=\int_0^1|\eps a_{\vec
r}(t\eps)-\eps c_{\vec r}(t\eps)|\dd t=\int_0^\eps| a_{\vec
r}(t)- c_{\vec r}(t)|\dd t\leq\frac 1{12}\eps.$$
Using this, the Lipschitz condition for  $\Phi^{-1}$ (property (\textbf{E})), and the fact that $\Phi^{-1}(c_{\vec r})=(\vec r,0)$ (property (\textbf{A})), we obtain
$$\|\frac 1\eps\psi(\vec{r})-\vec{r}\|\leq
\frac 1\eps\|\psi(\vec{r})-\eps\vec{r}\|+\frac 1\eps\|\phi(\vec{r})-0\|=
\frac 1\eps\|\Phi^{-1}(\widetilde{a}_{\vec{r}})-\Phi^{-1}(c_{\eps\vec{r}})\|\leq\frac
6\eps\|\widetilde{a}_{\vec{r}}-c_{\eps
\vec{r}}\|\leq\frac{1}{2}.$$
In other words, $\wt\psi:=\frac 1\eps \psi$ maps a ball $B^m(0,1)$ continuously into $\R^m$ in such a way that $\|\wt\psi(\vec r)-\vec r\|\leq \frac 12$. It is now an easy exercise to check that point $0\in\R^m$ lies in the image of $\wt\psi$. However, that means that $\psi(\vec {r_0})=0$ for some $\vec{r_0}$, and hence $\wt a_{\vec{r_0}}=\Phi(0,d)$ for  some $d\in\BM_0(I,\R^m)$. 
By property (\textbf{C}), $\left[\wt a_{\vec{r_0}}\right]_{t\in[0,1]}=0$. Finally, by Lemma \ref{lem:reparam},
$$0=\left[\wt a_{\vec r_0}(t)\right]_{t\in[0,1]}=\left[a_{\vec r_0}(t)\right]_{t\in[0,\eps]},$$
which finishes the proof.
\end{proof}

\section{The proof of PMP}\label{sec:proof}
In this section we will finish the proof of Theorem \ref{thm:pmp_A}. In our considerations it is crucial to understand the geometry of the cone $\B K^u_\tau$ of infinitesimal variations along the optimal trajectory $\B f(\B x(t),u(t))$. We interpreted $\B K^u_\tau$ as the set of all
directions in $\bm A_{\bm x(t_1)}$ in which one can move the
point $\bm x(t_1)$ by performing needle variations of the control
$u$ associated with symbols $\A=(\tau_i, v_i,\tau,\delta t_i, \delta t)_{i=1,\hdots,k}$, where $\tau$ is fixed. Consequently, a movement in the direction of the ray 
$$\bm\Lambda_{\bm
x(t_1)}:=\theta_{x(t_1)}\oplus\R_+\cdot(-\pa_t)\subset
E_{x(t_1)}\oplus\sT_{\ul{x}(t_1)}\R=\bm A_{\bm x(t_1)}$$ 
would correspond to a variation which
decreases the total cost of the trajectory without making changes
in the $E$-evolution. Such a behaviour should not be
possible if $\B f(\B x(t),u(t))$ is a solution of the OCP \eqref{eqn:P_A}, so we may expect that the ray $\B\Lambda_{\B x(t_1)}$ can be separated from the cone $\B K^u_\tau$ in such a case. This result is formulated in Theorem \ref{thm:separation_K_Lambda}. In the proof we use technical Lemma \ref{lem:htp_of_a_r} to deduce the existence of $E$-paths realising certain $E$-homotopy classes from the infinitesimal picture expressed in the language of the cone $\B K^u_\tau$ and the ray $\B\Lambda_{\B x(t_1)}$. When Theorem \ref{thm:separation_K_Lambda} is proved, to finish the proof of Theorem \ref{thm:pmp_A} we need only to follow a few rather technical steps from the original proof of Pontryagin and his collaborators \cite{pontryagin}.

\subsection{The geometry of the cone $\bm K^u_\tau$ }

Throughout this section we assume that the controlled pair $(\bm x(t),u(t))$ is a solution of the OCP \eqref{eqn:P_A}. All results obtained in this section are valid under this assumption. 

\begin{theorem}\label{thm:separation_K_Lambda} Let $(\B{x}(t), u(t))$, for $t\in[t_0,t_1]$, be a solution of the optimal control
problem \eqref{eqn:P_A}. Then the ray
$\bm{\Lambda}_{\B x(t_1)}$ and the convex cone
$\B K^u_\tau$ can be separated for any $\tau\in(t_0,t_1)$, which is a regular point of $u$.
\end{theorem}

The idea of the proof is the following.  Assuming the contrary we will
construct a family of symbols
$\A(\vec r)$, where $\vec r\in E_{x(t_1)}$, such that the associated infinitesimal variations
$\bm d^{\A(\vec r)}(s)$ are uniformly regular w.r.t. $\vec r$ at $s=0$ and  point into the 
directions $\vec r-\pa_t\in E_{x(t_1)}\oplus\T_{\ul x(t_1)}\R=\B A_{\B x(t_1)}$.
For such a family we will be able to use Lemma \ref{lem:htp_of_a_r} to deduce 
that for some $\vec r_0$ the variation $\B d^{\A(\vec r_0)}(s)$ has special properties. Next we will show that in such a case the pair $(\B x(t),u(t))$ cannot be a solution of the OCP \eqref{eqn:P_A}.

\begin{proof}
Assume the contrary, i.e., that the convex cone $\B K^u_\tau$ and the ray $\bm{\Lambda}_{\B{x}(t_1)}$ cannot be separated. Denote by  $\B \lambda:=\theta_{x(t_1)}-\pa_t\in E_{x(t_1)}\oplus\T_{\ul{x}(t_1)}\R$ a vector spanning $\B \Lambda_{\B x(t_1)}$. Then there exists a basis $\{\B e_1,\hdots,\B e_m\}$ of
$ E_{x(t_1)}\oplus\theta_{\ul{x}(t_1)}\subset
\B{A}_{\B x(t_1)}$ such that vectors
$\B\lambda,\B\lambda+\B e_i,\B\lambda-\B e_i$ lie in
$\B K^u_\tau$ for $i=1,\hdots, m$. 

Denote by $\A$, $\A_i$, and $\BB_i$ some symbols such that elements in $\bm K_\tau^u$ corresponding to these symbols are $\bm d^\A(0)=\B\lambda$, $\bm d^{\A_i}(0)=\B\lambda+\B e_i$, and $\bm d^{\BB_i}(0)=\B\lambda-\B e_i$, respectively. We deal with a finite set of symbols, hence we can assume that they all are of the form $(\tau_i, v_i,\tau, \delta t_i, \delta t)_{i=1,\hdots,k}$, where $\tau_i$, $v_i$, k, and $\tau$ are
fixed, and that they differ by $\del t_i$ and $\del t$ (we can
always add a triple $(\tau_i,v_i,\del t_i=0)$ to a symbol without
changing anything --- cf. Remark \ref{rem:var_dt=0}). For any
$\vec{r}=\sum_{i=1}^mr^i\B e _i\in B^m(0,1)=\{\vec r:\|\vec r\|=\sum_i|r_i|\leq 1\}$  we
may define a new symbol
$$\A(\vec{r})=\left(1-\sum_{i=1}^m|r^i|\right)\A+\sum_{i=1}^m
h^+(r^i)\A_i+\sum_{i=1}^m h^-(r^i)\BB_i,$$ where
$h^+(r)=\max\{r,0\}$ and $h^-(r)=\max\{-r,0\}$ are non-negative, and
the convex combination of symbols is defined using the natural
rule
\begin{align*}&\nu(\tau_i, v_i,\tau, \delta t_i, \delta
t)_{i=1,\hdots,k}+\mu(\tau_i, v_i,\tau, \delta t^{'}_i, \delta
t^{'})_{i=1,\hdots,k}\\&=(\tau_i, v_i,\tau, \nu\delta t_i+\mu\del
t_i^{'}, \nu\delta t+\mu\del t^{'})_{i=1,\hdots,k}\,.
\end{align*} 
We will now study the properties of $\B A$-paths $s\mapsto\B d^{\A(\vec r)}(s)$ corresponding to symbols $\A(\vec r)$ (see Remark \ref{rem:interpretation_K}).  

From \eqref{eqn:d_at_0} it is straightforward to verify that 
$$\B{d}^{\A(\vec{r})}(0)=\B\lambda+\vec{r}.$$

If follows from Theorem \ref{thm:1st_main} that, since the numbers $\del t(\vec r)$ and $\del t_i(\vec r)$ in the symbol $\A(\vec r)$ depend continuously on $\vec r$, which takes values in a compact set, we may choose $\theta>0$ such that $\bm d^{\A(\vec r)}(s)$ is well-defined for $s\in[0,\theta]$ and all $\vec r\in B^m(0,1)$,

Now consider $\ul d^{\A(\vec r)}(s)$ --- the projections of the family of $\B A$-paths $\B d^{\A(\vec r)}(s)$ to the algebroid $\T\R$.
Observe that since $\bm d^{\A(\vec r)}(s)$ are uniformly regular, so are  $\ul d^{\A(\vec r)}(s)$. Since in canonical coordinates on $\T\R$ we have $\ul d^{\A(\vec
r)}(0)=-1$, there exist a number $0<\eta\leq\theta$ such that 
\begin{equation}\label{eqn:cost_less_zero}
\int_0^\eps\ul d^{\A(\vec r)}(s)\dd s<0
\end{equation}  for all
$\eps\leq\eta$ and all $\vec r\in B^m(0,1)$. This property will be
used later.

After projecting $\B d^{\A(\vec{r})}(s)$ from $\B A $ onto $E$,
we obtain a family of bounded measurable admissible paths
$d^{\A(\vec{r})}(s)$, again uniformly
regular at $s=0$ w.r.t. $\vec{r}\in {B}^m(0,1)\subset \R^m$,
and such that
$d^{\A(\vec{r})}(0)=\vec{r}\in\R^m\approx E_{x(t_1)}$. In other words, the paths $d^{\A(\vec r)}(\cdot)$ satisfy the assumptions of Lemma \ref{lem:htp_of_a_r} and, consequently, there exists a vector $\vec r_0\in B^m(0,1)$ and a number $0<\eps< \eta$ such that 
\begin{equation}\label{eqn:htp_zero}
\left[d^{\A(\vec r_0)}(t)\right]_{r\in[0,\eps]}=\left[0\right].
\end{equation}

Properties \eqref{eqn:cost_less_zero} and \eqref{eqn:htp_zero} contradict the optimality of $\B f(\B x(t),u(t))$. Indeed, from \eqref{eqn:var_htp} we know that $\B f(\B x(t,s), u_s(t))$ --- the variation of the trajectory $\B f(\B x(t),u(t))$ associated with a needle variation $u_s(t):=u^{\A(\vec r_0)}_s(t)$ satisfies
$$\left[\B{f}(\B{x}(t,\eps),u_\eps(t))\right]_{t\in[t_0,t_1+\eps\del
t(\vec r_0)]} =\left[\B{f}(\B{x}(t),u(t))\right]_{t\in[t_0,t_1]}\left[\B{d}^{\A(\vec
r_0)}(s)\right]_{s\in[0,\eps]}.$$
Projecting the above equality to the algebroid $E$ and using \eqref{eqn:htp_zero} we get
$$\left[f(x(t,\eps),u_\eps(t))\right]_{t\in[t_0,t_1+\eps\del
t(\vec r_0)]}=\left[f(x(t),u(t))\right]_{t\in[t_0,t_1]},$$
and hence the $E$-homotopy classes agree.

What is more, the $\T\R$-projection gives
$$\left[L(x(t,\eps),u_\eps(t))\right]_{t\in[t_0,t_1+\eps\del
t(\vec r_0)]}=\left[L(x(t),u(t))\right]_{t\in[t_0,t_1]}\left[\ul d^{\A(\vec
r_0)}(s)\right]_{s\in[0,\eps]}.$$ From \eqref{eqn:cost_less_zero} we deduce
that the total costs satisfy the following inequality:
\begin{align*}\ul x(t_1+\eps\del t(\vec r_0),\eps)=&\int_{t_0}^{t_1+\eps\del
t(\vec r_0)}L(x(s,\eps),u_\eps(s))\dd s=\int_{t_0}^{t_1}L(x(s),u(s))\dd
s+\int_0^\eps \ul d^{\A(\vec r_0)}(s)\dd
s\\ \overset{\eqref{eqn:cost_less_zero}}<
&\int_{t_0}^{t_1}L(x(s),u(s))\dd s=\ul x(t_1).
\end{align*}
The above inequality proves that $\B f(\B x(t),u(t))$ cannot be a solution of the OCP \eqref{eqn:P_A}, which stays in a contradiction to our assumptions. 
\end{proof}

To finish the proof of Theorem \ref{thm:pmp_A} we will now follow the steps of the original result of \cite{pontryagin}. All
the important information is contained in Theorem
\ref{thm:1st_main} telling us that the set of infinitesimal
variations $\bm K_\tau^u$ is a convex cone with elements defined by means of a
local one-parameter group $\bm B_{tt_0}$ (see
\eqref{eqn:d_at_0}) and in Theorem \ref{thm:separation_K_Lambda} describing the geometry of this cone. The structure of an AL
algebroid, necessary to prove the above results, will now play no
essential role.

\subsection{The construction of $\bm{\xi}(t)$ and the ``maximum principle''}
Fix an element $\tau\in(t_0,t_1)$ to be a regular point of $u$. In
view of Theorem \ref{thm:separation_K_Lambda} there exists a non-zero covector
$\bm\xi(t_1)\in\bm A^\ast_{\bm x(t_1)}$ separating $\bm K^u_\tau$ and $\bm\Lambda_{\bm x(t_1)}$; that is,
\begin{equation}
\label{eqn:K_cov} 
\<\bm d,\bm\xi(t_1)>_{\bm{\tau}}\leq 0\leq\<\bm \lambda,\bm \xi(t_1)>_{\bm \tau} \text{
\ for every $\bm d\in\bm K^u_\tau$}.
\end{equation}
Let us define $\bm\xi(t):=\bm B^\ast_{t t_1}(\bm\xi(t_1))\in\bm
A^\ast_{\bm x(t)}$ for $t\in[t_0,t_1]$.

\begin{lemma}\label{lem:HM} For every $t\in[t_0,\tau]$ which is a regular point of the control $u(t)$ the following ``maximum principle'' holds:
$$\bm{H}(\bm x(t),\bm \xi(t),u(t))=\sup_{v\in U}\bm{H}(\bm x(t),\bm\xi(t),v).$$
Moreover, $\bm{H}(\bm x(\tau),\bm\xi(\tau),u(\tau))=0$.
\end{lemma}
\begin{proof}
Choose a regular point $t\in[t_0,\tau]$, take an arbitrary element
$v\in U$ and a number $\del t_1>0$, and consider a symbol
$\A=(\tau_1=t,\delta t_1,v_1=v,\tau,\delta t=0)$. The
corresponding element $\bm{d}^\A(0)\in \bm K^u_\tau$ equals
$\bm{B}_{t_1 t}[\bm f (\bm x (t),v)-\bm f (\bm x (t),u(t))]\delta
t_1$ (cf. \eqref{eqn:d_at_0}). From \eqref{eqn:K_cov} we obtain
\begin{align*}
&\:\,0\geq\<\bm{B}_{t_1
t}[\bm{f}(\bm{x}(t),v)-\bm{f}(\bm{x}(t),u(t))],\bm\xi(t_1)>_{\bm{\tau}}\delta
t_1\\
&\overset{\text{rem. } \ref{rem:B_paring}}{=}
\<\bm{B}_{t t_1}\bm{B}_{t_1
t}[\bm{f}(\bm{x}(t),v)-\bm{f}(\bm{x}(t),u(t))],\bm{B}^\ast_{t t_1}(\bm\xi(t_1))>_{\bm{\tau}}\delta t_1\\
&\phantom{X}=\<\bm f(\bm{x}(t),v)-\bm{f}(\bm{x}(t),u(t)),\bm\xi(t)>_{\bm\tau}\delta t_1\\
&\phantom{X}=\Big(\bm{H}(\bm x(t),\bm{\xi}(t),v)-\bm{H}(\bm
x(t),\bm{\xi}(t),u(t))\Big)\delta t_1.
\end{align*}
Since $\del t_1>0$, we have $\bm{H}(\bm
x(t),\bm{\xi}(t),v)\leq\bm{H}(\bm x(t),\bm{\xi}(t),u(t))$
for arbitrarily chosen $v\in U$.

To prove the second part of the assertion, consider a symbol
$\mathfrak{v}=(\tau,\del t)$. The associated element
$\bm{d}^{\mathfrak{v}}(0)$ is $\bm B_{t_1\tau}\left(\bm f (\bm x (\tau),u(\tau))\right)\del
t\in\bm K^u_\tau$. Consequently, from \eqref{eqn:K_cov}, we obtain
\begin{align*}
&\:\,0\geq\<\bm B_{t_1\tau}\left(\bm f (\bm{x}(\tau),u(\tau))\right),\bm\xi(t_1)>_{\bm{\tau}}\delta t\\
&\overset{\text{rem. } \ref{rem:B_paring}}{=}
\<\bm f (\bm{x}(\tau),u(\tau)),\bm\xi(\tau)>_{\bm{\tau}}\delta t=\bm{H}(\bm x(\tau),\bm{\xi}(\tau),u(\tau))\del t.
\end{align*}
Since $\del t$ can be arbitrary, we deduce that $\bm{H}(\bm
x(\tau),\bm{\xi}(\tau),u(\tau))=0$.
\end{proof}

\subsection{The condition $\bm{H}(\bm{x}(t),\bm{\xi}(t),u(t))=0$} 
To finish the proof just two more things are left.
We have to check that the Hamiltonian
$\bm{H}(\bm{x}(t),\bm\xi(t),u(t))$ is constantly 0 along the
optimal trajectory, and we have to extend the ''maximum principle'' to all regular  $t\in[t_0,t_1]$ (so far it holds only on the interval $[t_0,\tau]$, where
$\tau<t_1$ is a fixed regular point).

\begin{lemma}\label{lem:M=0}
For $\bm\xi(t)$ defined as above, the equality
$\bm{H}(\bm{x}(t),\bm\xi(t),u(t))=0$ holds at every regular point
$t\in[t_0,\tau]$ of the control $u$.
\end{lemma}

\begin{proof}
Denote by $P$ the closure of the set $\{u(t):t\in[t_0,\tau]\}$.
Since $u(t)$ is bounded, $P$ is a compact subset of $U$. Define a
new function $\m:\B{A}^*\lra\R$ by the formula
$$\m(\bm x,\bm{\xi}):=\max_{v\in P}\bm{H}(\bm x,\bm{\xi},v).$$
It follows from the previous lemma that
$\m(\bm{x}(t),\bm\xi(t))=\bm{H}(\bm x(t),\bm \xi(t),u(t))$ at
every regular point $t$ of $u$. We shall show that
$\m(\bm{x}(t),\bm\xi(t))$ is constant on $[t_0,\tau]$, and hence
equals $\bm{H}(\bm x(\tau),\bm\xi(\tau),u(\tau))=0$ (confront Lemma
\ref{lem:HM}). Observe that the function $\bm H\left(\bm
x(t),\bm\xi(t),v\right)$ is uniformly (for all $v\in P$) Lipschitz
w.r.t. $t$. Indeed, in local coordinates $(x^a,\ul
x,y^i,\ul y)$ on $\bm A$ and $(x^a,\ul x,\xi_i,\ul \xi)$ on $\bm
A^\ast$ we have $\bm H\left(\bm
x(t),\bm\xi(t),v\right)=\xi_i(t)f^i(x(t),v)+\ul{\xi}L(x(t),v)$.
Note that, by assumption, functions $f^i(x,v)$ and $L(x,v)$ are $C^1$ w.r.t. $x$ and their $x$-derivatives are continuous functions of both variables. Since $x(t)$ is an AC path with bounded derivative, functions
$\frac{\pa f^i}{\pa x^a}(x(t),v)$ and $\frac{\pa L}{\pa
x^a}(x(t),v)$, as well as functions $f^i(x(t),v)$ and $L(x(t),v)$,
are bounded in $[t_0,\tau]\times P$. As the evolution of $\xi(t)$
is governed by \eqref{eqn:par_trans_*A} and $u(t)$ is
bounded on $[t_0,\tau]$, the derivatives $\pa_t\xi^k(t)$ are also
bounded on $[t_0,\tau]$. Consequently, since the path $(\bm
x(t),\bm\xi(t))\in\bm A^\ast$ can be covered by a finite number of
coordinate charts, the $t$-derivative of $\bm H\left(\bm
x(t),\bm\xi(t),v\right)$ is bounded on $[t_0,\tau]\times P$. As a
result there exists a number $C$ such that
$$\left|\bm{H}(\bm{x}(t),\bm\xi(t),v)-\bm{H}(\bm{x}(t^{'}),\bm\xi(t^{'}),v)\right|\leq C|t-t^{'}|$$
for all $t,t^{'}\in [t_0^{'},\tau]$ and for any $v\in P$. Observe
also that
\begin{align*}
&\frac\pa{\pa t}\bm H\left(\bm
x(t),\bm\xi(t),v\right)|_{v=u(t)}=\pa_t\xi_i(t)f^i(x(t),u(t))+\\
&\phantom{==}+
\xi_i(t)\pa_t f^i(x(t),v)|_{v=u(t)}+\ul{\xi}\pa_tL(x(t),v)|_{v=u(t)}\overset{\eqref{eqn:par_trans_*A}}{=}\\
&=
\left[-\rho^a_i\left(x\right)\left(\frac{\pa f^k}{\pa
x^a}\left(x,u(t)\right)\xi_k(t)+\frac{\pa L}{\pa
x^a}\left(x,u(t)\right)\ul\xi(t)\right)+c^k_{ji}\left(x\right)f^j\left(x,u(t)\right)\xi_k(t) \right]f^i(x,u(t))\\
&\phantom{==} +\xi_i(t)\frac{\pa f^i}{\pa x^a}(x(t),u(t))\rho^a_k(x)f^k(x(t),u(t))+\ul{\xi}\frac{\pa L}{\pa x^a}(x(t),u(t))\rho^a_k(x)f^k(x(t),u(t))\\
&\ =\xi_k(t)c^k_{ij}(x)f^i(x(t),u(t))f^j(x(t),u(t))=0
\end{align*}
by the skew-symmetry of $c^k_{ij}(x)$.

Now take any regular points $t,t^{'}\in[t_0,\tau]$. Since
$u(t),u(t^{'})\in P$, we have
\begin{align*}
-C|t-t^{'}|&\leq\bm{H}(\bm{x}(t),\bm\xi(t),u(t^{'}))-\bm{H}(\bm{x}(t^{'}),\bm\xi(t^{'}),u(t^{'}))\\&\leq\m(\bm{x}(t),\bm\xi(t))-
\m(\bm{x}(t^{'}),\bm\xi(t^{'}))\\
&=\bm{H}(\bm{x}(t),\bm\xi(t),u(t))-\bm{H}(\bm{x}(t^{'}),\bm\xi(t^{'}),u(t^{'}))\\
&\leq\bm{H}(\bm{x}(t),\bm\xi(t),u(t))-\bm{H}(\bm{x}(t^{'}),\bm\xi(t^{'}),u(t))\leq C|t-t^{'}|;
\end{align*}
i.e., $\m(\bm x(t),\xi(t))$ satisfies the Lipschitz condition on
the set of regular points (dense in $[t_0,\tau]$). It is also a
continuous map (since $x(t)$ and $\xi(t)$ are AC and the maps $L(x,v)$,
$f(x,v)$ are continuous in both variables), therefore it is
Lipschitz on the whole interval $[t_0,\tau]$. By Rademacher{'}s
theorem, $\m(\bm x(t),\xi(t))$ is almost everywhere differentiable on $[t_0,\tau]$. Now take any point $t$ of
differentiability of $\m(\bm{x}(t),\bm\xi(t))$ which is also a
point of the regularity of the control $u$. We have
$$\m(\bm{x}(t^{'}),\bm\xi(t^{'}))-\m(\bm{x}(t),\bm\xi(t))\geq
\bm{H}(\bm
x(t{'}),\bm{\xi}(t^{'}),u(t))-\B{H}(\B{x}(t),\bm\xi(t),u(t)).$$
For $t^{'}>t$, we get
$$\frac{\m(\B{x}(t^{'}),\bm\xi(t^{'}))-\m(\B{x}(t),\bm\xi(t))}{t^{'}-t}\geq\frac{\B{H}(\B{x}(t^{'}),\bm\xi(t^{'}),u(t))-
\B{H}(\B{x}(t),\bm\xi(t),u(t))} {t^{'}-t}.$$ Consequently,
$$\frac{d}{dt}\m(\B{x}(t),\bm\xi(t))\geq\frac{\partial}{\partial
t^{'}}\Big|_{t^{'}=t}\B{H}(\B{x}(t^{'}),\bm\xi(t^{'}),u(t))=0.$$
Similarly, for $t^{'}<t$, we would get
$\frac{d}{dt}\m(\B{x}(t),\bm\xi(t))\leq 0$. We deduce that
$\frac{d}{dt}\m(\B{x}(t),\bm\xi(t))=0$ a.e. in
$[t_0,\tau]$; hence $\m(\B{x}(t),\bm\xi(t))$ is constant and equals $\bm H(\bm x(\tau),\bm \xi(\tau),u(\tau))=0$.
\end{proof}

\subsection{Extending the ''maximum principle'' to $[t_0,t_1]$}
\begin{lemma}\label{lem:K_increase} Let $t\in[t_0,\tau]$ be any regular point of
the control $u$. Then 
$$\B{K}^u_t\subset
\cl\left(\B{K}^u_\tau \right).$$
\end{lemma}

\begin{proof}
Consider an element $\bm d\in\B{K}_t^u$ of the
form
$$\B{d}=\bm B_{t_1t}[\B{f}(\B{x}(t),u(t))]\del
t+\sum_{i=1}^s\B{B}_{t_1\tau_i}\Big[\B{f}(\B{x}(\tau_i),v_i)-\B{f}(\B{x}(\tau_i),u(\tau_i))\Big]\del
t_i.$$ 
Since $\cl\left(\bm K_\tau^u\right)$ is a convex cone it is enough to show that $\bm B_{t_1t}\left(\bm{f}(\B{x}(t),u(t))\del
t\right)$ and \\$\sum_{i=1}^s\B{B}_{t_1\tau_i}\left[\B{f}(\B{x}(\tau_i),v_i)-\B{f}(\B{x}(\tau_i),u(\tau_i))\right]\del
t_i$ belong to $\cl\left(\B{K}^u_\tau\right)$. The later clearly belongs to $\bm K_\tau^u\subset\cl (\bm K_\tau^u)$ since $\tau_i\leq t<\tau$.

Assume now that  $\B{B}_{t_1
t}\left[\B{f}(\B{x}(t),u(t))\del t\right]$ does not belong to $\cl\left(\bm{K}^u_\tau\right)$. Since this set is a closed convex cone, there exists a covector $\ol{\bm{\xi}}(t_1)\in \bm A^\ast_{\bm x(t_1)}$ strictly separating  $\cl\left(\bm{K}^u_\tau\right)$ from $\left\{\B{B}_{t_1 t}\left(\B{f}(\B{x}(t),u(t))\delta
t\right)\right\}$; i.e., 
$$\<\B{d},\ol{\B{\xi}}(t_1)>_{\bm{\tau}}\leq 0<\<\B{B}_{t_1 t}\left(\B{f}(\B{x}(t),u(t))\del
t\right),\ol{\B{\xi}}(t_1)>_{\bm{\tau}}\quad \text{ for any $\B{d}\in
\cl(\bm{K}^u_\tau$).}$$ 
Define $\ol{\B{\xi}}(t):=\B{B}^\ast_{t t_1}(\ol{\B{\xi}}(t_1))$. Lemmas
\ref{lem:HM} and \ref{lem:M=0} hold for $\ol{\B{\xi}}(t)$ (we needed only $\<\bm d,\ol{\bm\xi}(t_1)>\leq 0$ for $\bm d\in \bm K_\tau^u\subset\cl(\bm K_\tau^u)$ in the proofs), hence,
in particular,
$$\<\B{f}(\B{x}(t),u(t)),\ol{\B{\xi}}(t)>_{\bm{\tau}}=\B{H}(\bm
x(t),\ol{\B{\xi}}(t),u(t))=0,$$ as $t$ is a regular point of $u$.
On the other hand,
$$0<\<\B{B}_{t_1 t}\big(\B{f}(\B{x}(t),u(t))\del
t\big),\ol{\B{\xi}}(t_1)>_{\bm{\tau}}=\<\B{f}(\B{x}(t),u(t)),\ol{\B{\xi}}(t)>_{\bm{\tau}}\del
t=\bm H(\bm x(t),\ol{\bm\xi}(t),u(t))\del t,$$
and hence $\B{H}(\bm
x(t),\ol{\B{\xi}}(t),u(t))\neq 0$. This contradiction finishes the proof.
\end{proof}

Note that, so far, we could define the set $\B{K}^u_\tau$ only for
a regular point $\tau<t_1$. With help of the above lemma we can
also define $\B{K}^u_{t_1}$ as the direct limit of the increasing
family of sets $\cl\left(\B{K}^u_\tau\right)$,
$$\B{K}^u_{t_1}:=\bigcup_{\tau<t_1, \text{ $\tau$ regular}}\cl\left(\B{K}^u_\tau\right).$$
It is clear that
$\B{K}^u_{t_1}$  is a convex cone in
$\B{A}_{\B{x}(t_1)}$.
It has geometric properties analogous to the
properties of $\B{K}^u_\tau$ described in Theorem \ref{thm:separation_K_Lambda}.

\begin{lemma}\label{lem:K_at_t1}
The ray $\bm{\Lambda}_{\B{x}(t_1)}$ and the convex cone
$\B{K}^u_{t_1}$ are separable.
\end{lemma}
\begin{proof}
Assume the contrary. Then there exists a vector $k\in \bm K^u_{t_1}\cap\bm\Lambda_{\bm x(t_1)}$ and a basis $\{e_1,\hdots,e_m\}$ of $ E_{x(t_1)}\oplus\theta_{\ul x(t_1)}$ such that
\begin{align}
k\pm e_i\in \bm K_{t_1}^u\quad \text{for $i=1,\hdots,m$}\label{cond:sep_B}.
\end{align} 
Since $\bm K^u_{t_1}$ is a limit of an increasing family of sets, there exists a regular $t<t_1$ such that \eqref{cond:sep_B} holds for $\cl\left(\bm K_t^u\right)$. In other words, $\bm\Lambda_{\bm x(t_1)}$ and $\cl\left(\bm K_t^u\right)$ are not separable. But then also $\bm\Lambda_{\bm x(t_1)}$ and $\bm K_t^u$ are not separable (separability of convex sets is a closed condition). This contradicts
Theorem \ref{thm:separation_K_Lambda}.
\end{proof}

Now choose a non-zero covector $\B{\xi}(t_1)\in\B{A}^\ast_{\B{x}(t_1)}$
separating $\bm K_{t_1}^u$ and $\bm\Lambda_{\bm x(t_1)}$ and define $\bm\xi(t):=\bm B^\ast_{tt_1}(\bm \xi(t_1))$.  We have 
$$\<\B{d},\B{\xi}(t_1)>_{\bm{\tau}}\leq 0\leq\<\bm\lambda,\B{\xi}(t_1)>_{\bm{\tau}}\quad\text{for every $\B{d}\in \B{K}^u_{t_1}$}.$$ 
Since, by construction, $\bm K^u_\tau\subset \bm K_{t_1}^u$, the covector $\bm\xi(t_1)$ separates also 
$\bm K^u_\tau$ and $\bm\Lambda_{\bm x(t_1)}$. This is enough for Lemmas \ref{lem:HM} and \ref{lem:M=0} to hold for $\bm \xi(t)$.
As a consequence, for every regular point of the control $u$, we have
$$
\bm{H}(\bm x(t),\bm \xi(t),u(t))=\sup_{v\in U}\bm{H}(\bm
x(t),\bm\xi(t),v)=0\,.
$$
Finally, since $\<\bm \lambda, \bm \xi(t_1)>\geq 0$, we have 
 $\ul\xi(t)\equiv\ul\xi(t_1)\leq 0$. This finishes the proof
of Theorem \ref{thm:pmp_A}. \hfill$\qed$

\bibliographystyle{amsplain}
\bibliography{books,articles}


\end{document}